\documentclass[a4paper]{amsart}
\date{December 9, 2004}
\usepackage[dvips]{graphicx}
\usepackage{verbatim,enumerate}
\usepackage{amssymb}
\usepackage[usenames]{color}
\usepackage{amsthm}
\theoremstyle{plain}
 \newtheorem{theorem}{Theorem}[section]
 \newtheorem*{theorem*}{Theorem}
 \newtheorem*{lemma*}{Lemma}
 \newtheorem{proposition}[theorem]{Proposition}
 \newtheorem{fact}[theorem]{Fact}
 \newtheorem{fact*}{Fact}

 \newtheorem{lemma}[theorem]{Lemma}
 
\theoremstyle{remark}
 \newtheorem{definition}[theorem]{Definition}
 \newtheorem{remark}[theorem]{Remark}
 \newtheorem*{remark*}{Remark}
 \newtheorem*{acknowledgements}{Acknowledgements}
 \newtheorem{example}[theorem]{Example}

\numberwithin{equation}{section}

\newcommand{\vect}[1]{\boldsymbol{#1}}
\newcommand{\Z}{\boldsymbol{Z}}
\newcommand{\R}{\boldsymbol{R}}
\newcommand{\C}{\boldsymbol{C}}

\newcommand{\Lor}{\boldsymbol{L}}

\newcommand{\trace}{\operatorname{trace}}
\newcommand{\mb}[1]{{\mathbf #1}}
\renewcommand{\phi}{\varphi}
\newcommand{\ord}{\operatorname{ord}}

\newcommand{\SL}{\operatorname{SL}}
\newcommand{\SU}{\operatorname{SU}}
\newcommand{\U}{\operatorname{U}}

\newcommand{\PSL}{\operatorname{PSL}}

\newcommand{\Herm}{\operatorname{Herm}}
\newcommand{\A}{\mathcal{A}}
\renewcommand{\O}{\mathcal{O}}
\newcommand{\fO}{\mathfrak{O}}

\renewcommand{\L}{\mathcal{L}}

\renewcommand{\Re}{\operatorname{Re}}
\renewcommand{\Im}{\operatorname{Im}}

\newcommand{\inner}[2]{\left\langle{#1},{#2}\right\rangle}
\title[Singularities of flat fronts]{
Singularities of flat
 fronts in Hyperbolic 3-space
}
\author{M.~Kokubu, W.~Rossman, K.~Saji, M.~Umehara and K.~Yamada}

\address[Masatoshi Kokubu]{%
   Department of Natural Science,
   School of Engineering,
   Tokyo Denki University,
   2-2 Kanda-Nishiki-Cho,
   Chiyoda-Ku, Tokyo 101-8457,
   Japan
}
\email{kokubu@cck.dendai.ac.jp}
\address[Wayne Rossman]{%
   Department of Mathematics, Faculty of Science,
   Kobe University,
   Rokko, Kobe 657-8501, Japan%
}
\email{wayne@math.kobe-u.ac.jp}
\address[Kentaro Saji]{%
   Department of Mathematics,
   Graduate School of Science,
   Hiroshima University,
   Higashi-Hiroshima 739-8526,
   Japan
}
\email{sajik@hiroshima-u.ac.jp}
\address[Masaaki Umehara]{%
   Department of Mathematics, Graduate School of Science,
   Osaka University,
   Toyonaka, Osaka 560-0043,
   Japan
}
\email{umehara@math.wani.osaka-u.ac.jp}
\address[Kotaro Yamada]{%
   Faculty of Mathematics,
   Kyushu University 36,
   Higashi-ku, Fukuoka 812-8581, Japan%
}
\email{kotaro@math.kyushu-u.ac.jp}
\thanks{
  The third author has been supported by JSPS
  Research Fellowships for Young Scientists.
}
\begin{document}
\begin{abstract}
It is well-known that the unit cotangent bundle of any Riemannian
manifold has a canonical contact structure.
A surface in a Riemannian $3$-manifold is  called a (wave) front if it
is the projection of a Legendrian immersion into the unit cotangent
bundle.
We shall give easily-computable criteria for a singular point
on a front to be a cuspidal edge or a swallowtail.
Using this, we shall prove that generically flat fronts in the hyperbolic
$3$-space admit only cuspidal edges and swallowtails.
Moreover, we will show that 
every complete flat front (which is not rotationally symmetric)
has associated parallel surfaces whose singularities consist
of only cuspidal edges and swallowtails.
\end{abstract}
\maketitle
\section{Introduction}
\label{sec:intro}
It is well-known that the unit cotangent bundle of any Riemannian
$n$-manifold $N^n$ has a canonical contact structure.
Let $M^2$ be a $2$-manifold and $f\colon M^2\to N^3$ a $C^\infty$-map.
Then $f$ is  called a {\it $($wave$)$ front} if it
is the projection of a Legendrian immersion into 
the unit cotangent
bundle of $N^3$.
Now let  $f\colon M^2\to \widetilde M^3(c)$  be a
front, where $\widetilde M^3(c)$ is the space form of constant curvature
$c$. 
Then the associated parallel front $f_t\colon{}M^2\to\widetilde M^3(c)$, 
i.e.\ the surface that is equi-distant from $f$ at a distance $t$ 
(called the {\it parallel surfaces\/} of $f(=f_0)$), is well-defined.
Moreover, if $f$ is a flat immersion, so is $f_t$
for $t$ close to zero.
Using this fact, we shall define a flat front:
A front $f$ is {\it flat} in a neighborhood of $p\in M^2$ 
if either
\begin{enumerate}
\item  $p$ is a regular point of $f$ and the Gaussian curvature
       of $f$ near $p$ vanishes, or 
\item $p$ is a singular point of $f$ and $f_t$ is a flat immersion
       around $p$ for all $t\ne 0$ close to zero.
\end{enumerate}
A front $f:M^2\to \widetilde M^3(c)$ is called a {\it flat front\/} if it
 is flat everywhere on  $M^2$.

For the case $c=0$, there have 
appeared several articles concerning the singularities
of developable surfaces in $\R^3$.
In particular, Izumiya and Takeuchi \cite{IT}
proved that the set of developable surfaces whose singularities
are only cuspidal edges,
swallowtails or cuspidal cross caps are open and dense in the set of
non-cylindrical developable surfaces,
where
$(u,v)\mapsto (u,v^2,v^3)$ represents a cuspidal edge,
$(u,v)\mapsto (3u^4+u^2v,4u^3+2uv,v)$ a swallowtail,
and $(u,v)\mapsto (u, uv^3, v^2)$ a cuspidal cross cap.
Recently, geometric inequalities for complete flat fronts in hyperbolic
$3$-space and complete maximal surfaces with certain singularities in
Minkowski $3$-space were found in \cite{KUY2} and \cite{UY}.  
We also note that Kitagawa has made a deep investigation of 
flat tori in the $3$-sphere (\cite{Ki}, \cite{Kii}, \cite{Kiii}).
The study of global properties of surfaces with singularities 
is a newly-developing 
research area in differential geometry.

In this paper, we shall investigate singularities of flat surfaces in
the hyperbolic $3$-space $H^3=\widetilde M^3(-1)$.
The geometry of flat fronts in $H^3$ has been investigated in
\cite{KUY1}, \cite{KUY2}.
 In particular, an analogue of the Osserman inequality for minimal 
surfaces in $\R^3$ was given in \cite{KUY2}.
Like the case of constant mean curvature one surfaces in $H^3$,
flat surfaces have a representation formula in terms of holomorphic data,
found by J. A. G\'alvez, A. Mart\'\i{}nez and F. Mil\'an
\cite{GMM}:
Let $\omega$ and $\theta$ be holomorphic $1$-forms on a simply-connected
Riemann surface $M^2$ such that $|\omega|^2+|\theta|^2$ is positive
definite.
Then there exists a holomorphic immersion 
$E_f\colon{}M^2 \longrightarrow \SL(2,\C)$ such that
\begin{equation}\label{eq:canonical}
  E_f^{-1}dE_f^{}
  =\begin{pmatrix}
	    0 & \theta \\
	    \omega & 0
   \end{pmatrix},
\end{equation}
and its projection to $H^3$ gives a flat front $f=E_f^{}E_f^*$ in $H^3$,
where we regard $H^3$ as
\begin{equation}\label{eq:hyperbolic}
 H^3 = \SL(2,\C)/\SU(2)=\{aa^*\,;\,a\in\SL(2,\C)\}
  \qquad (a^*=\,^{t}\bar a).
\end{equation}
Moreover, any simply-connected flat front has such a representation
with respect to the complex structure induced by the
second fundamental form (see \cite{GMM}, \cite{KUY1} and \cite{KUY2}).
We call $\omega$ and $\theta$ in \eqref{eq:canonical} 
the {\em canonical forms\/} of $f$.
In Section~\ref{sec:local} of this paper, we will show the
following:
\begin{theorem}\label{thm:criterion}
 Let $f\colon{}M^2\to H^3$ be a flat front with canonical forms
 $(\omega=\hat \omega\, dz,\theta=\hat\theta\,dz)$,
 where $z$ is a local complex coordinate.
 \begin{enumerate}
  \item A point $p\in M^2$ is a singular point if and only if
	$|\hat \omega(p)|=|\hat \theta(p)|$ holds.
  \item The image of $f$ around a  singular point $p$ is locally
	diffeomorphic to a cuspidal edge if and only if
	\[
 	  \hat\omega'\hat \theta-\hat \theta'\hat \omega\ne 0
	  \qquad \text{and}\qquad
	  \Im\left(
	    \frac{(\hat \theta'/\hat \theta)-
	          (\hat \omega'/\hat \omega)}{
	           \sqrt{\hat \omega\hat \theta }}
	  \right) \ne 0
	\]
	hold at $p$, where ${}'=d/dz$.
  \item The image of $f$ around a singular point $p$ is locally
	diffeomorphic to a swallowtail if and only if
	\begin{align*}
	 \hat\omega'\hat\theta-\hat\theta\hat\omega'\neq 0,\qquad
	 \Im&\left(
	   \frac{(\hat\theta'/\hat\theta)-
	         (\hat\omega'/\hat\omega)}{
	          \sqrt{\hat \omega\hat\theta}}
	    \right) = 0 \\
	 &\hphantom{and}\text{and}\qquad
	 \Re\left(\frac{s(\hat\theta)-s(\hat\omega)}{
	         \hat\omega\hat\theta}\right)
	 \ne 0
	\end{align*}
	hold at $p$,
	where $s(\hat\omega)$ is the Schwarzian derivative
	$\{h,z\}$ of the function $h(z):=\int_{z_0}^z \omega$ with
	respect to $z$, 
        that is, 
	\begin{equation}\label{eq:schwarz}
	  s(\hat\omega)=
	  \{h,z\}=
	  \left(\frac{h''}{h'}\right)'-
	  \frac{1}{2}\left(\frac{h''}{h'}\right)^2=
	  \left(\frac{\hat\omega'}{\hat\omega}\right)'-\frac12
	  \left(\frac{\hat\omega'}{\hat\omega}\right)^2.
	\end{equation}
 \end{enumerate}
\end{theorem}
Consequently, cuspidal edges and swallowtails are stable under 
perturbations of $(\omega,\theta)$.
It is well-known that generic fronts (which might not be flat) admit only
cuspidal edges or swallowtails (see Arnol'd \cite[Section 21.6]{A}).
However, density of the set of such fronts 
within the smaller set of  flat fronts 
does not immediately follow.
Using Theorem~\ref{thm:criterion}, 
we shall prove a similar assertion for flat fronts in $H^3$
(Theorem~\ref{thm:genericity}).
Moreover, we shall prove 
the following global result in Section~\ref{sec:global}. 
A front $f\colon{}M^2\to H^3$ is called complete if there exist a
compact set $C\subset M^2$ and a symmetric $2$-tensor $T$ on $M^2$ such
that $T$ is identically $0$ outside $C$ and $ds^2+T$ is a complete
Riemannian metric of $M^2$, where $ds^2$ is the first fundamental form 
of $f$.
\begin{theorem}\label{thm:global}
 Let $f \colon M^2\to H^3$ be a complete flat front which is not a covering of
 an hourglass {\rm (}hourglasses are rotationally symmetric, 
 see Example~{\rm \ref{ex:revolution})}, 
 and let 
 $\{f_t\}$ be the family of parallel fronts of $f$.
 Then, except for only finitely many values of $t$, all the 
 singular points of $f_t$ are locally diffeomorphic to
 cuspidal edges or swallowtails.
\end{theorem}
 The image of the singular points under an hourglass 
 is a single point in $H^3$,
 namely the so-called ``cone-like singularity'' 
 (see Figure~\ref{fig:revolution} in Section~\ref{sec:example}),
 and any parallel front of an
 hourglass has the same singularity.
 Thus the assumption of Theorem~\ref{thm:global} is necessary.

Fronts which admit only cuspidal edges and swallowtails are called
A-mersions, and their topological properties have 
been investigated by Langevin, Levitt and Rosenberg \cite{LLR}.
The above theorem implies that complete flat fronts in $H^3$ 
are generically included in this category.

The union of singular sets for the entire parallel family of a given
flat front is called a {\em caustic}.  
We note that Roitman \cite{R} very recently studied the geometric
properties of flat surfaces, motivated by a classical
result of L. Bianchi (see Section~\ref{sec:caustic}).

To prove Theorem~\ref{thm:criterion}, we shall give criteria for a
singular point on a front to  be a cuspidal edge or a swallowtail, as
follows:
Let $N^3$ be a Riemannian $3$-manifold,
$U$ a domain on $(\R^2;u,v)$, and let
\[
  f=f(u,v) \colon U\longrightarrow N^3
\]
be a $C^\infty$-map with a singular point $p\in U$.
Then there exist three functions $a,b,c\in C^\infty(U)$ such that
\[
  df\left(\frac{\partial}{\partial u}\right)\wedge
  df\left(\frac{\partial}{\partial v}\right)
  =
      a(u,v)\, \frac{\partial}{\partial x}\wedge \frac{\partial}{\partial y}
 +    b(u,v)\, \frac{\partial}{\partial y}\wedge \frac{\partial}{\partial z}
 +    c(u,v)\, \frac{\partial}{\partial z} \wedge \frac{\partial}{\partial x},
\]
where $(x,y,z)$ is a local coordinate system of $N^3$.
The rank of 
a map defined by
\[
\mathcal G \colon  (u,v)\longmapsto \bigl(a(u,v),b(u,v),c(u,v)\bigr) \in \R^3
\]
does not depend on the choice of local coordinate
$(x,y,z)$ nor on the choice of coordinate $(u,v)$.
Now we assume $f$ is a front.
A singular point $p\in U$ of $f$ is called {\it non-degenerate} 
if the Jacobian matrix of $\mathcal G$ is of rank one at $p$. 
There exists a regular curve near a non-degenerate 
singular point $p$
\[
  \gamma=\gamma(t) \colon (-\varepsilon,\varepsilon)
   \longrightarrow U
\]
(called  a singular curve) such that $\gamma(0)=p$, and so that the image of
$\gamma$ coincides with the set of singularities of $f$ near $p$.
The tangential direction of $\gamma(t)$ is called the
{\it singular direction}, and a non-zero vector  
$\eta \in T_{\gamma(t)}U$ 
such that $df(\eta)=0$ represents the {\it null direction}. 
For each point $\gamma(t)$, vectors in the null direction 
$\eta(t)$ are 
uniquely determined up to non-zero scalar multiplication.
\begin{proposition}\label{prop:criterion-front}
 Let  $N^3$ be a Riemannian $3$-manifold and $p=\gamma(0) \in U$ be a
 non-degenerate singular point of a front $f \colon U\to N^3$.
 \begin{enumerate}
  \item The germ of the front $f$ at $p$ is locally diffeomorphic
	to a cuspidal edge if and only if $\eta(0)$ is not proportional
	to $\dot\gamma(0)$, where $\dot\gamma(t)=d\gamma(t)/dt$.
  \item The germ of the front $f$ at $p$ is locally diffeomorphic
	to a swallowtail if and only if $\eta(0)$ is proportional to
	$\dot\gamma(0)$ and
	\[
	  \left.\frac{d}{dt}\right|_{t=0}\!\!
	  \det{\bigl(\dot\gamma(t),\eta(t)\bigr)}
	  \neq 0.
        \]
 \end{enumerate}
\end{proposition}
We shall prove this proposition in Section~\ref{sec:criterion}.
These criteria are useful in other situations.
In fact, this proposition is applicable for the study of singularities 
of maximal surfaces in Minkowski space (see \cite{UY}).
\begin{acknowledgements}
The authors thank Shyuichi Izumiya and Go-o Ishikawa for
fruitful discussions and valuable comments, especially
Ishikawa for pointing out an error in a preliminary version
of our paper.
\end{acknowledgements}
\section{Criteria for singular points}
\label{sec:criterion}
\subsection{Preliminaries}
\label{sec:prelim}
First, we recall well-known properties for singular points from
\cite{BG}.
Let
\[
    \varphi=\varphi(w)\colon{} I \longrightarrow \R
\]
be a $C^{\infty}$-function defined on an open interval $I$ containing
the origin such  that $\varphi(0)=0$.
Then $\varphi$ has an {\em $\A_k$-singularity\/} at $0$ if
\[
    \varphi'(0)=\varphi''(0)=\dots=\varphi^{(k)}(0)=0,
    \qquad\text{and}\qquad
    \varphi^{(k+1)}(0)\neq 0,
\]
where $\varphi'=d\varphi/dw$ and $\varphi^{(j)}=d^j\varphi/dw^j$.
Here, we shall consider the cases $k=2$ and $3$.

Let $\Omega$ be an open subset of $(\R^3;x,y,z)$ containing the origin
$\vect{0}$.
A map
\[
    \Phi\colon{}I\times\Omega\longrightarrow \R
\]
is called an {\em unfolding\/} of $\varphi$ if
\[
    \varphi(w)=\Phi(w,\vect{0})
\]
holds.
Moreover, if $\varphi$ has an $\A_k$-singularity at $0$ 
and the matrix
\[
    \begin{pmatrix}
     \Phi_x(0,\vect{0}) &
     \Phi'_x(0,\vect{0}) &\dots &
     \Phi_x^{(k-1)}(0,\vect{0}) \\
     \Phi_y(0,\vect{0}) &
     \Phi'_y(0,\vect{0}) &\dots &
     \Phi_y^{(k-1)}(0,\vect{0}) \\
     \Phi_z(0,\vect{0}) &
     \Phi'_z(0,\vect{0}) &\dots &
     \Phi_z^{(k-1)}(0,\vect{0}) \\
    \end{pmatrix}
\]
is of rank $k$, then $\Phi$ is called a {\em versal unfolding\/} 
 of $\varphi$, 
where, for example,
\[
     \Phi'_x=\frac{\partial^2\Phi}{\partial x\partial w},\qquad
     \Phi_x^{(j)}=\frac{\partial^{j+1}\Phi}{\partial x\partial^j w}.
\]
The set
\[
     D_{\Phi}:=
     \left\{
        \vect{x}\in\Omega\,
	;\,
	\text{there exists a $w \in I$ with
	$\Phi(w,\vect{x})=\Phi'(w,\vect{x})=0$}
     \right\}
\]
is called the {\em discriminant set\/} of $\Phi$.
The following fact is useful: 
\begin{fact}[{\cite[Section 6]{BG}}]\label{fact:descriminant}
 Suppose $\varphi\colon{}I\to\R$ has an $\A_k$-singularity 
 {\rm(}$k=2$ or $3${\rm )} at $0$
and $\Phi\colon{}I\times\Omega\to\R$ a versal
 unfolding of $\varphi$.
 Then
 \begin{enumerate}
  \item $D_{\Phi}$ is locally diffeomorphic to a cuspidal edge at
	$\vect{0}$ if $k=2$.
  \item $D_{\Phi}$ is locally diffeomorphic to a swallowtail at
	$\vect{0}$ if $k=3$.
 \end{enumerate}
\end{fact}

\subsection{Non-degenerate singular points}
\label{subsec:front}
Let $N^3$ be a Riemannian $3$-manifold and $T^*_1N^3$ the unit cotangent
bundle.
A $C^{\infty}$-map $f \colon{}M^2  \to N^3$ is called a ({\em wave})
{\em front\/} if there exists a Legendrian immersion
$L_f\colon{}M^2 \to T^*_1 N^3$ such that $f=\pi\circ L_f$,
where $\pi\colon{}T^*_1 N^3\to N^3$ is the projection.
We call $L_f$ the {\em Legendrian lift\/} of $f$.
We shall use the following lemma, 
first pointed out by Zakalyukin \cite{Z}: 

\begin{lemma}[Zakalyukin \cite{Z}]\label{fact:Z}
 Let $U(\subset\R^2)$  be a neighborhood of
 the origin, and let $f_j \colon U\to\R^3 \ (j=1,2)$ be fronts.
 Suppose that $(0,0)$ is a singular point of $f_j$ and
 the set of regular points of $f_j$ is dense 
 in $U$ for each $j=1,2$.
 Then the following two statements are equivalent{\rm :}
 \begin{enumerate}
  \item\label{item:fact-Z-1} 
        There exist neighborhoods $V_1,V_2(\subset \R^2)$ of the origin
	$(0,0)$ and a local diffeomorphism on $\R^3$ which maps
	the image $f_1(V_1)$ to $f_2(V_2)$, namely the image of $f_1$ is
	locally diffeomorphic to that of $f_2$.
  \item\label{item:fact-Z-2} 
        There exists a local diffeomorphism $h$ on $\R^2$ and
	a local contact diffeomorphism $\Phi$ on $T^*_1 \R^3$
	which sends fibers to fibers such that
	$\Phi\circ L_{f_1}=L_{f_2} \circ h$, namely the lift
	$L_{f_1}$ is {\it Legendrian equivalent\/} to the lift
	$L_{f_2}$.
 \end{enumerate}
\end{lemma}

We shall prove the lemma in the appendix.
Again, we shall return to the general setting:
Since any contact structure is locally equivalent to the canonical
contact structure on $T_1^*\R^3$, we may restrict our attention to
fronts in the Euclidean $3$-space $\R^3$.
Let $(U;u,v)$ be a domain in $\R^2$ 
and $f\colon{}U\to\R^3$ a front.
Identifying the unit cotangent bundle $T^*_1\R^3$ 
with the unit tangent bundle $T_1\R^3\simeq \R^3\times S^2$, 
there exists a unit vector field
\[
     \nu\colon{}U \longrightarrow S^2\subset \R^3
\]
such that the Legendrian lift $L_f$ is expressed as $(f,\nu)$.
Since $L_f=(f,\nu)$ is Legendrian,
\[
    \inner{df}{\nu} = 0\qquad
    \text{and}\qquad
    \inner{\nu}{\nu} =1
\]
hold, where $\inner{~}{~}$ is the Euclidean inner product of $\R^3$.
We call $\nu$ the {\em unit normal vector field\/} of the front $f$.
Then there exists a $C^{\infty}$-function $\lambda\in C^{\infty}(U)$
such that
\begin{equation}\label{eq:lambda}
     \frac{\partial f}{\partial u}(u,v) \times
     \frac{\partial f}{\partial v}(u,v)  =
     \lambda (u,v)\, \nu(u,v),
\end{equation}
where $\times$ denotes the cross product of $\R^3$.
Obviously, $(u,v)\in U$ is a singular point of $f$ if and only if
$\lambda(u,v)=0$.
\begin{proposition}\label{prop:non-deg}
 A singular point $p\in U$ of a front $f\colon{}U\to\R^3$
 is non-degenerate if and only if $d\lambda\neq 0$ holds at $p$.
\end{proposition}
\begin{proof}
Differentiating \eqref{eq:lambda} at $p$, we have
$d(f_u \times f_v)(p)= d \lambda(p) \nu(p)$.
This implies that the rank of 
$d(f_u \times f_v) \colon U \to \R^3$ at $p$
is at most 1, and that $d \lambda(p) \ne 0$ is equivalent to
$d(f_u \times f_v)(p) \ne 0$. Therefore, $d \lambda(p) \ne 0$ if
and only if $d(f_u \times f_v)$ has rank 1 at $p$, 
that is, the map $\mathcal G$ in the introduction has rank 1 at $p$.
\end{proof}
Let $p$ be a non-degenerate singular point of a front $f\colon{}U\to
\R^3$.
Since the set of singular points is the set $\{\lambda=0\}$,
Proposition~\ref{prop:non-deg} implies that
the set of singular points is parametrized by a smooth curve
\begin{equation*}
    \gamma\colon{}(-\varepsilon,\varepsilon)\longrightarrow U
\end{equation*}
in a neighborhood of $p$, so that $\gamma(0)=p$.
We call the curve $\gamma(t)$ a {\em singular curve\/} passing
through $p$, and the direction $\dot\gamma(0)$ the {\em singular
direction\/} at the singular point $p$, where the dot denotes
the derivative with respect to $t$.
Since $p$ is a non-degenerate singular point, so is any point $\gamma(t)$
for sufficiently small $t$.
Then there exists a unique direction $\eta(t)\in T_{\gamma(t)}U$
up to scalar multiplication such that $df\bigl(\eta(t)\bigr)=0$ for each
$t$.
We call $\eta(t)$, which is smooth in $t$, the {\em null direction}.
\begin{definition}\label{def:C-S-type}
 Let $p$ be a non-degenerate singular point of a front
 $f\colon{}U\to\R^3$,
 $\gamma(t)$ the singular curve with $\gamma(0)=p$, 
 and $\eta(t)$ the null direction.
 Then
 \begin{enumerate}
  \item $p$ is of {\em type C\/} if $\eta(0)$ is not proportional
	to $\dot\gamma(0)$.
  \item $p$ is of {\em type S\/} if $\eta(0)$ is proportional to
	$\dot\gamma(0)$ and
	\[
	    \left.\frac{d}{dt}\right|_{t=0}
	         \det \bigl(\dot\gamma(t),\eta(t)\bigr)\neq 0
	\]
	holds,
	where $\dot\gamma(t)$ and $\eta(t)$ are considered as
	column vectors in $\R^2$.
 \end{enumerate}
\end{definition}
This definition does not depend on the choices of $\gamma$ and $\eta$.

\begin{example}\label{ex:canonical-singularities}
 The map
 \[
   \hat f_{\mathrm C}(z,w):=(2 w^3,-3 w^2,z)
 \]
 gives a cuspidal edge along the $z$-axis.
 The null direction is perpendicular to the $z$-axis, and it has a type C
 singularity at $(0,0)$.  
 
 The map 
 \[
    \hat f_{\mathrm S}(z,w):=(3w^4+zw^2,4w^3+2wz,z)
 \]
 gives a swallowtail at $(0,0)$.
 The singular curve is $6w^2+z=0$, and the null direction is parallel to
 the $w$-axis.
 So $(0,0)$ is a singularity of type S.
\end{example}
The above 
$\hat f=(\hat f_1,\hat f_2,\hat f_3)=\hat f_{\mathrm C}$,
$\hat f_{\mathrm S}$ satisfy that
$\hat f_z(0,0)=(0,0,1)$,
$\hat f_3(z,w)=z$ and
the derivative $\hat f_w$ vanishes identically on the singular
curve.
We shall now prove that any front $f(u,v)$ can
be given such a parameterization $(z,w)$ near a non-degenerate 
singular point, as follows:
We assume that the origin $(0,0)$ of the $uv$-plane is an 
arbitrarily given non-degenerate singular point of $f$, namely
\begin{equation*}
 \lambda(0,0)=0\qquad \text{and}\qquad
 d\lambda(0,0)\neq 0,
\end{equation*}
and set
\[
     f(0,0)=\mb 0.
\]
Then we have: 
\begin{proposition}\label{prop:normalization}
 Suppose that $(0,0)$ is a non-degenerate singular point of a front
 $f\colon{}U\to \R^3$.
 Then there is a diffeomorphism
 \[
     \Psi\colon{}(V;z,w)\longmapsto (U;u,v)
 \]
 with $\Psi(0,0)=(0,0)$ and a rotation at the origin
 \[
     \Theta\colon{}\R^3\longrightarrow \R^3
 \]
 such that
 \[
    \hat f(z,w)
         =\left(\hat f_1(z,w),\hat f_2(z,w),\hat f_3(z,w)\right)
    =\Theta\circ f\circ\Psi(z,w)\colon V\longrightarrow\R^3
 \]
 satisfies the following properties{\rm :}
 \begin{enumerate}
  \item\label{item:normalize-1}
        $\hat f_z(0,0)=(0,0,1)$,
  \item\label{item:normalize-3}
        $\hat f_3(z,w)=z$,
  \item\label{item:normalize-4}
        the derivative $\hat f_w$ vanishes identically along the singular
	curve. In particular $\hat f_w(0,0)=(0,0,0)$ holds.
  \item\label{item:normalize-5}
        If $(0,0)$ is of type C, the tangent vector 
        $\frac{\partial}{\partial z}\in T_{(0,0)}V$
	can be chosen to be the singular direction at the origin of $V$.
 \end{enumerate}
\end{proposition}
\begin{proof}
 Let $\gamma(t)$ be the singular curve passing through $(0,0)$.
 The null direction $\eta(t)$ can be extended to a vector field
 $\tilde\eta$ on $U$, that is, 
 \[
     \eta (t) = \tilde\eta\circ\gamma(t).
 \]
 On the other hand, we take a vector $\xi_0\in T_{(0,0)}U$ which
 is not proportional to $\eta(0)$ and satisfies 
 \[
       |df(\xi_0)|=1.
 \]
 If $(0,0)$ is of type C, we choose $\xi_0$ to be 
proportional to $\dot\gamma(0)$.
 Then there exists a vector field $\xi$ on $U$ such that
 \[
      \xi(0,0) = \xi_0.
 \]
 The vector fields $\xi$ and $\tilde\eta$ are linearly independent
 in a neighborhood of the origin.
 Hence by a lemma in \cite[page 182]{KN}, there exists a new coordinate
 system $(\tilde u,\tilde v)$ such that
 $\tilde u(0,0)=\tilde v(0,0)=0$ 
 and $\partial/\partial \tilde u$ (resp.\ $\partial/\partial \tilde v$) 
 is proportional to $\xi$ (resp.\ $\tilde\eta$).
 Scaling $\xi$ and $\tilde \eta$, we may assume
 \[
    \frac{\partial}{\partial\tilde u} = \xi\qquad
    \text{and}\qquad
    \frac{\partial}{\partial\tilde v} = \tilde\eta,
 \]
 without loss of generality.
 From now on, we use the coordinates $(\tilde u,\tilde v)$.
 However, for notational simplicity, we drop the overhead tilde's 
 and write $(\tilde u,\tilde v)$ as just $(u,v)$. So we may assume:
 \begin{itemize}
  \item The derivative $f_v$ vanishes identically on the singular curve
	$\gamma(t)$.
  \item If $(0,0)$ is of type C, the tangent vector 
        $\frac{\partial}{\partial u}\in T_{(0,0)}U$
	points in the singular direction at the origin.
 \end{itemize}
 Since $f_u(0,0)$ has unit length, we can take a 
 rotation at the origin $\Theta\colon{}\R^3\to\R^3$ which
 maps $f_u(0,0)$ to $(0,0,1)$, and set
 \[
    \tilde f(u,v)=\left(\tilde f_1(u,v),
                        \tilde f_2(u,v),\tilde f_3(u,v)\right)
                =\Theta\circ f(u,v).
 \]
 Then we have
 \[
     \tilde f_u(0,0)=(0,0,1),\qquad \tilde f_v(0,0)=(0,0,0) .
 \]
 We set
 \[
      g(u,v,z):= \tilde f_3(u,v)-z.
 \]
 Since
 \[
      g_u(0,0,0) = (\tilde f_3)_u(0,0)=1\neq 0,
 \]
 there exists a function $u=u(z,v)$ such that $u(0,0)=0$ and
 $g\bigl(u(z,v),v,z\bigr)=0$,
 namely,
 \begin{equation}\label{eq:implicit}
  \tilde f_3\bigl(u(z,v),v\bigr)=z.
 \end{equation}
 Then by
 \[
     u=u(z,w),\qquad v=w,
 \]
 $(z,w)$ gives a new coordinate system.
 We now set
 \begin{equation}\label{eq:def-normal}
  \hat f(z,w):=\tilde f\bigl(u(z,w),w\bigr).
 \end{equation}
 Then \ref{item:normalize-3} follows immediately.
 By differentiating \eqref{eq:implicit}, we have
 \[
     u_z(0,0)(\tilde f_3)_u(0,0)=1,
 \]
 and we get
 \[
       u_z(0,0)=1.
 \]
 Thus, by differentiating \eqref{eq:def-normal}, we have
 \[
     \hat f_z(0,0)=u_z(0,0)\tilde f_u(0,0)=(0,0,1),
 \]
 which implies \ref{item:normalize-1}.

 On the other hand,
 \[
    \hat f_w(z,w)=
    \left(\tilde f\bigl(u(z,w),w\bigr)\right)_w =
    \tilde f_u \bigl(u(z,w),w\bigr)u_w(z,w)+
    \tilde f_w \bigl(u(w,z),w\bigr).
 \]
 Since $\tilde f_v$ vanishes on $\gamma$, so does
 $\tilde f_w\bigl(u(z,w),w\bigr)$.
 Thus we have
 \[
     \hat f_w(z,w) = \tilde f_u\bigl(u(z,w),w\bigr)u_w(z,w)
\quad \mbox{on $\gamma$}.
 \]
 By differentiating $\tilde  f_3\bigl(u(z,w),w\bigr)=z$ with
 respect to $w$, we have
 \[
     u_w(z,w)(\tilde f_3)_u\bigl(u(z,w),w\bigr)=0.
 \]
 Here, $(\tilde f_3)_u\bigl(u(z,w),w\bigr)$ does not vanish near 
 $(0,0)$, since $(\tilde f_3)_u(0,0)=1$.
 Then we have
 \[
     u_w(z,w)=0
 \]
 and thus $\hat f_w$ vanishes on the singular curve,
 which proves \ref{item:normalize-4}.
 If $(0,0)$ is of type C, then $(0,0,1)$ is proportional to the singular
 direction of $\hat f$.
 Since $\hat f_z(0,0)=(0,0,1)$, we have \ref{item:normalize-5}.
\end{proof}

\begin{remark}
 In the proof above,
 \[
     \hat \nu(z,w):=\nu\bigl(u(z,w),w\bigr)
 \]
 gives the unit normal vector field of the (normalized) front $\hat f(z,w)$.
\end{remark}
In addition to the case of surfaces, we shall define fronts for
plane curves.

\begin{definition}\label{def:planar-front}
 Let $I\subset\R$ be an interval.
 A map
 \[
    \sigma=\sigma(w)\colon I\longrightarrow \R^2
 \]
 is called a {\em {\rm(}planar{\rm )} front\/} if there exists a
 map
 \[
     n=n(w)\colon{} I\longrightarrow S^1\subset \R^2
 \] 
 such that $n(w)$ is perpendicular to $\sigma(w)$ and
 $w\mapsto \bigl(\sigma(w),n(w)\bigr)$ is an immersion. 
 A point $w=w_0$ with $\sigma'(w_0)=0$ is called a singular point of
 the planar front $\sigma(w)$,
 where $'=d/dw$.
 At a singular point $w_0$, $n'(w_0)\neq 0$ holds by definition.
\end{definition}
A planar front is a projection of a Legendrian immersion in the unit
cotangent bundle $T_1^*\R^2$ with respect to the canonical contact
structure.

In the cases of $\hat f_{\mathrm C}$ and $\hat f_{\mathrm S}$
in Example~\ref{ex:canonical-singularities}, one can easily check
that their slices $\sigma^z:w \mapsto \hat f(z,w)$ 
perpendicular to the $z$-axis give planar fronts.
The tangent line of $\sigma^z(w)$ is given by
\begin{equation}\label{eq:tangents}
 \Phi(w,x,y,z):=n_1(z,w)\bigl(x-\hat f_1(z,w)\bigr)+
 n_2(z,w)\bigl(y-\hat f_2(z,w)\bigr)=0, 
\end{equation}
where $n(z,w)=\left(n_1(z,w),n_2(z,w)\right)$ 
is the unit normal vector of
$\sigma^z(w)$.
Then $\sigma^z$ is the envelope of this family of
tangent lines, and the discriminant set
$D_\Phi$
characterizes the image of $\hat f_{\mathrm C}$ and $\hat f_{\mathrm S}$.
According to this observation,
we shall prove that $\hat f(z,w)$ as in 
Proposition \ref{prop:normalization}, 
which has type C or type S singularities at $(0,0)$, 
also satisfies that
\begin{enumerate}
\item[(a)] the slice perpendicular to the $z$-axis
gives a planar front,
\item[(b)] the set $D_\Phi$ of $\Phi$ given by \eqref{eq:tangents} 
is a discriminant
set and is locally diffeomorphic to a cuspidal edge or
a swallowtail, by applying Fact \ref{fact:descriminant}.
\end{enumerate}
Now we shall prove (a) for non-degenerate singular points as follows:
\begin{proposition} \label{prop:slice}
 Let $(0,0)$ be a non-degenerate singular point of the front
 \[
    \hat f = \hat f(z,w)\colon{}V\longrightarrow \R^3
 \]
 satisfying \ref{item:normalize-3} and \ref{item:normalize-4}
 in Proposition~\ref{prop:normalization}.
 Then there exists an $\varepsilon>0$ such that the map defined by
 \[
   \sigma^z \colon  w\longmapsto \bigl(\hat f_1(z,w),\hat f_2(z,w)\bigr)
        \qquad \bigl(|z|<\varepsilon\bigr)
 \]
 is a planar front.
 Moreover, $w$ is a singular point of $\sigma^z$ if $(z,w)$
 is  a singular point of $\hat f$.
\end{proposition}
To prove Proposition~\ref{prop:slice}, we need the following:
\begin{lemma}\label{lem:after}
 Under the assumptions of Proposition~\ref{prop:slice},
 the derivative $\hat \nu_w(0,0)$ is non-zero and
 perpendicular to $\vect{e}_3:=(0,0,1)$,
 where $\hat\nu$ is the unit normal vector field of $\hat f$.
\end{lemma}
\begin{proof}
 Since $\hat f$ is a front, it follows from \ref{item:normalize-4}
 in Proposition~\ref{prop:normalization}
 that $\hat \nu_w(0,0)$ does not vanish.
 Since $\inner{\hat f_w}{\hat\nu}=0$, we have
 \[
    0 = \inner{\hat f_w}{\hat\nu}_z=
      \inner{\hat f_{wz}}{\hat\nu}+\inner{\hat f_w}{\hat\nu_z}.
 \]
 Since $\hat f_w(0,0)=0$ by \ref{item:normalize-4} 
in Proposition~\ref{prop:normalization}, 
we have
 \[
     \inner{\hat f_{wz}(0,0)}{\hat\nu(0,0)}=0.
 \]
 Thus,
 \begin{align*}
    \inner{\hat\nu_w(0,0)}{\vect{e}_3}&=
    \inner{\hat\nu_{w}(0,0)}{\hat f_z(0,0)} \\
  &=
    \left.\frac{\partial}{\partial w}\right|_{(z,w)=(0,0)}
           \inner{\hat \nu}{\hat f_z}-
            \inner{\hat\nu(0,0)}{\hat f_{wz}(0,0)}=0 \; , 
 \end{align*}
 which is the desired conclusion.
\end{proof}
\begin{proof}[Proof of Proposition~\ref{prop:slice}]
 We fix $z$ and let
 \[
   \sigma(w):=\bigl(\hat f_1(z,w),\hat f_2(z,w),0\bigr).
 \]
 Then $\sigma(w)$ is a map into the $xy$-plane.
 By \ref{item:normalize-3}  in Proposition~\ref{prop:normalization},
 we have
 \[
    \sigma(w)=\hat f(z,w) - \inner{\hat f(z,w)}{\vect{e}_3}\vect{e}_3
             = \hat f(z,w)-\hat f_3(z,w)\vect{e}_3
             =\hat f(z,w)-z\vect{e}_3
 \]
 and
 \[
     \frac{d}{dw}\sigma(w)=\hat f_w(z,w).
 \]
 This implies that a singular point of $\sigma$ is a singular point of
 $\hat f$.

 On the other hand, we set
 \[
     n(w):=
    \frac{\hat \nu-\inner{\hat\nu(z,w)}{\vect{e}_3}\vect{e}_3}{%
         \left(1-\inner{\hat\nu(z,w)}{\vect{e}_3}^2\right)^{1/2}}.
 \]
 Since
 \[
      \inner{\hat\nu(0,0)}{\vect{e}_3}=
      \inner{\hat\nu(0,0)}{\hat f_z(0,0)}=0,
 \]
 $n(w)$ is a well-defined unit vector field near $(0,0)$.
 Moreover,
 \begin{align*}
  \inner{\sigma_w(w)}{n(w)}&=
  \inner{\hat f_w(z,w)}{n(w)}=
  \frac{\inner{\hat\nu}{\hat f_w}-\inner{\hat\nu}{\vect{e}_3}
                                   \inner{\vect{e}_3}{\hat f_w}}{%
              \left(1-\inner{\hat\nu(z,w)}{\vect{e}_3}^2\right)^{1/2}}\\
    &=-\frac{\inner{\hat\nu}{\vect{e}_3}}{%
              \left(1-\inner{\hat\nu(z,w)}{\vect{e}_3}^2\right)^{1/2}}
       (\hat f_3)_w =0,
 \end{align*}
 where we used the fact that $\hat f_3(z,w)=z$.
 Thus, $n(w)$ is a normal vector of $\sigma(w)$.
 By Lemma~\ref{lem:after}, we have
 $\inner{\hat\nu_w(0,0)}{\vect{e}_3}=0$,
 and
 \[
      \frac{d}{dw}n(0)=\hat \nu_w(0,0)\neq 0.
 \]
 Hence $n'(w)\neq 0$ for sufficiently small $(z,w)$, and the
 map $w\mapsto (\sigma(w),n(w))$ is an immersion.
\end{proof}

\subsection{Proof of the criteria}
\label{subsec:criterion}
In this section, we shall prove Proposition~\ref{prop:criterion-front}
in the introduction.
As pointed out in the beginning of the previous section, it is sufficient
to prove the assertion for fronts in the Euclidean $3$-space $\R^3$.
The idea of the proof is as follows:
Let $(0,0)$ be a non-degenerate singular point of a front
$f=f(u,v)\colon{}U\to\R^3$.
Then by Proposition~\ref{prop:normalization}, we have a normalized front
$\hat f(z,w)$.
We set
\[
    \sigma(z,w):=\bigl(\hat f_1(z,w),\hat f_2(z,w)\bigr).
\]
By Proposition~\ref{prop:slice}, there exist positive numbers
$\varepsilon_1$ and $\varepsilon_2$ such that
\[
    (-\varepsilon_1,\varepsilon_1)\ni w \longmapsto
       \sigma(z,w)\in \R^2
\]
gives a planar front for $|z|<\varepsilon_2$; 
that is, there exists a unit normal vector field
\[
     n= n(z,w)\colon(-\varepsilon_2,\varepsilon_2)\times
                    (-\varepsilon_1,\varepsilon_1)
		    \longrightarrow \R^2
\]
such that
$
     \inner{\sigma_w(z,w)}{n(z,w)}=0.
$
If we set $n=(n_1,n_2)$, the equation
\[
    n_1(z,w)\bigl(x-\hat f_1(z,w)\bigr)+
    n_2(z,w)\bigl(y-\hat f_2(z,w)\bigr)=0
\]
gives the tangent line of the planar front $w\mapsto \sigma(z,w)$, and
the image of the planar front is the envelope of these tangent lines.
On the other hand, it is well-known that the envelope generated by a
family of lines
\[
   \left\{ F(w,x,y)=0\,;\, w\in\R\right\}
\]
is given by $\{(x,y)\,;\,F(w,x,y)=F_w(w,x,y)=0, w\in\R\}$.
So if we set
\begin{equation}\label{eq:Phi}
 \Phi(w,x,y,z):= n_1 (z,w)\bigl(x-\hat f_1(z,w)\bigr)+
                 n_2 (z,w)\bigl(y-\hat f_2(z,w)\bigr),
\end{equation}
the discriminant set
\[
   D_{\Phi}:=\{\vect{x}\in\Omega\,;\,
    \text{there exists a $w\in\R$ with
    $\Phi(w,\vect{x})=\Phi_w(w,\vect{x})=0$}\}
\]
coincides with the image of the front $\hat f$.
Now we set
\[
    \varphi(w)=\Phi(w,0,0,0).
\]
Then if $\varphi(w)$ has an $\A_k$-singularity ($k=2,3$) and $\Phi$
is a versal unfolding, we can conclude (b), that is,
the image of $\hat f$ is locally diffeomorphic to a cuspidal edge or a
swallowtail, by Fact~\ref{fact:descriminant}.
According to this plan, we shall first prove the criterion for cuspidal 
edges.  First, we prepare three lemmas:
\begin{lemma}\label{lem:cusp-1}
 $\sigma(0,0)=\sigma'(0,0)=\sigma_z(0,0)=(0,0)$ and 
$n'(0,0)\neq (0,0)$ hold,
where $'$ denotes the derivative with respect to $w$.
\end{lemma}
\begin{proof}
 These are easily computed from Proposition \ref{prop:normalization}
 and \ref{prop:slice}. 
\end{proof}
\begin{lemma}\label{lem:cusp-2}
 $\Phi_z(0,0,0,0)=\Phi'_{z}(0,0,0,0)=0$.
\end{lemma}
\begin{proof}
These are computed by differentiating \eqref{eq:Phi} and by 
using Lemma \ref{lem:cusp-1}.
\end{proof}
\begin{lemma}\label{lem:cusp-3}
 $\varphi(w)$ has an $\mathcal A_2$-singularity at $w=0$ if and only if
 $\sigma''(0,0)\ne (0,0)$.
\end{lemma}
\begin{proof}
 By differentiating $\inner{\sigma'}{n}=0$ and using
$\sigma'(0,0)=(0,0)$, we have \\
$\inner{\sigma''(0,0)}{n(0,0)}=0$.
Since $n'(0,0)\neq (0,0)$ by Lemma~\ref{lem:cusp-1},
 $\{n,n'/|n'|\}$ forms an orthonormal basis for $\R^2$.
Therefore,
 \begin{equation}\label{eq:sigma''}
     \sigma''(0,0)=
     \inner{\sigma''(0,0)}{n'(0,0)}n'(0,0)/|n'(0,0)|^2.
 \end{equation}

On the other hand, by differentiating \eqref{eq:Phi} and
by using Lemma \ref{lem:cusp-1}, we have
 \[
     \varphi(0)=\varphi'(0)=\varphi''(0)=0,\quad
     -\varphi'''(0)=\inner{\sigma''(0,0)}{n'(0,0)}.
 \]
 Hence $\varphi$ has an $\A_2$-singularity at the origin if and only if
 $\inner{\sigma''(0,0)}{n'(0,0)}\neq 0$.
This and \eqref{eq:sigma''} prove the assertion.
\end{proof}

\begin{proposition}\label{prop:C-type}
 Suppose $(0,0)$ is a non-degenerate singular point.
Then
the germ of the image of the front is locally diffeomorphic to a 
 cuspidal edge if and only if $(0,0)$ is of type C.
\end{proposition}

\begin{proof}
By Lemma~\ref{fact:Z}, local diffeomorphic equivalence between singular
points on fronts implies Legendrian equivalence.
Since a cuspidal edge itself is of type C, any singular point
locally diffeomorphic to a cuspidal edge is of type C.
Conversely, we shall show that a singularity of type C 
is locally diffeomorphic to a cuspidal edge. 
To prove this, it is sufficient to show that  $\varphi$ as above
has an $\A_2$-singularity and
 $\Phi$ is versal.
 By Lemma~\ref{lem:cusp-2} and \eqref{eq:Phi}, we have
 \[
    \begin{pmatrix}
     \Phi_x(0,\vect{0}) & \Phi'_x(0,\vect{0}) \\
     \Phi_y(0,\vect{0}) & \Phi'_y(0,\vect{0}) \\
     \Phi_z(0,\vect{0}) & \Phi'_z(0,\vect{0})
    \end{pmatrix}=
    \begin{pmatrix}
     n_1(0,0) & n_1'(0,0) \\
     n_2(0,0) & n_2'(0,0) \\
        0     & 0
    \end{pmatrix}.
 \]
 This matrix is of rank $2$, since $n(0,0)$ and $n'(0,0)$ are linearly
 independent.

 Next we prove that $\varphi$ has an $\A_2$-singularity.
 We set
 \[
     \lambda =\det \bigl(\hat f_w,\hat f_z,\hat\nu\bigr),
 \]
 where $\hat\nu$ is the unit normal vector of the front $\hat f$.
 Here $\lambda=0$ on the singular curve.
 Since we have assumed that $(0,0)$ is of type C,
 \ref{item:normalize-5}
 in Proposition~\ref{prop:normalization}
 implies that
 $\partial/\partial z$ is
 the singular direction at the origin.
 So we have
 \[
     \lambda_z(0,0)=0.
 \]
 On the other hand,
 \[
   0 \ne     \lambda_w=
     \det\bigl(\hat f_{ww},\hat f_{z},\hat\nu\bigr) +
     \det\bigl(\hat f_{w},\hat f_{zw},\hat\nu\bigr) +
     \det\bigl(\hat f_{w},\hat f_{z},\hat\nu_w\bigr)
     = \det\bigl(\hat f_{ww},\hat f_{z},\hat\nu\bigr) 
 \]
holds at $(0,0)$, because $\hat f_w(0,0)=\boldsymbol 0$.
 Since $\hat f_z(0,0)\times \hat\nu(0,0)$ is parallel to the $xy$-plane,
 we have
 \[
  \inner{\sigma''(0,0)}{\hat f_z(0,0)\times\hat\nu(0,0)}=
  \inner{\hat f_{ww}(0,0)}{\hat f_z(0,0)\times \hat \nu(0,0)}\neq 0.
 \]
 In particular we have $\sigma''(0,0)\neq (0,0)$, and by
 Lemma~\ref{lem:cusp-3}, $\varphi(w)$ has an $\A_2$-singularity at
 $w=0$.
\end{proof}

Next, we prove the criterion for swallowtails:
\begin{proposition}\label{prop:S-type}
 Suppose $(0,0)$ is a non-degenerate singular point.
 Then the germ of the image of the front is locally diffeomorphic
 to a swallowtail if and only if $(0,0)$ is of type S.
\end{proposition}

To prove this, we prepare a lemma:
\begin{lemma}\label{lem:swallow-1}
 Suppose $(0,0)$ is a non-degenerate singular point of $\hat f(z,w)$, 
 but not of type C.
 Then
 \begin{enumerate}
  \item $\hat f_{ww}(0,0)=\vect{0}$, in particular $\sigma''(0,0)=(0,0)$, and 
  \item $\Phi_z''(0,0,0,0)\ne 0$.  
 \end{enumerate}
\end{lemma}
\begin{proof}
 By \ref{item:normalize-4} of Proposition~\ref{prop:normalization},
 $\hat f_w$ vanishes identically on the singular curve.
 Since $(0,0)$ is not of type C,
 the singular direction is equal to the null direction
 $\partial/\partial w$, thus $\hat f_{ww}(0,0)=\vect{0}$.
 In particular, we have $\sigma''(0,0)=(0,0)$.
 Differentiating \eqref{eq:Phi} by $w$ and $z$ and substituting
 the relation
 \[
    \sigma(0,0)=\sigma'(0,0)=\sigma_z(0,0)=\sigma''(0,0)=(0,0),
 \]
 we have
 \[
   - \Phi''_z(0,0,0,0)=\inner{n'(0,0)}{\sigma'_z(0,0)}.
 \]
 Since $\inner{n}{\sigma'}=0$, we have
 $\inner{n}{\sigma'_z}=0$.
 Since $n'$ is orthogonal to $n$,
 $\sigma'_z$ is proportional to $n'$.
 To show $\Phi''_z(0,0,0,0)\neq 0$, it is sufficient to show
 $\sigma'_z(0,0)\neq (0,0)$.
 Moreover, $(\hat f_3)_{wz}$ vanishes identically, so 
 $\sigma'_z(0,0)\neq (0,0)$ is equivalent to
 $\hat f_{zw}(0,0)\neq \vect{0}$.

 Differentiating 
 \[
    \lambda = \det \bigl(\hat f_w,\hat f_z,\hat\nu\bigr)
 \]
 with respect to $z$ and using the relation $\hat f_w(0,0)=\vect{0}$,
 we have
 \begin{equation}\label{eq:pr-sw-1}
    \lambda_z(0,0) = \det\bigl(\hat f_{wz}(0,0),
                          \hat f_{z}(0,0),
                          \hat\nu(0,0)\bigr).
 \end{equation}
 Since $(0,0)$ is not of type C, $\partial/\partial w$ is the singular
 direction.
 In particular $\lambda_w(0,0)=0$ holds.
 Since $d\lambda\neq 0$ at $(0,0)$, we have $\lambda_z(0,0)\neq 0$.
 Hence by \eqref{eq:pr-sw-1}, we have $\hat f_{zw}(0,0)\neq \vect{0}$.
\end{proof}
\begin{proof}[Proof of Proposition~\ref{prop:S-type}]
For the same reason as in the proof of Proposition~\ref{prop:C-type},
any singular point
locally diffeomorphic to a swallowtail is of type S.
Conversely, we shall show that a singularity of type S 
is locally diffeomorphic to a swallowtail.
To prove this, it is sufficient to show that $\varphi$ as above
has an $\A_3$-singularity and
$\Phi$ is versal.
By Lemma~\ref{lem:cusp-2} and \eqref{eq:Phi}, we have
 \[
   \begin{pmatrix}
    \Phi_x(0,\vect{0}) &
    \Phi'_x(0,\vect{0}) &
    \Phi''_x(0,\vect{0}) \\
    \Phi_y(0,\vect{0}) &
    \Phi'_y(0,\vect{0}) &
    \Phi''_y(0,\vect{0}) \\
    \Phi_z(0,\vect{0}) &
    \Phi'_z(0,\vect{0}) &
    \Phi''_z(0,\vect{0})
   \end{pmatrix}
   =
   \begin{pmatrix}
    n_1(0,0) & n_1'(0,0) & * \\
    n_2(0,0) & n_2'(0,0) & * \\
       0     & 0         &\Phi''_z(0,\vect{0})
   \end{pmatrix}.
 \]
 By Lemma~\ref{lem:swallow-1}, $\Phi_z''(0,\vect{0})\neq 0$, and then
 the rank of this matrix is $3$.
 By Lemma~\ref{lem:cusp-1} and Lemma~\ref{lem:swallow-1}, we have
 \[
    \varphi(0)=\varphi'(0)=\varphi''(0)=\varphi'''(0)=0,\qquad
    \varphi^{(4)}(0)=\inner{\sigma'''(0)}{n'(0)}.
 \]
 By differentiating $\inner{\sigma'}{n}=0$ twice, we have
 \[
    \inner{\sigma'''(0,0)}{n(0,0)}=0.
 \]
 Thus, $\varphi(w)$ has an $\A_3$-singularity if and only if
 $\sigma'''(0)\neq 0$,
 which is equivalent to $\hat f_{www}(0,0)\neq \vect{0}$, 
 since $\hat f_3(z,w)=z$.

 Since $\lambda(z,w)=\det\bigl(\hat f_w,\hat f_z,\hat \nu\bigr)$ and
 the singular curve $\gamma$ is given by $\lambda(z,w)=0$,
 the singular direction is given by
 \[
    \xi(z,w)=\bigl(\lambda_w(z,w), -\lambda_z(z,w)\bigr) 
        \quad \mbox{on $\gamma$}.
 \]
 On the other hand, the null direction $\eta$ is given by
 \[
     \eta(z,w)=\frac{\partial}{\partial w} \quad
\mbox{on $\gamma$}.
 \]
 Since $(0,0)$ is of type S, we have
 \begin{equation}\label{eq:s-type-cond}
    0\neq \left.\frac{d}{dw}\right|_{w=0}
    \det\bigl(\xi(z,w),\eta(z,w)\bigr)=\lambda_{ww}(0,0).
 \end{equation}
 By the definition of $\lambda$, we have
 \[
    \lambda_w =
       \det \bigl(\hat f_{ww},\hat f_z,\hat \nu\bigr) +
       \det \bigl(\hat f_{w},\hat f_{zw},\hat \nu\bigr) +
       \det \bigl(\hat f_{w},\hat f_{z},\hat \nu_w\bigr) .
 \]
 Here, since $\hat f_w\times \hat f_z$ is proportional to $\hat \nu$,
 we have
 \[
     \det \bigl(\hat f_w,\hat f_z,\hat\nu_w\bigr)
     =\inner{\hat f_w\times\hat f_z}{\hat\nu_w}=0.
 \]
 Thus, we have
 \[
   \lambda_w =
       \det \bigl(\hat f_{ww},\hat f_z,\hat \nu\bigr) +
       \det \bigl(\hat f_{w},\hat f_{zw},\hat \nu\bigr) .
 \]
 Then we get
 \begin{multline*}
    \lambda_{ww}=
       \det \bigl(\hat f_{www},\hat f_z,\hat \nu\bigr) +
       \det \bigl(\hat f_{ww},\hat f_z,\hat \nu_w\bigr) \\
      +
       2\det \bigl(\hat f_{ww},\hat f_{wz},\hat \nu\bigr) +
       \det \bigl(\hat f_{w},\hat f_{wwz},\hat \nu\bigr) +
       \det \bigl(\hat f_{w},\hat f_{wz},\hat \nu_w\bigr).
 \end{multline*}
 Using $\hat f_w(0,0)=\hat f_{ww}(0,0)=\vect{0}$ and 
 \eqref{eq:s-type-cond}, we have
 \[
    0\neq \lambda_{ww}(0,0) =
    \det\bigl(\hat f_{www}(0,0),\hat f_z(0,0),\hat\nu(0,0)\bigr),
 \]
 which proves $\hat f_{www}(0,0)\neq \vect{0}$.
\end{proof}
Izumiya and Takeuchi \cite{IT} gave criteria for the
singularities of a non-cylindrical flat ruled front
in the Euclidean $3$-space $\R^3$ to be cuspidal edges, swallowtails and 
cuspidal cross caps.  
One can prove this for the case of cuspidal edges and swallowtails by 
directly applying our criteria.  
\section{Local properties of flat fronts in $H^3$}
\label{sec:local}
In this section, we give a proof of Theorem~\ref{thm:criterion}
in the introduction, and
show that, generically, singular points of flat fronts are
cuspidal edges or swallowtails.
\subsection{Preliminaries}

We denote by $\Lor^4$ the Minkowski $4$-space with
the inner product $\inner{~}{~}$ of signature $(-,+,+,+)$.
The hyperbolic $3$-space $H^3$ is considered as the upper
half component of the two sheet hyperboloid in $\Lor^4$
with the metric induced by $\inner{~}{~}$.
Identifying $\Lor^4$ with $\Herm(2)$, the set of $2\times 2$-hermitian
matrices, as
\[
      \Lor^4\ni (x_0,x_1,x_2,x_3)\leftrightarrow
      \begin{pmatrix}
        x_0+x_3 & x_1+\sqrt{-1}x_2 \\
        x_1-\sqrt{-1}x_2 & x_0-x_3
      \end{pmatrix}\in\Herm(2),
\]
one has $\inner{X}{X}=-\det X$ for $X\in\Herm(2)$, and
$H^3$ is represented as
\begin{align*}
     H^3 &= \{x=(x_0,x_1,x_2,x_3)\in\Lor^4\,;\,
             \inner{x}{x}=-1,x_00\}\\
         &= \{ X\in\Herm(2)\,;\,\det X=1,~\trace X0\}\\
         &= \{ aa^*\,;\,a\in\SL(2,\C)\}=\SL(2,\C)/\SU(2).
\end{align*}
The tangent space of $H^3$ at $p\in H^3$ is the set of 
vectors in $\Lor^4$ which are perpendicular to $p$:
\[
     T_p H^3 = \{Y\in \Lor^4\,;\,\inner{p}{Y}=0\}.
\]
We define a bilinear, skew-symmetric product $\times$ as
\begin{equation}\label{eq:hyp-cross}
     X\times Y:=\frac{\sqrt{-1}}{2} \bigl(Xp^{-1}Y-Yp^{-1}X\bigr)
     \qquad \text{for}\quad
     X,Y\in T_pH^3,
\end{equation}
where $X$, $Y$ and $p$ are considered as matrices in $\Herm(2)$,
and the products of the right-hand side are matrix multiplications.
It is easy to show that $X\times Y$ is a vector in $T_p H^3$ and 
perpendicular to both $X$ and $Y$.
We call ``$\times$'' the {\em cross product\/} of $H^3$.

Let $M^2$ be an oriented simply-connected Riemannian $2$-manifold, and let 
\[
    f \colon M^2\longrightarrow H^3=\SL(2,\C)/\SU(2)
\]
be a front whose Legendrian lift is
\[
    L_f \colon M^2\to T_1^* H^3=\SL(2,\C)/\U(1).
\]
Identifying $T_1^*H^3$ with $T_1H^3$, we can write $L_f=(f,\nu)$,
where $\nu(p)$ is a unit vector in $T_pH^3$ such that
$\inner{df(p)}{\nu(p)}=0$ holds for each $p\in M^2$.
We call $\nu$ the {\em unit vector field\/} of the front $f$.

Suppose that $f$ is flat, then there is a (unique) complex structure on
$M^2$ and a holomorphic Legendrian immersion
\begin{equation}\label{eq:lift-eq}
  E_f\colon M^2\longrightarrow \SL(2,\C)
\end{equation}
such that $f$ and $L_f$ are projections of $E_f$, where being a 
holomorphic Legendrian map means that 
$E_f^{-1}dE_f$ is off-diagonal.
In particular, $f=E_f^{}E_f^*$, with $H^3$ considered to be as in 
\eqref{eq:hyperbolic}.
(See \cite{GMM}, \cite{KUY1} and \cite{KUY2} for details.)
If we set
\begin{equation}\label{eq:lift-can}
 E_f^{-1}dE_f=
 \begin{pmatrix}
  0 & \theta \\
  \omega & 0
 \end{pmatrix},
\end{equation}
the first and the second fundamental forms $ds^2$ and $dh^2$
are given by
\begin{equation}\label{eq:fund-forms}
\begin{array}{rl}
 ds^2&=\omega\theta+\bar \omega \bar \theta+
  (|\omega|^2+|\theta|^2), \\
 dh^2&=|\theta|^2-|\omega|^2.
\end{array}
\end{equation}
We call $\omega$ and $\theta$ the {\em canonical forms\/} of the 
front $f$.
The holomorphic $2$-differential
\begin{equation}\label{eq:hopf}
   Q:=\omega\theta \; , 
\end{equation}
which appears in the $(2,0)$-part of $ds^2$, is called the 
{\it Hopf differential\/} of $f$.
By definition, the umbilic points of the front $f$ coincide with the
zeroes of $Q$.  We remark that the $(1,1)$-part of the first
fundamental form 
\begin{equation}\label{eq:one-one-part}
     ds^2_{1,1}:=|\omega|^2+|\theta|^2
\end{equation}
is positive definite on $M^2$.

Conversely, the following assertion holds (see \cite{KUY2} for the front
case and \cite{GMM} for the regular case):
\begin{fact}
 Let $\omega$ and $\theta$ be holomorphic $1$-forms on a simply-connected 
 Riemann surface $M^2$  such that $|\omega|^2+|\theta|^2$ is
 positive definite.
 Then the solution of the ordinary differential equation
 \[
    E^{-1}dE=
     \begin{pmatrix}
      0 & \theta \\
      \omega & 0
     \end{pmatrix},\qquad
   E(z_0)=     \begin{pmatrix}
      1 & 0 \\
      0 & 1
     \end{pmatrix}
 \]
 gives a holomorphic Legendrian immersion of $M^2$ into $\SL(2,\C)$,
 where $z_0\in M^2$ is a base point, and its projection into $H^3$ gives
 a flat front. Conversely, 
 any flat front is locally expressed in this manner.
\end{fact}
\begin{remark}\label{rem:parallel}
 If we identify $H^3$ with the upper-half component of the hyperboloid
 in Minkowski $4$-space $\Lor^4$,
 the parallel surface of $f$ is written as
 \[
    f_t=(\cosh t)f +(\sinh t)\nu \colon M^2\to H^3 \subset \Lor^4 \; , 
 \]
 where $t$ is the signed distance from $f$ and $\nu$ is the unit normal
 vector of $f$ in $H^3$.
 As pointed out in \cite{GMM} and \cite{KUY2},
 \begin{equation}\label{eq:parallel-lift}
    E_{f_t}=
       E_f \begin{pmatrix}
	    e^{t/2} & 0 \\ 0 & e^{-t/2}
	   \end{pmatrix}.
 \end{equation}%
 Then the canonical forms $\omega_t$ and $\theta_t$ of $f_t$ are written
 as
 \begin{equation}\label{eq:parallel-canonical}
    \omega_t = e^{t} \omega,\qquad
    \theta_t = e^{-t}  \theta.
 \end{equation}
\end{remark}
\subsection{Proof of Theorem~\ref{thm:criterion}}
 Let $f\colon{}M^2\to H^3$ be a flat front.
 Then, on a neighborhood of $p$, we can take a holomorphic Legendrian
 immersion $E_f$ as in \eqref{eq:lift-eq}.
 Since $ds^2_{1,1}=|\omega|^2+|\theta|^2$
 is positive definite,
 it holds that either $\omega(p)\neq 0$ or $\theta(p)\neq 0$.
 So, by \eqref{eq:lift-eq} and the fact that $f=E_f^{}E_f^*$, we have
 \[
       f^{-1}df=
          (E_f^*)^{-1}\left\{
                   \begin{pmatrix}
		    0 & \theta \\
		    \omega & 0
                   \end{pmatrix}+
                   \begin{pmatrix}
		    0 & \bar\omega \\
		    \bar\theta & 0
                   \end{pmatrix}
                 \right\}E_f^*.
 \]
 Thus, if we write $\omega=\hat\omega\,dz$ and $\theta=\hat\theta\,dz$ 
 in a complex coordinate $z$, we have
 \[
   f^{-1}f_z=
          (E_f^*)^{-1}
                   \begin{pmatrix}
		    0 & \hat\theta \\
		    \hat\omega & 0
                   \end{pmatrix}\,
                 E_f^* \;\;\;\; \text{and} \;\;\;\; 
   f^{-1}f_{\bar z}=
          (E_f^*)^{-1}
                   \begin{pmatrix}
		    0 & \bar{\hat\omega} \\
		    \bar{\hat\theta} & 0
                   \end{pmatrix}\,
                 E_f^* \; , 
 \]
 and then
 \[
      (f^{-1}f_z)\times (f^{-1}f_{\bar z})=
           \left(|\hat\theta|^2-|\hat\omega|^2\right)
           \begin{pmatrix}
                      1 & \hphantom{-}0 \\
                      0 & -1
           \end{pmatrix},
 \]
 where ``$\times$'' is the cross product as in \eqref{eq:hyp-cross}.
 Thus, the singular set is the set of zeroes of the function
 \[
         \lambda=|\hat\theta|^2-|\hat\omega|^2.
 \]
 Then $p$ is a singular point if and only if
 \begin{equation}\label{eq:singular}
  |\hat\omega(p)|=|\hat\theta(p)|.
 \end{equation}
 Hence (1) is proven.
 Since $f$ is a front, $ds^2_{1,1}$ as in \eqref{eq:one-one-part}
 is positive definite.
 Hence $|\hat\omega(p)|=|\hat\theta(p)|\neq 0$ holds on a
 singular point $p$.

 Moreover, at a  singular point $p$, we have
 \begin{align*}
    d\lambda &= d\left(\hat\theta\bar{\hat\theta}-
         \hat\omega\bar{\hat\omega}\right)
      =\left(\hat{\theta}'\bar{\hat\theta}-{\hat\omega}'\bar{\hat\omega}\right)
        \,dz +
      \left(\hat\theta\bar{\hat\theta}'-\hat\omega\bar{\hat\omega}'\right)
       \,d\bar z
      \\
     &=\left\{
            \hat\theta'\frac{\hat\theta\bar{\hat\theta}}{\hat\theta}-
            \hat\omega'\bar{\hat\omega}
       \right\}\,dz +
       \overline{
      \left\{
            \hat\theta'\frac{\hat\theta\bar{\hat\theta}}{\hat\theta}-
            \hat\omega'\bar{\hat\omega}
       \right\}}\,d\bar z\\
     &=\left\{
            \hat\theta'\frac{\hat\omega\bar{\hat\omega}}{\hat\theta}-
            \hat\omega'\bar{\hat\omega}
       \right\}\,dz +
       \overline{
          \left\{
            \hat\theta'\frac{\hat\omega\bar{\hat\omega}}{\hat\theta}-
            \hat\omega'\bar{\hat\omega}
       \right\}}\,d\bar z \qquad (\text{by \eqref{eq:singular}})\\
     &= \frac{\bar{\hat\omega}}{\hat\theta}
        \left(\hat\theta'\hat\omega-\hat\omega'\hat\theta\right)\,dz+
        \frac{\hat\omega}{\bar{\hat\theta}}
       \overline{\left(\hat\theta'\hat\omega
         -\hat\omega'\hat\theta\right)}
        \,d\bar z.
 \end{align*}
 Hence a singular point $p$ is non-degenerate if and only if
 \begin{equation}\label{eq:non-deg}
  \hat\theta'\hat\omega-\hat\omega'\hat\theta\neq 0
 \end{equation}
 holds at $p$.

 Let $p$ be a non-degenerate singular point, that is, 
 \eqref{eq:singular} and \eqref{eq:non-deg} hold at $p$.
 Let $\gamma(t)$ be a singular curve such that $\gamma(0)=p$.
 Since $|\hat\theta|^2-|\hat\omega|^2=0$ holds on $\gamma(t)$,
 \begin{align*}
   0&=
    \left(\hat\theta'\bar{\hat\theta}-\hat\omega'{\bar{\hat\omega}}\right)
    \dot\gamma +
    \overline{
    \left(\hat\theta'\bar{\hat\theta}-\hat\omega'{\bar{\hat\omega}}\right)}
    \dot{\bar \gamma}\\
    &=
    |\hat\theta|^2
     \left(\frac{\hat\theta'}{\hat\theta}-
     \frac{\hat\omega'}{\hat\omega}\right)\dot\gamma +
    |\hat\theta|^2
    \overline{
     \left(\frac{\hat\theta'}{\hat\theta}-
     \frac{\hat\omega'}{\hat\omega}\right)}\dot{\bar{\gamma}}
    =2|\hat\theta|^2
     \left\langle
       \overline{
       \left(\frac{\hat\theta'}{\hat\theta}-
       \frac{\hat\omega'}{\hat\omega}\right)},\dot\gamma
     \right\rangle
 \end{align*}
 holds on $\gamma(t)$, because of \eqref{eq:singular},
 where $\langle~,~\rangle$ is the Hermitian inner product on $\C$.
 Hence
 $\sqrt{-1}\overline{
         \left({\hat\theta'}/{\hat\theta}-
         {\hat\omega'}/{\hat\omega}\right)}$
 gives the singular direction.
 Thus, by a suitable choice of the parameter $t$,
 the singular curve $\gamma(t)$ can be parametrized as
 \begin{equation}\label{eq:singular-direction}
    \dot\gamma(t)=\sqrt{-1}\overline{
         \left(\frac{\hat\theta'}{\hat\theta}-
         \frac{\hat\omega'}{\hat\omega}\right)}.
 \end{equation}

 The first fundamental form $ds^2$
 is written as
 \begin{align*}
   ds^2 =
     \left(\hat\omega\,dz+\overline{\hat\theta}\,d\bar z\right)
     \left(\hat\theta\,dz+\overline{\hat\omega}\,d\bar z\right)
 \end{align*}
 on the curve $\gamma(t)$.
 Now we set $\rho = \hat\theta/\hat\omega$.
 Since $\rho(p)\neq 0$, there exists a holomorphic function $g$
 defined on a neighborhood of $p$  such that $g^2=\rho$.
 Since $|g|=1$ on the singular curve $\gamma(t)$, we have
 \[
    \hat\omega\,\left(\frac{\sqrt{-1}}{g\hat\omega}\right)
    +
    \overline{\hat\theta}\,
    \overline{\left(\frac{\sqrt{-1}}{g\hat\omega}\right)}
    =\sqrt{-1}
      \left(\frac{1}{g}-\bar g\right)=0.
 \]
 Thus the null direction $\eta(t)$ is given by
 \begin{equation*}
    \eta(t) = \frac{\sqrt{-1}}{g\hat\omega}=
        \frac{\sqrt{-1}}{\sqrt{\hat\omega\hat\theta}}.
 \end{equation*}
 So we have
 \begin{equation}\label{eq:det}
    \det\bigl(\dot\gamma,\eta\bigr)
     =\Im \overline{\dot\gamma}\eta
     =\Im\left\{
       \left(\frac{\hat\theta'}{\hat\theta}-\frac{\hat\omega'}{\hat\omega}
       \right)\frac{1}{\sqrt{\hat\omega\hat\theta}}
      \right\}.
 \end{equation}
 Here, by Proposition~\ref{prop:criterion-front}, $p$ is diffeomorphic
 to a cuspidal edge
 if $\det\bigl(\dot\gamma,\eta\bigr)\neq 0$ at $t=0$.
 Hence, we have (2).

Next, let us prove (3). Using \eqref{eq:singular-direction} 
and \eqref{eq:det}, we can compute that 
\begin{equation*}
 \frac{d}{dt}\det(\dot\gamma,\eta) = 
\Re\left[
        \frac{s(\hat\theta)-s(\hat\omega)}{\sqrt{\hat\omega\hat\theta}}
        \vert \hat\omega\hat\theta \vert     
 \overline{\left(\frac{\hat\theta'/\hat\theta-\hat\omega'/\hat\omega}
 {\sqrt{\hat\omega \hat \theta}}\right)}
       \right].
\end{equation*}
Hence, if $\bigl(\hat\theta'/\hat\theta-\hat\omega'/\hat\omega\bigr) /%
\sqrt{\hat\omega \hat \theta} \in \R$ at $p$, then 
\begin{equation*}
 \left. \frac{d}{dt}\det(\dot\gamma,\eta) \right\vert_{t=0}= 
\left. \left\{
\vert \hat\omega\hat\theta \vert         
 \left(\frac{\hat\theta'/\hat\theta-\hat\omega'/\hat\omega}
 {\sqrt{\hat\omega \hat \theta}}\right)
\Re \left(
\frac{s(\hat\theta)-s(\hat\omega)}{\sqrt{\hat\omega\hat\theta}}
\right) 
\right\} \right\vert_p \; .
\end{equation*}
This proves (3), because of 
Proposition~\ref{prop:criterion-front}. 
\qed

\begin{remark}\label{rem:DescriptionG}
 We set the two hyperbolic Gauss maps to be 
 \[
   G=\frac{A}{C}, \quad G_*=\frac{B}{D}, \quad
 \text{ where }
    E_f=\begin{pmatrix} A & B \\ C & D \end{pmatrix}.
 \]
 As shown in \cite{KUY1}, we have the following expression
 \begin{equation}\label{eq:G-repr}
   E_f=\begin{pmatrix} 
	G/\varDelta & 
	 \varDelta G_*/(G-G_*) \\
	 1/\varDelta & \varDelta/(G-G_*)
       \end{pmatrix} \; , 
   \qquad \varDelta:=c \exp\left( \int_{z_0}^z\frac{dG}{G-G_*} \right),
 \end{equation}
 where $c=e^{-t/2}\in \R$ determines which member of the parallel
 family $f_t$  of $f$ we have.  
 In this $(G,G_*)$-construction of flat fronts,
 it is convenient to rewrite the conditions in
 Theorem~\ref{thm:criterion}  in terms of $(G,G_*)$.
 We have the following identities, which will be useful
 for an application of Theorem~\ref{thm:criterion} (see \cite{KUY2}):
 \begin{align}
  &Q:=\omega\theta=-\frac{dG dG_*}{(G-G_*)^2}, \qquad
      \omega = -\varDelta^{-2}\,dG,\quad
      \theta = \frac{\varDelta^2 dG_*}{(G-G_*)^2},
  \label{eq:Q-Gauss}\\
  &\frac{\hat\omega'}{\hat\omega}=\frac{G''}{G'}-2\frac{G'}{G-G_*}, \quad
  \frac{\hat\theta'}{\hat\theta}=\frac{G''_*}{G'_*}-2\frac{G'_*}{G_*-G},
  \label{eq:can-der-Gauss}
  \\
  &s(\hat\omega)=2\hat Q+\{G,z\}, 
  \quad s(\hat\theta)=2\hat Q+\{G_*,z\},
  \label{eq:schwarz-Gauss}
 \end{align}
 where ${}':=d/dz$, $Q=\hat Q\,dz^2$, $s(\cdot)$ is as in
 \eqref{eq:schwarz} and 
 $\{G,z\}$ represents
 the Schwarzian derivative of $G$ with respect to $z$:
 \begin{equation}\label{eq:schwarz2}
    \{G,z\} = \left(
		     \frac{G''}{G'}
		    \right)'-
    \frac{1}{2}
    \left(
     \frac{G''}{G'}
    \right)^2.
 \end{equation}
\end{remark}

\subsection{Generic singularities of flat fronts}
As an application of Theorem~\ref{thm:criterion}, we shall now show that
generic singularities of flat fronts are cuspidal edges or swallowtails.
Let $U$ be a simply-connected domain in $\C$, and 
$\O(U)$ the set of holomorphic functions on $U$.
Then, for each $h\in \O(U)$, we can construct a flat front
\[
  f_h \colon U\longrightarrow H^3
\]
which is represented by a pair of holomorphic $1$-forms
$(\omega, \theta)=(dz,e^h\, dz)$.

Conversely, if $p$ is not an umbilic point
(i.e. $Q(p)=\omega(p)\theta(p)\neq 0$), both $\omega(p)$ and
$\theta(p)$ are not equal to zero, and 
we can choose a complex coordinate such that $\omega=dz$ and
$\theta=e^h\,dz$.
Thus, any flat front is locally congruent to some $f_h$, except in 
neighborhoods of umbilic points.
We remark that an umbilic point cannot be a singular point, since
$ds^2_{1,1}=|\omega|^2+|\theta|^2$ is positive definite.
By Theorem~\ref{thm:criterion}, we have the following 
\begin{enumerate}
 \item The zeroes of $\Re h$ correspond to singular sets.
       Moreover, a singular point $p\in U$ is non-degenerate
       if and only if $h'(p)\neq 0$.
 \item A singular point $p$ is a cuspidal edge if and only if
       it is non-degenerate and $e^{-h(p)/2}h'(p)\not\in\R$.
 \item A singular point $p$ is a swallowtail if and only if
       it is non-degenerate, $e^{-h(p)/2}h'(p)\in\R$
       and 
       $\Re\left[e^{-h(p)}\left(h''(p)
       -\frac{1}{2}h'(p)^2\right)\right]\neq 0$.
\end{enumerate}
We let $J_H^k(U)$ be the space of $k$-jets of holomorphic functions on
$U$. 
Then $J_H^k(U)$ is canonically identified with the product space 
$U \times \C^{k+1}$: 
\[
   J_H^k(U)\ni j^kh \longleftrightarrow 
      (p,h(p),h'(p),h''(p),\dots,h^{(k)}(p))
      \in U \times \C^{k+1} \; . 
\]
In particular, $J_H^k(U)$ can be considered as
a $C^\infty$-manifold of dimension $2(k+2)$
as well as a complex manifold of dimension $k+2$.
For a compact set $K$ of $U$ and an open subset $O$
in $J_H^k(U)$, we set
\[
   [K,O]_k:=\{h\in \mathcal O(U)\,;\,
              {j^kh}(K)\subset O\}.
\]
Let $\fO_k$ be the topology of $\O(U)$ generated by such $[K,O]_k$,
which is called {\it the compact open $C^k$-topology.}
If 
$\pi:J^{k+1}_H(U)\to J^k_H(U)$
is the canonical projection, it can be easily seen that
$\pi$ is a continuous map and
$[K,\pi^{-1}(O)]_{k+1}=[K,O]_k$.
In particular, 
$\fO_k\subset \fO_{k+1}$ holds,
and 
\[
   \fO:=\bigcup_{k=0}^\infty \fO_k
\]
gives a topology on $\O(U)$, 
called the {\it compact open $C^\infty$-topology}. 
A holomorphic function $h\in \O(U)$ is an interior point of a given 
subset $\mathcal{V}(\subset \O(U))$
 (with respect to the compact open $C^\infty$-topology)  
if and only if there exist a non-negative integer $l$ and a finite
number of compact sets $C_1,\dots,C_s$ in $U$ and open subsets
$O_1,\dots,O_s \subset J^l_H(U)$ such that
\[
   h \in \bigcap_{r=1}^s [C_r,O_r]_l \subset \mathcal{V}\;.
\]
Now we give a topology on the family of flat fronts
$\{f_{h}\}_{h\in \O(U)}$ induced from the compact open 
$C^\infty$-topology on $\O(U)$.
We shall prove the following: 
\begin{theorem}\label{thm:genericity}
 Let $K$ be an arbitrary compact set of $U$ and
 $S(K)$ the subset of $\{f_h\}_{h\in \mathcal O(U)}$
 which consists of $f_h$ whose singular points on $K$ are locally
 diffeomorphic to cuspidal edges or swallowtails.
 Then $S(K)$ is an open and dense subset of $\{f_h\}_{h \in \O(U)}$.
\end{theorem}

\begin{remark}
 Generic properties of $C^\infty$-maps are usually described in
 terms of the Whitney $C^\infty$-topology (cf.\ \cite{GG}), 
 because it is suitable for the technique of multiplying by a
 cut-off function.
 However, generic properties of analytic functions are different 
 in the Whitney $C^{\infty}$-topology.
 In the above theorem, we use the compact open $C^\infty$-topology.
 The two topologies are the same when the source space is compact.
 However, they are different on $\O(U)$.
 In fact, when the source space is non-compact, the 
 compact open $C^\infty$-topology satisfies the second axiom of
 countability,
 but the Whitney $C^\infty$-topology
 on $\O(U)$ does not satisfy even the first axiom of countability
 and cannot be treated by sequence convergence.
 We do not know if the set $S(U)$ 
 (which consists of $f_h$ whose singular points on $U$ are locally
 diffeomorphic to cuspidal edges or swallowtails)
 is an open dense subset with respect to the Whitney $C^\infty$-topology.
\end{remark}
\begin{proof}[Proof of Theorem~\ref{thm:genericity}]
 We set
\begin{align*}
   A_1&:=\{j^2h(p)\in J^2_H(U)\,|\,
        \Re h(p)=0\text{ and }h'(p)=0 \}, \\
   A_2&:=
   \left\{
    j^2 h(p)\in J^2_H(U)\,\left|
   \begin{array}{l}
    \Re h(p)=0, \quad\Im \bigl(e^{-h(p)/2}h'(p)\bigr)=0,\\[2pt]
    \Re\left[e^{-h(p)}\left(h''(p)-\frac{1}{2}h'(p)^2\right)\right]=0
   \end{array}
   \right.
   \right\}.
\end{align*}
Then $A_1$ and $A_2$ are both closed subsets of $J^2_H(U)$.
The set given by
\[
  \hat S(K)
    =\{h\in \O(U)\,;\, {j^2h}(K)\subset J^2_H(U)\setminus(A_1\cup A_2)\}
    =[K,J^2_H(U)\setminus(A_1\cup A_2)]_2
\]
corresponds to $S(K)$ under the identification $h\leftrightarrow f_h$,
which is by definition an open subset in $\O(U)$.

So it is sufficient to show that $\hat S(K)$ is a dense subset.
Obviously, $A_1$ is a real closed submanifold of $J^2_H(U)$ with
codimension three.
We remark that $J^2_H(U)\setminus A_1$ is
an open submanifold of $J^2_H(U)$.
The following lemma holds: 
\renewcommand{\qed}{\relax}
\end{proof}
\begin{lemma}\label{submfd}
  $A_2\setminus A_1$ is a submanifold of  $J^2_H(U)$
  with codimension three.
\end{lemma}
\begin{proof}
  We define a $C^\infty$-map
  $\zeta=(\zeta^1,\zeta^2,\zeta^3) 
   \colon J^2_H(U)\setminus A_1\rightarrow\R^3$ by
\begin{align*}
    \zeta(j^2 h(p))&=\left(
      \Re h(p),
      \Im \bigl(e^{-h(p)/2}h'(p)\bigr),
      \Re\left[e^{-h(p)}\left(h''(p)-
           \frac{1}{2}\bigl(h'(p)\bigr)^2\right)\right]
      \right)\\
    &=\left(u, e^{-u/2}(v_1\cos\frac{v}2-u_1\sin\frac{v}2), \right. \\
   &\phantom{=(u, e^{-u/2}(v_1\cos} \left.
    \frac{e^{-u}}2\left((-u_1^2+v_1^2+2u_2)\cos v+2(-u_1v_1+v_2)\sin v
    \right)
 \right),
\end{align*}
where we set
\[
       h(z)=u(z)+\sqrt{-1} v(z),~ 
       h'(z)=u_1(z)+\sqrt{-1}v_1(z),~
       h''(z)=u_2(z)+\sqrt{-1}v_2(z).
\]
Then $(z,u,v,u_1,v_1,u_2,v_2)$ gives the canonical coordinate
system on $J^2_H(U)$. 
By a direct calculation, we have $\zeta^{-1}(0,0,0)=A_2\setminus A_1$.
We show that $(0,0,0)$ is a regular value of $\zeta$.
To determine the rank of the Jacobian matrix of $\zeta$
at any point in $\zeta^{-1}(0,0,0)$,
we calculate the derivative of $\zeta$ with respect to $u, u_1$ and $v_1$:
\begin{align*}
 (\zeta^1_u,\zeta^2_u,\zeta^3_u)
 &=\left(1, \frac{e^{-u/2}}2(-v_1\cos\frac{v}2+u_1\sin\frac{v}2),\right. \\
 &\phantom{\frac12e^{-u/2}(\cos\frac{v}2}\left.
-\frac{e^{-u}}2\left((-u_1^2+v_1^2+2u_2)\cos v+2(-u_1v_1+v_2)\sin v
 \right)\right), \\
 (\zeta^1_{u_1},\zeta^2_{u_1},\zeta^3_{u_1})
 &=\left(0, -e^{-u/2}\sin\frac{v}2,
 -{e^{-u}}\left(u_1\cos v+v_1\sin v\right)\right), \\
 (\zeta^1_{v_1},\zeta^2_{v_1},\zeta^3_{v_1})
 &=\left(0, e^{-u/2}\cos\frac{v}2,
 {e^{-u}}\left(v_1\cos v-u_1\sin v\right)\right).
\end{align*}
Then we have
\[
 \frac{\partial(\zeta^1,\zeta^2,\zeta^3)}{\partial (u,u_1,v_1)}
    =e^{-3u/2}(u_1\cos \frac{v}2+v_1\sin \frac{v}2) \; .
\]
We now suppose $\zeta(z,u,v,u_1,v_1,u_2,v_2)=0$.
Then
\[
   u=0,\quad
   e^{-u/2}(v_1\cos\frac{v}2-u_1\sin\frac{v}2)=0
\]
hold and thus
\begin{equation}\label{eq:star}
     v_1\cos\frac{v}2-u_1\sin\frac{v}2=0.
\end{equation}
Then \eqref{eq:star} and 
${\partial(\zeta^1,\zeta^2,\zeta^3)}/{\partial (u,u_1,v_1)}=0$
imply that $u=u_1=v_1=0$, namely, that 
$(z,u,v,u_1,v_1,u_2,v_2)$ belongs to $A_1$.
Hence $d\zeta$ is of rank $3$ at $\zeta^{-1}(0,0,0)$ in
 $J^2_H(U)\setminus A_1$.
By the implicit function theorem,
$A_2\setminus A_1$ is a submanifold of codimension $3$.
\end{proof}
\begin{proof}[Proof of Theorem~\ref{thm:genericity} {\rm(}continued\/{\rm )}]
We shall prove that $\hat S(K)$ is a dense subset.
We now fix a function $h\in \O(U)$.
Let $B$ be the set of polynomials of degree at most $2$ in $z$
and define a map as follows
\[
   G\colon{}U\times B\ni (z,\phi)\longmapsto j^2(h+\phi)(z)\in J^2_H(U).
\]
Obviously the map $G$ is a diffeomorphism. 
In particular, 
$G^{-1}(A_1)$ and $G^{-1}(A_2\setminus A_1)$
are both submanifolds of dimension $5$ diffeomorphic to
$A_1$ and $A_2\setminus A_1$ respectively.
Let
\[
  \pi\colon{}U\times B\longrightarrow B
\]
be the canonical projection.
Since $B$ is a $C^\infty$-manifold of dimension $6$,
Sard's theorem yields that
$\pi(G^{-1}(A_1))$ and $\pi(G^{-1}(A_2\setminus A_1))$
are measure zero sets in $B$.
Thus 
\[
  \pi(G^{-1}(A_1\cup A_2))
  =\pi(G^{-1}\bigl(A_1\cup (A_2\setminus A_1))\bigr)
  =\pi(G^{-1}(A_1))\cup \pi(G^{-1}(A_2\setminus A_1))
\]
is also a measure zero set.
Thus there is a sequence $\{\phi_n\}$ in $B$ such that
$\phi_n$ converges to the zero polynomial and
$\phi_n\not \in \pi(G^{-1}(A_1\cup A_2)).$
We set $h_n:=h+\phi_n$.
Then
$(j^2h_n)(U)\not\in A_1\cup A_2$,
that is, 
\[
  h_n\in [U,J^2_H(U)\setminus(A_1\cup A_2)]_2
      \subset \hat S(K).
\]

Let $d$ be a distance function on $J^{l}_H(U)$ which is
compatible with respect to its topology.
Then a sequence $\{g_n\}$ in $C^0\bigl(U,J^{l}_H(U)\bigr)$
{\it converges to $g$ uniformly}
on a given compact subset $K$ of $U$
if for any $\varepsilon0$, there exists a positive integer
$n_0$ such that
\[
   \sup_{z\in K }d\bigl(g_n(z),g(z)\bigr)<\varepsilon 
    \qquad (n\ge n_0)
\]
holds.  

We remark that
\[
   j^k\bigl(\O(U)\bigr)\subset C^0\bigl(U,J^{l}_H(U)\bigr)
\]
holds. 
Since the difference  $h_n-h$ is only a polynomial $\phi_n$ of
degree at most $2$ converging to the zero polynomial,
one can easily check that
for each non-negative integer $l$,
$j^l h_n$ converges to $j^l h$ uniformly on any
compact subset of $U$.

Let $\mathcal{V}$ be an open neighborhood of $h$ in $\O(U)$. 
Then by the definition of the compact open $C^\infty$-topology,  
there exist a non-negative integer $l$,
 a finite number of compact sets $C_1,\dots,C_s$
of $U$,
and open subsets $O_1,\dots,O_s$ of $J^l_H(U)$ such that
\[
   h\in \bigcap_{r=1}^s[C_r,O_r]_{l}\subset \mathcal{V}.
\]
We set
\[
   \Delta_r=d\bigl(j^l h(C_r),J^{l}_H(U)\setminus O_r\bigr)0
         \qquad (r=1,2,3,\dots, s).
\]
Note that
\[
   C:=C_1\cup C_2\cup\cdots \cup C_s
\]
is a compact set.
Since $j^l h_n$ converges to $j^l h$ uniformly on any compact subset of
$U$,
there exists an integer $n_00$ such that
\[
  \sup_{z\in C}d \left((j^{l} h)(z),(j^{l} h_n)(z)\right)
 <\frac{\min(\Delta_1,\dots,\Delta_s)}2
   \qquad (n\ge n_0).
\]
On the other hand,
\[
   d\left((j^{l} h)(z),J^{l}_H(U)\setminus O_r\right)
      \le d\left((j^{l} h)(z),(j^{l} h_n)(z)\right)+
      d\left(j^{l}(h_n)(z),J^{l}_H(U)\setminus O_r\right),
\]
then
\begin{align*}
   d\left(j^{l}(h_n)(z),J^{l}_H(U)\setminus O_r\right)
    &\ge d\left(j^{l} h(z),J^{l}_H(U)\setminus O_r\right)
          -d\left((j^{l} h)(x),(j^{l}h_n)(z)\right) \\
    &\ge \Delta_r-\frac{\min(\Delta_1,\dots,\Delta_s)}20
    \quad (z\in C_r).
\end{align*}
This implies that $j^l h_n(z)\in O_r$ if $z\in C_r$.
Thus $h_n \in [C_r,O_r]_{l}$ 
holds for all $r=1,2,\dots,s$ and
\[
    h_n\in  \bigcap_{r=1}^s[C_r,O_r]_{l}\subset \mathcal V
     \qquad (n\ge n_0) \; . 
\]
Since $h_n\in \hat S(K)$,
this implies that  $\hat S(K)$ is a dense subset.
\end{proof}
%
%
\section{Global Properties of Singular Points}
\label{sec:global}
In this section, we shall give a proof of Theorem~\ref{thm:global}
in the introduction.
\subsection{Preliminaries}
\label{subsec:global-prelim}
Let $f\colon{}M^2\to H^3$ be a flat front defined on a Riemann surface
$M^2$.
In this section, we do not assume that $M^2$ is 
simply-connected.
Thus the holomorphic lift $E_f$ of $f$ is defined only on the universal
cover $\widetilde M^2$ of $M^2$:
\[
    E_f\colon{}\widetilde M^2\longrightarrow\SL(2,\C),
\]
and then the canonical forms $\omega$ and $\theta$ as in
\eqref{eq:lift-can} are holomorphic $1$-forms defined on
$\widetilde M^2$.
Note that the first and second fundamental forms as in
\eqref{eq:fund-forms}, the Hopf differential as in \eqref{eq:hopf},
and the $(1,1)$-part $ds^2_{1,1}$ of the first fundamental form 
are all well-defined on $M^2$.
Moreover $ds^2_{1,1}$ is positive definite on $M^2$.
We have that 
\begin{align}
&\text{$\omega$ and $\theta$ have no common zeroes on $\widetilde M^2$, and}
\label{eq:can-nondeg}\\
&\text{$|\omega|^2$ and $|\theta|^2$ are well-defined pseudometrics
 on $M^2$.}
\label{eq:can-welldef}
\end{align}
From now on, we assume $f$ is {\em complete\/}, 
that is, there exist a
compact set $C\subset M^2$ and a symmetric $2$-tensor $T$ such that
$T$ is identically zero outside $C$ and $ds^2+T$ is a complete
Riemannian metric (see \cite{KUY2}).
We remark that $f$ is complete if and only if (see \cite{KRUY})
\begin{enumerate}
  \item\label{item:complete-1}
       The $(1,1)$-part $ds^2_{1,1}$ of the first fundamental form
       is complete
       (in this case, we say that $f$ is ``weakly-complete''),
  \item\label{item:complete-2}
       $ds^2_{1,1}$ has finite total absolute curvature, and
  \item\label{item:complete-3}
       the singular set is a compact set of $M^2$.
\end{enumerate}
In the proof of Theorem~\ref{thm:global}, we use only 
properties  \ref{item:complete-1} and \ref{item:complete-2};
that is, the conclusion of Theorem~\ref{thm:global} holds
for weakly-complete flat fronts such that $ds^2_{1,1}$ has 
finite absolute total curvature.

By completeness, we know that there exist a compact Riemann 
surface $\overline M^2$ and a finite number of points
$\{p_1,\dots,p_N\}$ in $\overline M^2$ such that 
\[
    M^2 \cong \overline M^2\setminus\{p_1,\dots,p_N\}
    \qquad\text{(i.e.\ they are biholomorphic)}
\]
(\cite[Lemma 3.3]{KUY2}, see also \cite{GMM}).
We call the points $\{p_j\}$ the {\em ends\/} of $f$.
Moreover, as shown in \cite[Lemma 2]{GMM}, the Hopf differential $Q$ can
be extended meromorphically on $\overline M^2$, and at each end $p_j$ 
there exists a complex coordinate $z$ around $p_j$ such that 
$z(p_j)=0$ and the canonical forms are written as
\begin{equation}\label{eq:can-order}
  \omega = \hat\omega(z)\,dz = z^{\mu} \hat\omega_0(z)\,dz,\qquad
  \theta = \hat\theta(z)\,dz = z^{\mu_*}\hat\theta_0(z)\,dz,
  \qquad
  (\mu,\mu_*\in\R),
\end{equation}
where $\hat\omega_0$ and $\hat\theta_0$ are holomorphic functions in $z$
which do not vanish at the origin.
Since $\mu$ and $\mu_*$ do not depend on the choice of complex
coordinates, we denote
\begin{equation*}
 \ord_{p_j}\omega:=\mu,\qquad
 \ord_{p_j}\theta:=\mu_*.
\end{equation*}
These are the orders of the pseudometrics $|\omega|^2$ and $|\theta|^2$,
respectively.
By \eqref{eq:hopf}, we have
\begin{equation*}
 \mu + \mu_* = \ord_{p_j}\omega + \ord_{p_j}\theta = \ord_{p_j}Q\in \Z, 
\end{equation*}
 where, by convention, $\ord_0 Q=k$ if $Q = z^k \,dz^2$.
Since $f$ is complete, $ds^2_{1,1}$ is a complete Riemannian metric on
$M^2$ \cite[Corollary 3.4]{KUY2}. 
Thus, we have
\begin{equation}\label{eq:end-order}
 \min\{\ord_{p_j}\omega , \ord_{p_j}\theta\} \leq -1.
\end{equation}
\begin{definition}\label{def:cylindrical}
 An end $p_j$ is called {\em cylindrical\/} if 
\[
    \ord_{p_j}\omega = \ord_{p_j}\theta.
\]
\end{definition}

Let $G$ and $G_*$  be the hyperbolic Gauss maps of $f$.
Then $G$ and $G_*$ are both meromorphic functions on $M^2$, and
$G(p)\neq G_*(p)$ for all $p\in M^2$.
\begin{fact}[{\cite[Theorem 4]{GMM}}]\label{fact:regular}
 At an end $p_j$, the following properties are equivalent{\rm :}
 \begin{enumerate}
  \item $G$ is meromorphic at $p_j$.
  \item $G_*$ is meromorphic at $p_j$.
  \item $\ord_{p_j}Q\geq -2$, that is,  $Q$ has at most a pole 
        of order $2$ at $p_j$.
 \end{enumerate}
\end{fact}
\begin{definition}\label{def:regular}
 An end $p_j$ is called {\em regular\/} if the three properties in
 Fact~\ref{fact:regular} hold.
 Otherwise, $p_j$ is called {\em irregular}.
\end{definition}
\begin{remark}\label{rem:cylindrical}
 The ends of the hyperbolic cylinders are regular and cylindrical.
 As a special case of  \cite[Theorem 6]{GMM}, a regular cylindrical end
 is asymptotic to a finite cover to a hyperbolic cylinder.
\end{remark}

An {\em umbilic point\/} $q\in M^2$ is a zero of the Hopf differential
$Q$.
When $Q$ is identically zero, that is, $f$ is totally umbilic, 
$f$ represents the horosphere.
In this section, we assume that $f$ is not totally umbilic.
Since $Q$ is meromorphic on the compact Riemann surface $\overline M^2$,
the number of umbilic points is finite.
As $ds^2_{1,1}$ is positive definite at $q$, \eqref{eq:hopf}
implies that either
\begin{multline}\label{eq:can-umbilic}
\bigl(
 \ord_q \omega=\ord_q Q\in \Z_+ \quad\text{and}\quad
  \ord_q\theta = 0
\bigr)\\
 \text{or}\qquad
\bigl(
 \ord_q \theta=\ord_q Q\in \Z_+\quad\text{and}\quad
  \ord_q\omega = 0
\bigr)
\end{multline}
holds at each umbilic point $q$.

Using a local complex coordinate $z$, we write
\begin{equation}\label{eq:can-local}
 \omega = \hat\omega\,dz,\qquad
 \theta = \hat\theta\,dz,\qquad
 Q = \hat Q\,dz^2.
\end{equation}
\subsection{%
  Global descriptions of the criteria for singular points}
\label{sub:criteria}

Let $\overline M^2$ be a compact Riemann surface and 
\[
   f\colon{}M^2=\overline {M}^2\setminus\{p_1,\dots,p_N\}\longrightarrow
   H^3
\]
a complete flat front which is not totally umbilic.
Using the canonical forms $\omega$ and $\theta$ in \eqref{eq:lift-can},
we define
\begin{equation*}
 \rho:=\frac{\theta}{\omega}.
\end{equation*}
Though $\rho$ might be defined only on the universal cover of $M^2$,
\eqref{eq:can-welldef} implies that $|\rho|$ is well-defined on $M^2$.
Moreover, by \eqref{eq:can-order}, $|\rho|$ can be extended on
$\overline M^2$ as a continuous map
\[
    |\rho|\colon{}\overline M^2 \longrightarrow [ 0, +\infty].
\]
As seen in Section~\ref{sec:local}, the set of singular points of the
flat front $f$ is given by
\begin{equation*}
 \Sigma(f):=\{p\in M^2\,;\,|\rho(p)|=1\}.
\end{equation*}
Using a local expression as in \eqref{eq:can-local}, we define
\begin{equation}\label{eq:criterion}
 \xi:=
  \left(\frac{\hat\theta'}{\hat\theta}-\frac{\hat\omega'}{\hat\omega}\right)
  \hat Q \,dz^3,\quad
 \zeta_c:=
 \left(\frac{\hat\theta'}{\hat\theta}-\frac{\hat\omega'}{\hat\omega}\right)^2
 \frac{1}{\hat Q},\quad\text{and}\quad
 \zeta_s:=
 \frac{s(\hat\theta)-s(\hat\omega)}{\hat Q},
\end{equation}
where $Q=\hat Q\,dz^2$, $'=d/dz$ and $s(\cdot)$ is as in
\eqref{eq:schwarz}.

\begin{lemma}\label{lem:criterion-global}
 The quantities in \eqref{eq:criterion} are independent on the choice 
 of complex coordinate.
 In particular,
 $\xi$ is a meromorphic $3$-form on $\overline M^2$, and 
 both $\zeta_c$ and $\zeta_s$ are meromorphic functions on 
 $\overline M^2$.
\end{lemma}
\begin{proof}
 Since $|\rho|$ is well-defined on $M^2$,
\[
 d(\rho\bar\rho)=|\rho|^2\left(\frac{d\rho}{\rho}+
                               \frac{d\bar \rho}{\bar\rho}\right)
\]
 is well-defined on $M^2$, and then so is its $(1,0)$-part.
 Hence
 \[
     \frac{d\rho}{\rho}=
           \left(\frac{\hat\theta'}{\hat\theta}-
                 \frac{\hat\omega'}{\hat\omega}
           \right)\,dz
 \]
 is a meromorphic $1$-form on $M^2$. 
 Moreover, by \eqref{eq:can-order} and \eqref{eq:can-umbilic}, 
 $d\rho/\rho$ is a meromorphic $1$-form on $\overline M^2$.
 Since $Q$ is a meromorphic $2$-differential on $\overline M^2$,
 $\xi=(d\rho/\rho)\cdot Q$ is a meromorphic $3$-differential.
 As the symmetric product $(d\rho/\rho)\cdot(d\rho/\rho)$ is a
 meromorphic $2$-differential, 
 \[
     \zeta_c = \frac{(d\rho/\rho)\cdot(d\rho/\rho)}{Q}
 \]
 is a meromorphic function on $\overline M^2$.
 
 Though the Schwarzian derivative as in \eqref{eq:schwarz2}
 depends on the choice of complex coordinates, the difference of
 two Schwarzian derivatives is  considered as a meromorphic
 $2$-differential; that is, if we write
 $S(G):=\{G,z\}\,dz^2$ in the complex coordinate $z$,
 \[
     S(G_*)-S(G)=\left[
                  \{G_*,z\}-\{G,z\}
                 \right]\,dz^2
 \]
 is independent of the choice of a coordinate $z$, as
 a meromorphic $2$-differential.
 Here, by \eqref{eq:schwarz-Gauss}, 
 \begin{equation}\label{eq:zeta-s-G}
      \zeta_s=\frac{\{G_*,z\}-\{G,z\}}{\hat Q}
             =\frac{S(G_*)-S(G)}{Q}
 \end{equation}
 holds.
 This shows that $\zeta_s$ is a well-defined meromorphic function
 on $M^2$.
 Moreover, by \eqref{eq:can-order} and the definition
 \eqref{eq:criterion}, $\zeta_s$ is meromorphic at each end.
\end{proof}

Using the invariants of \eqref{eq:criterion}, we define
\begin{equation*}
 \begin{array}{r@{\,}l}
  \Sigma(f)&:=\{p\in M^2\,;\,|\rho(p)|=1\},\\[6pt]
  Z_0(f)&:=\{p\in M^2\,;\,\xi(p)=0\},\\[6pt]
  Z_c(f)&:=\{p\in M^2\,;\,\Im\sqrt{\zeta_c(p)}=0\},\\[6pt]
  Z_s(f)&:=\{p\in M^2\,;\,\Re\zeta_s(p)=0\}.
 \end{array}
\end{equation*}
Though $\sqrt{\zeta_c}$ is multi-valued on $M^2$, the condition
$\Im\sqrt{\zeta_c}=0$ is unambiguous.

Then by Theorem~\ref{thm:criterion}, the following hold:
\begin{itemize}
  \item $p\in M^2$ is a singular point if and only if $p\in \Sigma(f)$.
  \item $p\in\Sigma(f)$ is a non-degenerate singular point if and only 
	if $p\in Z_0(f)^c$.
  \item $p\in\Sigma(f)$ is a cuspidal edge if and only if
	$p\in Z_0(f)^c\cap Z_c(f)^c$.
  \item $p\in\Sigma(f)$ is a swallowtail if and only if
	$p\in Z_0(f)^c\cap Z_c(f)\cap Z_s(f)^c$.
  \item A singular point $p\in\Sigma(f)$ is  neither
	a cuspidal edge nor a swallowtail if and only if 
	\[ 
              p\in Z_0 (f)\cup \bigl(Z_c(f)\cap Z_s(f)\bigr).
        \]
\end{itemize}
Here we denote the complementary set by the upper suffix ${^c}$.  
Since $Z_c(f)$ (resp.\ $Z_s(f)$) describes a criterion for a singular
point to be a cuspidal edge (resp.\ a swallowtail), we use the 
lower suffix ``$c$'' (resp.\ ``$s$'').

The sets $Z_0(f)$, $Z_c(f)$ and $Z_s(f)$ are the same for all the  
parallel fronts of $f$; that is, 
if $\{f_t\}_{t\in\R}$ is the family of parallel fronts of $f$, 
then we have:
\begin{lemma}\label{lem:parallel-criterion}
 \begin{align*}
    &\Sigma(f_t)=\{p\in M^2\,;\,|\rho(p)|=e^{2t}\},\\
    &Z_0(f_t)=Z_0(f),\quad
    Z_c(f_t)=Z_c(f),\quad
    Z_s(f_t)=Z_s(f).
 \end{align*}
\end{lemma}%
\begin{proof}
 By \eqref{eq:parallel-canonical}, we have
 the first assertion.
 Though the remaining parts can be proved by direct calculations, 
we give an alternative proof:
 Let $G$ and $G_*$ be the hyperbolic Gauss maps of $f$.
 Then by \eqref{eq:can-der-Gauss}, we have
 \begin{align*}
  \xi &= \left(
          \frac{G_*''}{G_*'}-
          \frac{G''}{G'}+2 \frac{G'+G_*'}{G-G_*}
          \right)\hat Q\,dz^3,\\
  \zeta_c &=
          \left(
          \frac{G_*''}{G_*'}-
          \frac{G''}{G'}+2 \frac{G'+G_*'}{G-G_*}
          \right)^2\frac{1}{\hat Q},
 \end{align*}
 and $\zeta_s$ is written as in \eqref{eq:zeta-s-G}.
 Since the hyperbolic Gauss maps and the Hopf differential are 
independent of the choice of parallel front $f_t$, we have the conclusion.
\end{proof}

By a direct calculation using the formulas in the proof of
Lemma~\ref{lem:parallel-criterion}
and \eqref{eq:Q-Gauss}, we have
\begin{equation}\label{eq:zeta-deriv}
 \zeta_s = \frac{\sqrt{\zeta_c}'}{\sqrt{\hat Q}}
  \qquad \left({~}'=\frac{d}{dz}\right).
\end{equation}
Using this, we can prove that:
\begin{proposition}\label{prop:revolution}
 Let $f$ be a complete flat front
 which is not totally umbilic.
 Then the function $\zeta_c$ is constant if and only if $f$ is a
 covering of a front of revolution.%
\end{proposition}
\begin{proof}
 If $\zeta_c=0$, $d\rho=0$ holds on $M^2$.
 Hence $\rho$ is constant.
 In this case, one can conclude that $f$ is  a covering of a
 hyperbolic cylinder,
 which is a surface of revolution. %
 
 On the other hand, assume $\zeta_c$ is a non-zero constant.
 By \eqref{eq:can-order},
 $\hat\theta'/\hat\theta-\hat\omega'/\hat\omega$
  can only have simple poles.
 Then by the definition of $\zeta_c$ in \eqref{eq:criterion},
 the order of $Q$ is at least $-2$.
 Thus, by Fact~\ref{fact:regular}, all ends must be regular.
 
 By \eqref{eq:zeta-deriv}, $\zeta_s=0$ holds.
 Then by \eqref{eq:zeta-s-G}, we have
 $\{G,z\}=\{G_*,z\}$ with respect to any complex coordinate $z$.
 Then it holds that
 \begin{equation}\label{eq:gauss-differ}
    G_*=b\star G=\frac{b_{11}G+b_{12}}{b_{21}G+b_{22}},\qquad
    b = \begin{pmatrix}b_{11} & b_{12} \\ b_{21} & b_{22}\end{pmatrix}
      \in\SL(2,\C),
 \end{equation}
 where $\star$ denotes the M\"obius transformation.
 Here, the group $\SL(2,\C)$ acts isometrically on $H^3$ as
 \begin{equation}\label{eq:isom}
      H^3\ni x \longmapsto axa^*\in H^3\qquad a\in\SL(2,\C),
 \end{equation}
 where we consider $H^3$ as in \eqref{eq:hyperbolic}.
 Under the isometry \eqref{eq:isom}, the hyperbolic Gauss maps transform
 as $(G,G_*)\mapsto (a\star G,a\star G_*)$.
 Hence we may assume $b$ in \eqref{eq:gauss-differ} is  a Jordan
 normal form.
 
 When $b$ is diagonal, we have $G_*=\mu G$, where $\mu$ is constant.
 Here, since $f$ is a flat front, $G$ and $G_*$ have no common branch
 points (see \cite{KUY2}).
 Thus $G$ has no branch point, and then we can take $z=G$ as a local
 coordinate.
 Hence 
 $f$ is locally congruent to a front of
 revolution (see Example~\ref{ex:revolution} in
 Section~\ref{sec:example}).
 Thus we have the conclusion.
 If $b$ is not diagonal, the eigenvalue of $b$ is $\pm 1$,
 which is a double root.
 Then we have $G_*=G-1$.
 Since the ends of $f$ are the points where $G=G_*$
 (\cite[Lemma~4.10]{KUY2}),
 the ends are common poles of $G$ and $G_*$.
 In this case, by \eqref{eq:Q-Gauss}
 we have $Q=-dG\,dG_*=-dG^2$.
 Then the $\ord_p Q$ at a pole $p$ of $G$ is less than or equal to $-4$,
 which contradicts the  fact that all ends are regular.
\end{proof}
\subsection{Proof of Theorem~\ref{thm:global}}
Let 
\begin{equation}\label{eq:flat-front-ass}
   f\colon{}M^2=\overline M^2\setminus\{p_1,\dots,p_N\}\longrightarrow
   H^3
\end{equation}
be a complete flat front which is not totally umbilic, and
$\{f_t\}$ its parallel family.
For simplicity, we write
\begin{multline*}
   \Sigma_t:=\Sigma(f_t),\quad
   Z_0:=Z_0(f_t)=Z_0(f),\quad
   Z_c:=Z_c(f_t)=Z_c(f),\\
  \text{and}\quad
   Z_s:=Z_s(f_t)=Z_s(f).
\end{multline*}

A point $p\in M^2$ is a singular point of $f_t$ which is neither a
cuspidal edge nor a swallowtail if and only if
\[
    p\in \Sigma_t\cap\bigl(Z_0\cup(Z_c\cap Z_s)\bigr).
\]
Then by Lemma~\ref{lem:parallel-criterion}, $f_t$ admits such a singular
point if and only if
\begin{equation}\label{eq:singular-singular}
    \{|\rho(p)|\,;\,p\in Z_0\cup(Z_c\cap Z_s)\} \ni e^{2t}.
\end{equation}%
Since $\xi$ in \eqref{eq:criterion} is a meromorphic
$3$-differential on the compact Riemann surface $\overline M^2$
and $Z_0$ is the set of zeroes of $\xi$, $Z_0$ is a finite set of
points.
Thus, to prove Theorem~\ref{thm:global}, it is sufficient to show 
the following proposition:
\begin{proposition}\label{prop:global}
 Let $f$ be a complete flat front such that $\zeta_c$ defined in 
 \eqref{eq:criterion} is not constant.
 Then 
 $\{|\rho(p)|\,;\,p\in Z_c\cap Z_s\}\subset\R_+$ is a finite set.
\end{proposition}
Before proving this proposition, we shall give a proof of 
Theorem~\ref{thm:global}.
\begin{proof}[Proof of Theorem~\ref{thm:global}]
 Assume a complete flat front $f$ is a front of revolution.
 Such a flat front is a horosphere, a finite cover of a hyperbolic
 cylinder, a snowman, or an hourglass (see Example~\ref{ex:revolution}
 in Section~\ref{sec:example}).
 Among these, the horospheres and hyperbolic cylinders do not have
 singular points, and all singularities of the snowman are cuspidal
 edges.
 Since we assumed $f$ is not a cover of an hourglass, 
 we have the conclusion for the case of fronts of revolution.

 Next, we assume $f$ is not a front of revolution.
 Then by Proposition~\ref{prop:revolution}, $\zeta_c$ is non-constant.
 Hence by Proposition~\ref{prop:global},
 $\{|\rho(p)|\,;\,p\in Z_c\cap Z_s\}$ is a finite set.
 On the other hand, the parallel front $f_t$ admits a singular point
 which is neither a cuspidal edge nor a swallowtail if and only if 
 \eqref{eq:singular-singular} holds.
 Hence we have the conclusion.
\end{proof}

To prove Proposition~\ref{prop:global}, we need the following lemma:
\begin{lemma}\label{lem:accumulate}
 Let $f$ be a complete flat front as in \eqref{eq:flat-front-ass}
 with non-constant $\zeta_c$.
 Assume $Z_c\cap Z_s$ accumulates at a point $p\in\overline M^2$.
 Then
\begin{enumerate}
 \item $p$ is a non-umbilic point in $M^2$ or an irregular cylindrical
       end, and
 \item there exists a neighborhood $U$ of $p$ such that the number
       of connected components of
       \[
           \bigl(U\setminus\{p\}\bigr)\cap(Z_c\cap Z_s)
       \]
       is finite, and each connected component is a level 
       set of $|\rho|$.
\end{enumerate}
\end{lemma}
This lemma will be proven in Section~\ref{sec:proof-of-lemma} later.
Using these, we shall prove Proposition~\ref{prop:global}.
\begin{proof}[Proof of Proposition~\ref{prop:global}]
 Assume
 \[
   \#\{|\rho(p)|\,;\,p\in Z_c\cap Z_s\}=+\infty.
 \]
 Then there exists an infinite sequence $\{z_n\}\subset Z_c\cap Z_s$
 such that $|\rho(z_n)|$ ($n=1,2,\dots$) are mutually distinct.
 Since $\overline M^2$ is compact, we can take a subsequence of
 $\{z_n\}$ which converges to $z_{\infty}\in \overline M^2$.
 Thus by Lemma~\ref{lem:accumulate}, $\#\{|\rho(z_n)|\}$
 is finite.
 This is a contradiction because the $|\rho(z_n)|$ $(n=1,2,\dots)$ are
 mutually distinct.
\end{proof}

\subsection{Proof of Lemma~\ref{lem:accumulate}}
\label{sec:proof-of-lemma}
\begin{proof}[Proof of the first part of Lemma~\ref{lem:accumulate}]
 Let $p$ be an accumulation point of $Z_c\cap Z_s$,
 and take a sequence $\{p_n\}$ consisting of mutually distinct
 points in $Z_c\cap Z_s$ such that $p_n\to p$ as $n\to\infty$.
 We show the first assertion of the lemma by way of contradiction:
 We assume
 \begin{itemize}
  \item $p\in M^2$ is an umbilic point, or
  \item $p\in \overline M^2$ is an end which is not an
	irregular cylindrical end,
 \end{itemize}
 and  set 
 \[
   \mu=\ord_p \omega ,\qquad 
   \mu_*=\ord_p\theta 
   \qquad\text{and}\qquad
   k=\ord_p Q =\mu+\mu_*\in\Z.
 \]
 If $p$ is an umbilic point, $\mu\neq \mu_*$ holds because of
 \eqref{eq:can-umbilic}.
 If $p$ is an end, $\mu=\mu_*$ holds when $p$ is cylindrical.
 So, we consider two cases:
\begin{itemize}
 \item[{\it Case 1}:] $\mu\neq\mu_*$,
	       that is, $p$ is an umbilic point or a 
	       non-cylindrical end.
 \item[{\it Case 2}:] $\mu=\mu_*$,
	       that is, $p$ is a cylindrical end.
	       In this case, $p$ is a regular end
	       because of our assumption.
               Then by Fact~\ref{fact:regular} and 
               \eqref{eq:end-order}, we have $\mu=\mu_*=-1$.
\end{itemize}
\paragraph{\it Case 1{\rm:}}
We assume $\mu\neq\mu_*$.
If we take a complex coordinate $z$ around $p$ such that $z(p)=0$,
we can write
\begin{equation}\label{eq:case-1-q}
    \frac{\hat\theta'}{\hat\theta}-
    \frac{\hat\omega'}{\hat\omega}
      = \frac{a}{z}\bigl(1+O(z)\bigr)\qquad
      \bigl(a:=\mu_*-\mu\bigr),
\end{equation}
where $O(z)$ denotes a higher-order term.
On the other hand, the Hopf differential $Q$ is written as
\[
    Q = z^k\,\bigl(q_0+O(z)\bigr)\,dz^2\qquad (q_0\neq 0).
\]
Thus, it follows from \eqref{eq:criterion} that
\begin{equation}\label{eq:case-1-zeta-c}
    \sqrt{\zeta_c} = z^{-(k+2)/2}\left(
                           \frac{a}{\sqrt{q_0}}+O(z)
                          \right).
\end{equation}

We assume $k\neq -2$.
Then by \eqref{eq:case-1-q}, \eqref{eq:case-1-zeta-c} and
\eqref{eq:zeta-deriv}, we have
\begin{equation}\label{eq:case-1-zeta-s-1}
    \zeta_s = z^{-k-2}\left(
                       \frac{a}{q_0}+O(z)
                      \right).
\end{equation}
Let $z_n=z(p_n)$.
Then $z_n$ tends to the origin as $n\to\infty$.
Since $p_n\in Z_c\cap Z_s$, $\Im\sqrt{\zeta_c(z_n)}=\Re\zeta_s(z_n)=0$
holds.
Since $a\in\R$, there exist sequences $\{\varepsilon_n\}$
and $\{\varepsilon_n'\}$ of real numbers such that
\begin{alignat}{3}
 0 &\equiv \arg\sqrt{\zeta_c(z_n)} &=&
            -\left(\frac{k}{2}+1\right)\arg z_n -
             \frac{1}{2}\arg q_0 + \varepsilon_n\qquad
      &&\pmod\pi,\label{eq:case-1-a-arg-1}\\
 \frac{\pi}{2} & \equiv \arg\zeta_s(z_n)
      &=&-(k+2)\arg z_n
    -\arg q_0 + \varepsilon'_n\qquad
      &&\pmod\pi,\label{eq:case-1-a-arg-2}
\end{alignat}
and $\varepsilon_n, \varepsilon'_n\to 0$ as $n\to\infty$.
Here, 
by \eqref{eq:case-1-a-arg-1} and \eqref{eq:case-1-a-arg-2},
we have
\[
   -\frac{\pi}{2}\equiv
   2\arg\sqrt{\zeta_c(z_n)}-\arg\zeta_s(z_n)\equiv
      2\varepsilon_n-\varepsilon_n'
      \qquad\pmod\pi,
\]
giving a contradiction.
Then the case $k\neq -2$ is impossible.

Assume $k=-2$.
In this case, \eqref{eq:case-1-zeta-c} is written as
$\sqrt{\zeta_c} = a q_0^{-1/2} + O(z)$.
Then by the assumption that $\zeta_c$ is non-constant, there 
exists a positive integer $l$ such that
\begin{equation}\label{eq:case-1-zeta-c-b}
    \sqrt{\zeta_c} = \frac{a}{\sqrt{q_0}}+b\,z^l + O(z^{l+1})
    \qquad (b\neq 0).
\end{equation}
In this case, by \eqref{eq:zeta-deriv}, we have
\begin{equation}\label{eq:case-1-zeta-s-b}
 \zeta_s = z^l \left(\frac{l\,b}{\sqrt{q_0}}+O(z)\right).
\end{equation}
Here $\Im\sqrt{\zeta_c(z_n)}=0$ holds on a sequence
$\{z_n=z(p_n)\}$
such that $z_n\to 0$ as $n\to\infty$, and $a\in\R$.
Hence \eqref{eq:case-1-zeta-c-b} implies that $\sqrt{q_0}\in\R$.
Thus, we have
\begin{alignat}{3}
  0 &\equiv \arg\sqrt{\zeta_c(z_n)} &=&
            \arg \left(b\,z_n^l+O(z^{l+1})\right)
            =l\arg z_n + \arg b + \varepsilon_n
      &&\pmod\pi,\label{eq:case-1-b-arg-1}\\
 \frac{\pi}{2} & \equiv \arg\zeta_s(z_n)
      &=&l \arg z_n + \arg b + \varepsilon_n',
      &&\pmod\pi,\label{eq:case-1-b-arg-2}
\end{alignat}
where $\varepsilon_n,\varepsilon_n'\to 0$ as $n\to\infty$.
Again, \eqref{eq:case-1-b-arg-1} and \eqref{eq:case-1-b-arg-2}
contradict each other.
\paragraph{\it Case 2{\rm:}}
We assume $\mu=\mu_*=-1$.
Taking a complex coordinate $z$ such that $z(p)=0$,
we can write $\hat Q = z^{-2}\bigl(q_0+O(z)\bigr)$.
Denote by $d\in \Z_+ \cup \{0\}$ the branch order of $G$ at $z=0$.
(for example, if $G=a+z^{d+1}$, the branch order of $G$ at $z=0$ is $d$.)
Since $\mu=-1$,  \eqref{eq:schwarz-Gauss} implies that 
\[
    q_0 = \frac{1}{4}(d+1)^20.
\]

Since $\mu=\mu_*=-1$, it holds that
$(\hat\theta'/\hat\theta)-(\hat\omega'/\hat\omega)=O(1)$.
Hence we have $\sqrt{\zeta_c(z)}=O(z)$.
Thus, we can write
\[
   \sqrt{\zeta_c} = b\,z^l+O(z^{l+1}),\qquad
  \hat Q = \frac{1}{z^2}\left(\frac{1}{4}(d+1)^2+O(z)\right),
\]
where $l\geq 1$ is an integer and $b\neq 0$.
Thus, 
\[
    \zeta_s=\frac{l\,b}{\sqrt{q_0}}z^l\bigl(1+O(z)\bigr)=
            \frac{l\,b}{d+1}z^l\bigl(1+O(z)\bigr). 
\]
As in Case 1, we set $z_n=z(p_n)$.
Then we have 
\[
    \arg b+l\arg z_n +\varepsilon_n\equiv 0 ,\qquad
    \arg b+l\arg z_n+\varepsilon_n'\equiv \frac{\pi}{2}\qquad\pmod\pi,
\]
where $\varepsilon_n,\varepsilon'_n\to 0$.
This is impossible.

Hence in any case, $Z_c\cap Z_s$ does not accumulate at $0$.
\end{proof}
\begin{proof}[Proof of the second part of Lemma~\ref{lem:accumulate}]
We consider two cases.
\paragraph{\it Case 1{\rm :}}
Suppose that $Z_c\cap Z_s$ accumulates at a non-umbilic point $p\in M^2$,
i.e.\ $Q(p)\neq 0$.

Take a complex coordinate $z$ around $p$ with $z(p)=0$.
Since $\hat Q(0)\neq 0$, there exists a holomorphic function
$\varphi(z)$ defined on a neighborhood of the origin such that
\[
    \bigl(\varphi(z)\bigr)^2 = \hat Q(z);
    \qquad\text{that is,}\quad
    \varphi(z)=\sqrt{\hat Q(z)}, \qquad \text{and}\qquad
    \varphi(0)\neq 0.
\]
On the other hand, both $\hat\theta$ and $\hat\omega$ have neither 
a zero nor a pole at $z=0$, so by \eqref{eq:criterion},
$\sqrt{\zeta_c}$ is a holomorphic function near $z=0$.
Since $\zeta_c$ is not a constant, 
there exists a positive integer $l$ such that
$\sqrt{\zeta_c}=a + b z^l + O(z^{l+1})$,
where $b\neq 0$.
Then by the Weierstrass preparation theorem, we can choose a
coordinate $z$ such that
\[
   \sqrt{\zeta_c(z)}=a+z^l,\qquad
   \varphi(z)=\sqrt{\hat Q(z)}=\varphi_0+O(z)
   \qquad\bigl(l\in\Z_+,\varphi_0\in\C\setminus\{0\}\bigr).
\]
Moreover, replacing $\varphi_0z$ by $z$, we can set
\begin{equation}\label{eq:accum-1}
   \sqrt{\zeta_c(z)}=a + b z^l,\qquad
   \varphi(z)=1+O(z)
   \qquad (l\in \Z_+,~b=\varphi_0^{-l}\in\C\setminus\{0\}).
\end{equation}
Here, since $Z_c$ accumulates at $0$, $a$ in \eqref{eq:accum-1} must be
real, and then
\begin{equation}\label{eq:accum-1-a}
 \Im\sqrt{\zeta_c} = \Im (b\,z^l)
\end{equation}
holds.
On the other hand, by \eqref{eq:zeta-deriv},
\begin{equation}\label{eq:accum-zeta-s}
   \zeta_s =\frac{l b z^{l-1}}{\varphi(z)}= 
    l b z^{l-1}\bigl(1+O(z)\bigr).
\end{equation}
We identify a neighborhood of $p$ with a neighborhood of the origin of
$z$-plane.
Since $Z_c\cap Z_s$ accumulates to the origin, we can take a 
sequence $\{z_n\}\subset Z_c\cap Z_s$ such that
$z_n\to 0$ as $n\to\infty$.
Then by \eqref{eq:accum-1-a} and \eqref{eq:accum-zeta-s},
there exists a sequence $\{\varepsilon_n\}\subset\R$
such that $\varepsilon_n\to 0$ and 
\begin{align}
 \arg b  + l \arg z_n &\equiv 0 \qquad\pmod\pi, \label{eq:accum-1-b}\\
 \arg b  + (l-1) \arg z_n +\varepsilon_n&\equiv \frac{\pi}{2}
 \qquad\pmod\pi\label{eq:accum-1-c}
\end{align}
hold.
Subtracting \eqref{eq:accum-1-c} from \eqref{eq:accum-1-b}, 
we have
\begin{equation}\label{eq:accum-1-d}
     \arg z_n-\varepsilon_n \equiv \frac{\pi}{2} \qquad\pmod\pi.
\end{equation}
On the other hand, subtracting \eqref{eq:accum-1-b} multiplied by $l-1$
from \eqref{eq:accum-1-c} multiplied by $l$, we have
\[
     \arg b \equiv l\frac{\pi}{2} - l\varepsilon_n \pmod\pi. 
\]
Here, since $\varepsilon_n\to 0$, we deduce that
\[
   \arg b \equiv l \frac{\pi} {2} \pmod \pi
   \qquad\text{and}\qquad
   \varepsilon_n=0.
\]
Substituting these into \eqref{eq:accum-1-d}, we have
$\arg z_n\equiv \pi/2 \pmod \pi$; that is, $z_n\in\sqrt{-1}\R$.
Since $\sqrt{\zeta_c(z_n)}\in\R$ holds for $n=1,2,\dots$, 
the imaginary part of $\sqrt{\zeta_c(z)}$ vanishes identically
on $\sqrt{-1}\R$, namely, 
\begin{equation}\label{eq:accum-1-e}
 \sqrt{\zeta_c(z)}\in\R\qquad
  \text{(if $z \in\sqrt{-1}\R$)}.
\end{equation}
Similarly, one can prove that 
\[
   \zeta_s(z)\in\sqrt{-1}\R \qquad
   \text{(if $z \in\sqrt{-1}\R$)}.
\]
Thus, on a neighborhood of the origin, $Z_c\cap Z_s$ is 
the imaginary axis in the $z$-plane.

Next, we shall prove that the imaginary axis is a level set of $|\rho|$.
By \eqref{eq:accum-1}, \eqref{eq:accum-zeta-s} and \eqref{eq:accum-1-e},
\begin{align*}
 \frac{\pi}{2} \equiv \arg \zeta_s(z)
              &= \arg \left(l b z^{l-1}\frac{1}{\varphi}\right)
               =\arg (bz^l)-\arg z\varphi(z)\\
              &= \arg{\sqrt{\zeta_c(z)}}-\arg z-\arg\varphi(z)\\
              &= -\arg z -\arg\varphi(z)=-\frac{\pi}{2}-\arg\varphi(z)
\pmod \pi
\end{align*}
holds on the imaginary axis.
Thus we have
\[
     \varphi(z)\in\R\qquad    \text{(if $z \in\sqrt{-1}\R$)}.
\]
As seen in \eqref{eq:singular-direction} in Section~\ref{sec:local},
the tangent vector field of a level set of $|\rho|$ is represented as
\[
     \sqrt{-1}\overline{
      \left(\frac{\hat\theta'}{\hat\theta}-
            \frac{\hat\omega'}{\hat\omega}\right)}.
\]
On the other hand, 
\[
    \sqrt{\zeta_c(z)}
    =\frac{1}{\varphi(z)}
      \left(\frac{\hat\theta'}{\hat\theta}-
            \frac{\hat\omega'}{\hat\omega}\right)\in\R,\quad
     \text{and}\quad
	    \varphi(z)\in\R\qquad (\text{if $z\in\sqrt{-1}\R$}).
\]
Without loss of generality, 
$d\rho\neq 0$ holds on $U\setminus\{0\}$, where $U$ is a neighborhood of
the origin, because a zero of $d\rho$ is isolated in $\overline M^2$.
Then the tangent vector of the level set of $|\rho|$ at a point on the
imaginary axis is parallel to the imaginary axis.
Hence the level set passing through a point of the imaginary axis
is the imaginary axis.
That is, $Z_c\cap Z_s$ coincides with the imaginary axis, which is 
a level set of $|\rho|$.

\paragraph{\it Case 2{\rm :}}
Suppose now that $Z_c\cap Z_s$ accumulates at an irregular
cylindrical end $p$.
Let $z$ be a complex coordinate with $z(p)=0$.
By irregularity, $\ord_p Q\leq -3$ holds.
Without loss of generality, we may assume $\ord_p Q$ is an even number.
In fact, if we set $z=w^2$, that is, we take the double cover of a
neighborhood of $p$, the order $Q$ at the origin with respect to the
coordinate $w$ will be an even number.

Hence, we assume
\begin{equation*}
     \ord_p Q = -2 k ,\qquad
     \ord_p \omega =\ord_p\theta=-k,
\end{equation*}
where $k\geq 2$ is an integer.
The second equality holds because $p$ is a cylindrical end.

Since $Q$ has even order at the origin, $\sqrt{\hat Q}$ 
is a meromorphic function on a neighborhood of $0$.
More precisely, we can write
\begin{equation*}
   \varphi:=\sqrt{\hat Q}=
    \frac{1}{z^k}\bigl(\varphi_0+O(z)\bigr)
    \qquad (\varphi_0\in\C\setminus\{0\}).
\end{equation*}
Since $\ord_p\omega=\ord_p\theta$, \eqref{eq:criterion} implies that
$\sqrt{\zeta_c(z)}=O(z^k)$, that is,
there exists an integer $l$ ($l\geq k$)
such that $\sqrt{\zeta_c(z)}=az^l+O(z^{l+1})$ $(a\in \C\setminus\{0\})$.
Then by the Weierstrass preparation theorem, we can choose a coordinate
$z$ such that
\begin{equation}\label{eq:accum-2-c}
 \sqrt{\zeta_c(z)} = z^l,\qquad
 \varphi=\frac{1}{z^k}\bigl(b+O(z)\bigr)
 \qquad(b\in\C\setminus\{0\}).
\end{equation}
Then by \eqref{eq:zeta-deriv}, $\zeta_s$ is written as
\begin{equation}\label{eq:accum-2-c1}
 \zeta_s(z)=l z^{l+k-1}\left(\frac{1}{b}+O(z)\right).
\end{equation}
As $\{z_n\}\subset Z_c\cap Z_s$ is a sequence with $z_n\to 0$,
we have
\begin{align}
 l \arg z_n &\equiv 0\qquad\pmod\pi,
 \label{eq:accum-2-d}\\
 -\arg b+
 (l+k-1) \arg z_n +\varepsilon_n&\equiv \frac{\pi}{2}\qquad\pmod\pi,
 \label{eq:accum-2-e}
\end{align}
where $\varepsilon_n\to 0$.
Subtracting  \eqref{eq:accum-2-e} from \eqref{eq:accum-2-d}, and
 \eqref{eq:accum-2-d} multiplied by $l+k-1$ from 
\eqref{eq:accum-2-e} multiplied by $l$, we have
\[
    -\arg b+(k-1)\arg z_n+\varepsilon_n\equiv \frac{\pi}{2},\qquad
    -l\arg b+l\varepsilon_n\equiv l\frac{\pi}{2}\qquad
    \pmod{\pi}.
\]
Since $\varepsilon_n\to 0$, this yields
$-l\arg b\equiv l\pi/2\pmod\pi$, and then 
$\varepsilon_n=0$ for sufficiently large $n$.
Thus, we have
\begin{equation}\label{eq:accum-2-f}
    (k-1)\arg z_n \equiv \frac{\pi}{2}+\arg b\qquad \pmod\pi.
\end{equation}

Let
\[
   \L_j:=\left\{
           z\in U\,;\, \arg z\equiv \frac{1+2j}{2(k-1)}\pi
	               +
		          \frac{\arg b}{k-1}
		          \pmod\pi
         \right\}
	 \qquad (j=0,\dots,k-1).
\]
Then  $\{\L_j\}$ is a set consisting of a finite number of lines 
in the $z$-plane through the origin, 
and by \eqref{eq:accum-2-f}, each $z_n$ lies on some
$\L_j$.
Hence there exists a subset $J\subset\{0,1,2,\dots,k-1\}$ such that
each $\L_j$ ($j\in J$) contains an infinite number of elements 
of $\{z_n\}$. We fix $j\in J$.
Then we can take a subsequence $\{z_m\}$ of $\{z_n\}$ such that 
$z_m\in \L_j\cap Z_c\cap Z_s$ and $z_m\to 0$.
 Since $\zeta_c(z_m)\in\R$ and $\zeta_s(z_m)\in\sqrt{-1}\R$, 
 we have
\begin{equation}\label{eq:accum-2-g}
 \sqrt{\zeta_c(z)} \in \R,\qquad
 \zeta_s  (z)\in\sqrt{-1}\R\qquad 
 \text{(if  $z\in\L_j$)}.
\end{equation}
This shows that, on a neighborhood of the origin,
$Z_c\cap Z_s$ coincides with the set of lines $\bigcup_{j\in J}\L_j$.

Next, we show that $\L_j$ ($j\in J$) is a level set of $|\rho|$
for each $j$.
By \eqref{eq:accum-2-c} and \eqref{eq:accum-2-g}, 
\begin{align*}
  \frac{\pi}{2}\equiv \arg\zeta_s(z)
  &=\arg\sqrt{\zeta_c(z)}-\arg z-\arg\varphi\\
  &=-\frac{1+2j}{2(k-1)}\pi
    -
      \frac{\arg b}{k-1}
   -\arg\varphi\pmod\pi
\end{align*}
holds on $\L_j$.
Hence we have
\begin{equation*}
    \arg\varphi\equiv \frac{\pi}{2}-\frac{1+2j}{2(k-1)}\pi
     -\frac{\arg b}{k-1}
 \pmod\pi 
 \qquad (\text{if $z\in \L_j$}).
\end{equation*}

At any point in $\L_j$, the argument of the tangent vector of the
level set of $|\rho|$ is
\begin{align*}
   \arg\sqrt{-1}\overline{\left(
                 \frac{\hat\theta'}{\hat\theta}-
		 \frac{\hat\omega'}{\hat\omega}
		 \right)}
    &\equiv
                \arg
                \left(\sqrt{-1}
                \overline{\left(\sqrt{\zeta_c}\varphi\right)}\right)
          \equiv
                 \frac{\pi}{2}-\arg\varphi-\arg\sqrt{\zeta_c}\\
         &\equiv \frac{\pi}{2}-\arg\varphi=
          \frac{1+2j}{2(k-1)}\pi+\frac{\arg b}{k-1}\pmod\pi,
\end{align*}
and then, the tangent vector is proportional to the line $\L_j$.
Hence each $\L_j$ ($j\in J$) is a level set of $|\rho|$.
Thus, we have the conclusion.
\end{proof}
\begin{remark}\label{rem:conical}
 Let $p\in \overline M^2$ be an accumulation point of $Z_c\cap Z_s$.
 Then by the second part of Lemma~\ref{lem:accumulate}, 
 $Z_c\cap Z_s$ is a level set of $|\rho|$ in a neighborhood of $p$;
 that is, by taking a suitable parallel front, we may assume that 
 a component of $Z_c\cap Z_s$ is a part of the singular set.
 Since the null direction and the singular direction coincide
 at each point in $Z_c\cap Z_s$, the image of such a singular set 
 is a single point in $H^3$.
 If the point $p$ is not an end, such a singularity seems to be a so-called
 {\em cone-like singularity},
 see, for example, the hourglass in Example~\ref{ex:revolution}.
 Another example is as follows:
 Set 
 \[
     \omega = \exp\left(z+\frac{1}{3}z^3\right)\,dz,\qquad
     \theta = \exp\left(-z-\frac{1}{3}z^3\right)\,dz
 \]
 on $\C$.
 Then by solving \eqref{eq:lift-eq}, we have a flat front 
 $f\colon{}\C\to H^3$.
 The singular set of $f$ contains the imaginary axis, which
 coincides with $Z_c\cap Z_s$ (see Figure~\ref{fig:conical}).
 However, this example is not complete because the canonical
 forms have an essential singularity at $z=\infty$.
\end{remark}
\begin{figure}
\begin{center}
  \includegraphics[width=2.5cm]{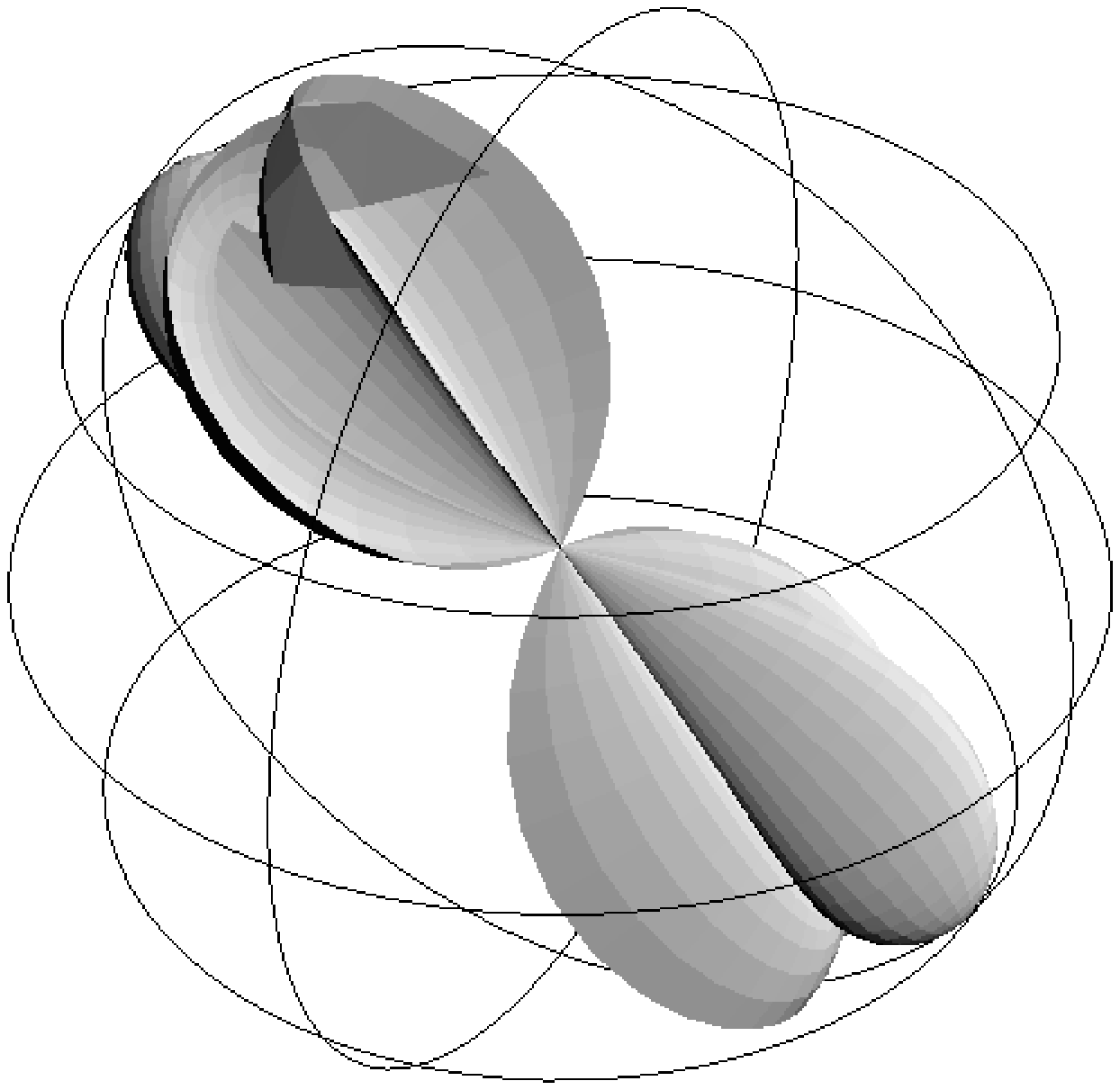} 
\end{center}
\caption{The example in Remark~\ref{rem:conical}}\label{fig:conical}
\end{figure}

\section{Caustics of flat fronts}
\label{sec:caustic}
In \cite{R}, Roitman investigated the caustic of (the parallel
family of) a flat front, which is considered as the locus of singular
points of the fronts in the parallel family.
In this section, we discuss caustics of flat fronts from our point of
view.

Let $U\subset \C$ be a simply connected domain and 
$f\colon{}U\to H^3$ a flat front without umbilic points.
We denote by $\omega$ and $\theta$ the canonical forms
of $f$, and $\rho:=\theta/\omega$, as in the previous section.
Since $f$ has no umbilic points, $\rho$ does not take values
$0$ and $\infty$.

For a point $z\in U$, we denote by $\kappa_1$ and
$\kappa_2$ the principal curvatures of $f$ at $z$.
Then by \eqref{eq:fund-forms}, we have
\begin{equation}\label{eq:principal}
    \kappa_1 = \frac{|\rho|+1}{|\rho|-1},\qquad
    \kappa_2 = \frac{|\rho|-1}{|\rho|+1}.
\end{equation}
Since $\rho\neq 0,\infty$, $|\kappa_1|>1$ holds.
Then there exists a real number $r_1$ such that 
$\coth r_1=\kappa_1$, which is called the {\em radius of curvature}.
By \eqref{eq:principal}, we have
\begin{equation}\label{eq:radius-rho}
    e^{2r_1}=|\rho|.
\end{equation}
The {\em caustic} $C_f$ of $f$ is defined as 
\begin{equation*}
 C_f\colon{}U \ni z \longmapsto
    \cosh r_1(z) f(z) + \sinh r_1(z)\nu(z)\in H^3 \subset {\Lor}^4,
\end{equation*}
where
$\Lor^4$ is the Minkowski $4$-space and
$\nu$ is the unit normal vector of $f$.
In other words, $C_f$ is the locus of the centers of the principal
curvature $\kappa_1$ of $f$.

Let $E_f\colon{}U\to \SL(2,\C)$ be the holomorphic lift of the front
$f$.
Then $f$ and the unit normal vector $\nu$ are given by 
\[
    f = E_fE_f^*,\qquad \nu = 
     E_f\begin{pmatrix} 1  & \hphantom{-}0 \\ 
                         0 & -1 \end{pmatrix}E_f^*.
\]
Thus, the caustic of $f$ is 
\begin{align*}
 C_f & = E_f\left[\cosh r_1 \begin{pmatrix} 1 & 0 \\ 0 & 1 \end{pmatrix}
               +\sinh r_1 \begin{pmatrix} 1 & \hphantom{-}0 \\ 0 &
			  -1\end{pmatrix}
           \right] E_f^*\\
     & =E_f\begin{pmatrix}e^{r_1} & 0 \\ 0 & e^{-r_1}\end{pmatrix}E_f^*
      = E_f\begin{pmatrix} |\rho|^{1/2} & 0 \\ 
                           0 & |\rho|^{-1/2}
           \end{pmatrix}E_f^*.
\end{align*}
Hence if we set
\begin{equation}\label{eq:caustic-lift}    
 E_c=
 	   E_f\begin{pmatrix}
	     \rho^{1/4} & 0 \\ 0 & \rho^{-1/4} 
            \end{pmatrix}
           P,\qquad
           P = \frac{1}{\sqrt{2}}
               \begin{pmatrix} 1 & \sqrt{-1}\\
		              \sqrt{-1} & 1 \end{pmatrix}\in\SU(2),
\end{equation}
we have
\begin{equation}\label{eq:caustic-legendre}
   C_f = E_c E_c^*,\qquad\text{and}\qquad
   E_{c}^{-1}d  E_c
   =\begin{pmatrix}
       0 & \theta_c \\
     \omega_c & 0
    \end{pmatrix},
\end{equation}
where
\begin{equation}\label{eq:caustic-fund-forms}
   \omega_c =
     \sqrt{\hat\omega\hat\theta}\,dz-\sqrt{-1}\frac{d\rho}{4\rho},
     \quad
   \theta_c =
     \sqrt{\hat\omega\hat\theta}\,dz+\sqrt{-1}\frac{d\rho}{4\rho}
   \qquad \bigl(\omega=\hat\omega\,dz,\theta=\hat\theta\,dz\bigr) .    
\end{equation}

Since $U$ contains no umbilic points of $f$,
both $\omega$ and $\theta$ have no zeroes.
Thus $\omega_c$ and $\theta_c$ have no common zero, which implies that: 
\begin{theorem}[Roitman~\cite{R}]\label{lem:caustic-front}
 The caustic $C_f=E^{}_cE^*_c\colon{}U\to H^3$ of a flat front 
 $f\colon{}U\to H^3$ without umbilic points is a flat front 
 with canonical forms $\omega_c$ and $\theta_c$ as in
 \eqref{eq:caustic-fund-forms}.
 Moreover, we have
 \[
    E_c = \frac{(-1)^{1/4}\alpha^{-1/4}}{\sqrt{2}\sqrt{G-G_*}}
          \begin{pmatrix}
	   G + \sqrt{\alpha} G_* & 
	   \sqrt{-1}(G-\sqrt{\alpha} G_{*}) \\
           1 + \sqrt{\alpha} 
	   & \sqrt{-1}(1-\sqrt{\alpha})
	  \end{pmatrix}
    \qquad
    \left(\alpha = \frac{dG\hphantom{_{*}}}{dG_*}\right),
 \]
 where $G$ and $G_*$ are the hyperbolic Gauss maps of $f$.
 In particular, the hyperbolic Gauss maps $(G_c,G_{c,*})$ 
 of $C_f$ are given by
 \[
   G_c = \frac{G+\sqrt{\alpha} G_*}{1+\sqrt{\alpha}},\qquad
   G_{c,*}
     = \frac{G-\sqrt{\alpha} G_*}{1-\sqrt{\alpha}}.
 \]
\end{theorem}
If $z$ is a singular point of $f$, $r_1(z)=0$ holds because
$|\rho(z)|=1$.
Therefore, the caustic of a parallel family $\{f_t\}$ of flat fronts
is the locus of singular points of the fronts $f_t$ for $t\in\R$.

Since the parallel family has a common caustic,
the sets $Z_0$, $Z_c$ and $Z_s$ in Section~\ref{sec:global} 
can be considered as well-defined on the caustic.
In particular, we have the following:
\begin{proposition}\label{prop:sing-caustic}
 Let $f\colon{}U\to H^3$ be a flat front
 without umbilic points, and  with caustic $C_f$,
 where $U\subset \C$ is a simply connected domain.
 Then
 \begin{enumerate}
  \item A point $p\in U$ is a singular point of the caustic $C_f$
	if and only if $p\in Z_c(f)$.
  \item A point $p\in Z_c(f)$ is a non-degenerate singular point of the
	caustic if and only if $S(G)-S(G_*)\neq 0$ holds at $p$, 
	where $G$ and $G_*$ are the hyperbolic Gauss maps of $f$.
  \item A point $p\in Z_c(f)$ where $S(G)-S(G_*)\neq 0$
	is a cuspidal edge of the caustic if and only if $p\not\in Z_s(f)$.
 \end{enumerate}
\end{proposition}
In other words, 
the locus of the cuspidal edges of $\{f_t\}_{t\in\R}$
is the set of regular points of the caustic.
Furthermore, the locus of the swallowtails of $\{f_t\}$ is
the set of cuspidal edges on the caustic, except the points 
at which $S(G)-S(G_*)=0$.
\begin{proof}[Proof of Proposition~\ref{prop:sing-caustic}]
 A point $p\in U$ is a singular point of $C_f$ if and only if
 $|\omega_c|^2=|\theta_c|^2$.
 By \eqref{eq:caustic-fund-forms}, this is equivalent to
 \[
   0=\Im \overline{\sqrt{\hat\omega\hat\theta}}\,
       \frac{\rho'}{\rho}=
   \Im \left[|\hat\omega\hat\theta|
       \frac{1}{\sqrt{\hat\omega\hat\theta}}
       \left(\frac{\hat\theta'}{\hat\theta}-
             \frac{\hat\omega'}{\hat\omega}\right)\right]
   =|\hat\omega\hat\theta|\Im\sqrt{\zeta_c}.
 \]
 Hence the first assertion holds.

 In this case, $p$ is a degenerate singular point of $C_f$ if and only
 if
 \[
    0=\hat\theta_c'\hat\omega_c-\hat\omega_c'\hat\theta_c
     =\sqrt{\hat\omega\hat\theta}\bigl(s(\hat\theta)-s(\hat\omega)\bigr).
 \]
 Then by \eqref{eq:schwarz-Gauss}, we have the second assertion.

 Finally, if $p$ is a non-degenerate singular point of $C_f$, 
 $p$ is a cuspidal edge if and only if 
 \begin{equation}\label{eq:caustic-cusp}
   \Im\frac{1}{\sqrt{\hat\omega_c\hat\theta_c}}
   \left(
    \frac{\hat\theta_c'}{\hat\theta_c} -
    \frac{\hat\omega_c'}{\hat\omega_c}
   \right)\neq 0.
 \end{equation}
 Here, by direct calculation,
 \begin{align*}
    \frac{1}{\sqrt{\hat\omega_c\hat\theta_c}}
    \left(
    \frac{\hat\theta_c'}{\hat\theta_c} -
    \frac{\hat\omega_c'}{\hat\omega_c} 
    \right) &= 
    \frac{1}{\sqrt{\hat\omega_c\hat\theta_c}^3}
    \frac{\sqrt{-1}}{4}\sqrt{\hat\omega\hat\theta}^3
    \bigl(s(\hat\theta)-s(\hat\omega)\bigr)\\
    &=\sqrt{-1}\frac{\zeta_s(z)}{16+\zeta_c(z)}
 \end{align*}
 holds.
 Since $\zeta_c(z)$ is a positive real number if $z\in Z_c$, 
 \eqref{eq:caustic-cusp} holds if and only if
 $\zeta_s(z)\not\in\sqrt{-1}\R$.
\end{proof}
\section{Examples}
\label{sec:example}

\begin{figure}
\begin{center}
\begin{tabular}{c@{\hspace{2em}}c@{\hspace{2em}}c@{\hspace{2em}}c}
 \includegraphics[width=2.5cm]{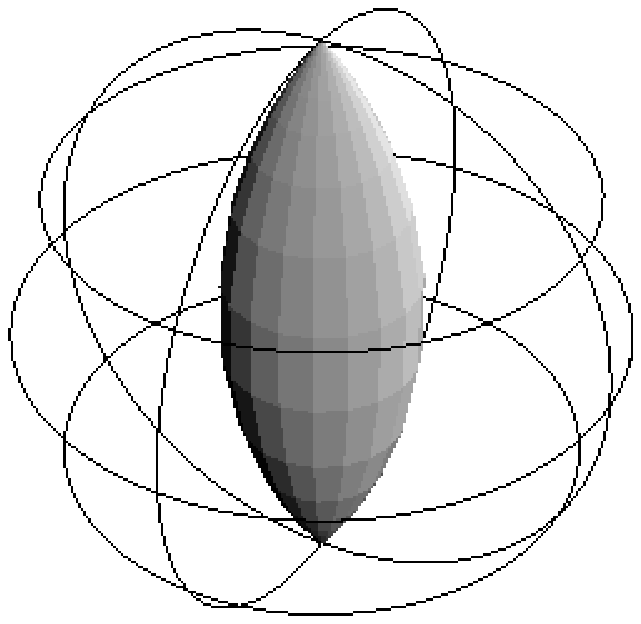} &
 \includegraphics[width=2.5cm]{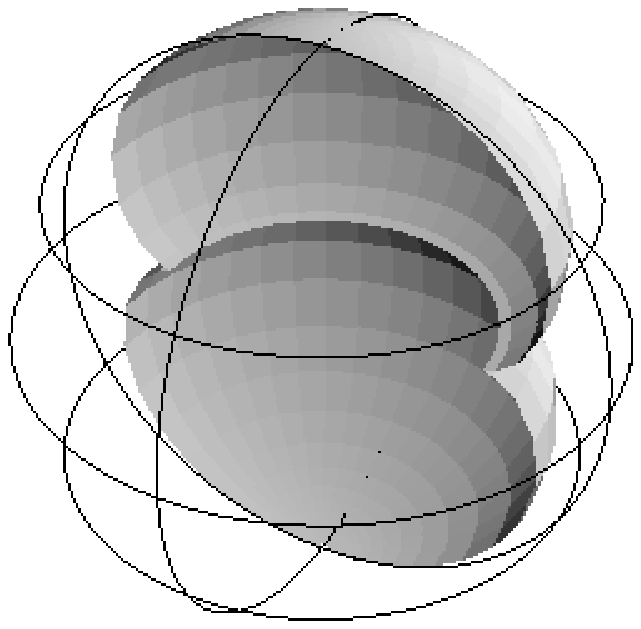} &
 \includegraphics[width=2.5cm]{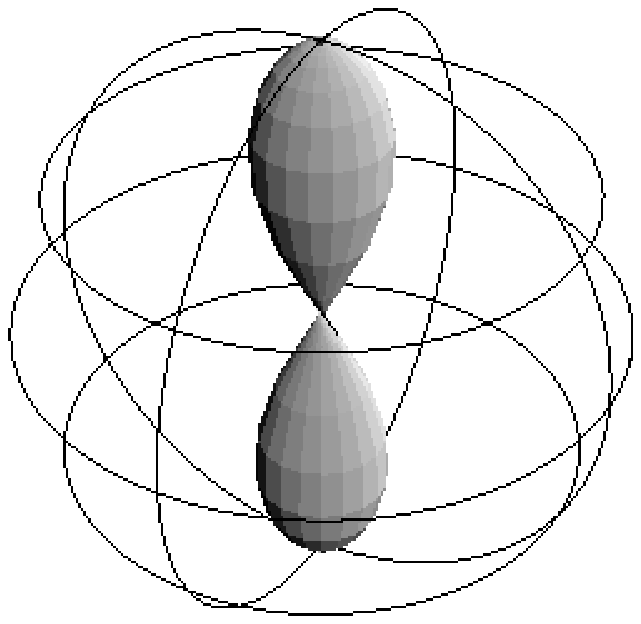} &
 \includegraphics[width=2.5cm]{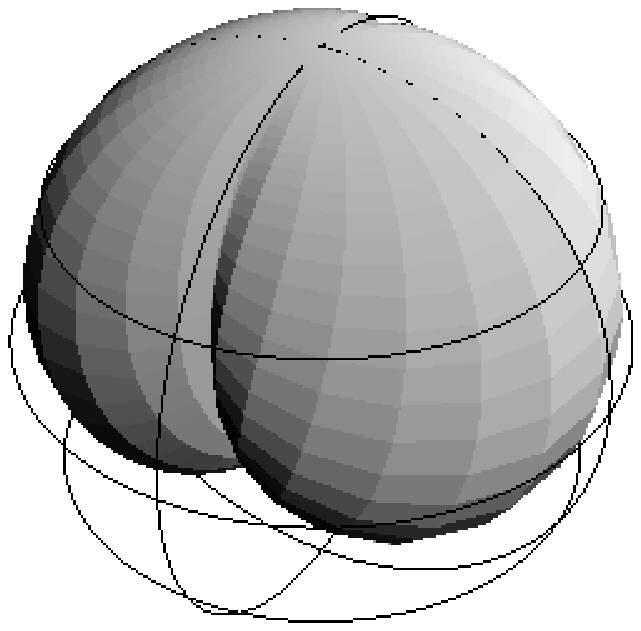}  \\
 \small (a) cylinder &
 \small (b) snowman &
 \small (c) hourglass &
 \small (d) peach front
\end{tabular}
\end{center}
 \caption{Examples~\ref{ex:cylinder}, \ref{ex:revolution} and
          \ref{ex:peach}}\label{fig:revolution}
\end{figure}

Here we give examples that reaffirm the properties of singularities
in Theorems~\ref{thm:criterion} and \ref{thm:global}.
We make examples of flat fronts by choosing hyperbolic Gauss maps
$G$ and $G_*$ as follows:
Let $G$ and $G_*$ be meromorphic functions on a compact Riemann
surface $\overline M^2$ such that $G$ is not identically equal to $G_*$,
and let
\[
    \{p_1,\dots,p_N\}=\{p\in\overline M^2\,;\,G(p)=G_*(p)\}
    \qquad
    \text{and}\qquad
    M^2=\overline M^2\setminus\{p_1,\dots,p_N\}.
\]
If the period condition
\begin{equation*}
 \oint_{\gamma}\frac{dG}{G-G_*}\in\sqrt{-1}\R
\end{equation*}
holds for any loop $\gamma$ on $M^2$, we have the parallel family of 
a complete flat front 
\[
    f_t\colon{}M^2 =\overline M^2\setminus\{p_1,\dots,p_N\} \longrightarrow H^3
\]
by substituting $G$ and $G_*$ in the representation formula
\eqref{eq:G-repr}
in Remark~\ref{rem:DescriptionG}
with $c=e^{-t/2}$.
Moreover, by \eqref{eq:Q-Gauss}, we have
 \[
    \rho = -\frac{\varDelta^4}{(G-G_*)^2}\frac{dG_*}{dG}.
 \]
For details, see \cite{KUY2}.
\begin{example}[Cylinders]\label{ex:cylinder}
 Let $G=z$ and $G_*=1/z$ on $\overline M^2 =\C \cup \{\infty\}$.
 Then $\varDelta=e^{-t/2}\sqrt{z^2-1}$ and 
 $\rho(z)=e^{-2t}$.
 So $|\rho|=1$ if and only if $t=0$ and then all points of the front
 are singular.
 When $t=0$, $\zeta_c=\zeta_s=0$ identically,
 and the surface degenerates to a single geodesic line.
 When $t\neq 0$, we have a cylinder with no singularities
 (Figure~\ref{fig:revolution} (a)).
\end{example}
\begin{example} [Flat fronts of revolution]\label{ex:revolution}
 Let $G=z$ and $G_*=\mu z$ on $\overline M^2=\C \cup \{\infty\}$, where
 $\mu\in\R\setminus\{1\}$.
 Then $\varDelta = e^{-t/2}z^{1/(1-\mu)}$.
 The set of singular points is
 \[ 
     \Sigma_t=\left\{e^{\sqrt{-1}\beta}
               \left(\frac{e^{t}|1-\mu|}{\sqrt{|\mu|}}
                \right)^{\frac{1-\mu}{\mu+1}}\,;\,\beta\in\R\right\}.
 \]
 Since $\sqrt{\zeta_c}|_{\Sigma_t}=\pm 2\sqrt{-1}(\mu+1)/\sqrt{\mu}$ 
 is constant, and real if and only if $\mu<0$, the 
 singularities are cuspidal edges when $\mu>0$.
 When $\mu<0$, $\zeta_s|_{\Sigma_t}=0$.
 In this case, the singular points are neither cuspidal edges nor 
 swallowtails,
 although they are nondegenerate.
 (The singular image is a single point, 
 but there are many singular points in the domain.)

 When $\mu>0$, the image of $\Sigma_t$ is a circular cuspidal edge
 centered about the surface's rotation axis (the {\em snowman}, see
 Figure~\ref{fig:revolution} (b)).
 When $\mu<0$, the image of $\Sigma_t$ is a single point on the rotation
 axis (the {\em hourglass}, Figure~\ref{fig:revolution} (c)).

 When $\mu=0$, the surface is a horosphere, and when $\mu=-1$, the
 surface is a hyperbolic cylinder.
 
 When $\mu\to +1$, the entire surface approaches the ideal boundary
 $\partial H^3$ of $H^3$. 
 When $\mu0$, the corresponding caustic is a cylinder.
\end{example}
\begin{example}[Peach fronts]\label{ex:peach}
 Let $G=z+\frac{1}{2}$ and 
 $G_*=z-\frac{1}{2}$ on $\overline M^2 =\C \cup \{\infty\}$. 
 Then one has a parallel family of flat fronts
 $f_t\colon{}\C\to H^3$ resembling peaches.
 Since $\varDelta=e^{-t/2}e^z$, the set of singular points
 is $\Sigma_t=\{t/2+\sqrt{-1}y\,;\,y\in\R\}$.
 Since $\sqrt{\zeta_c}|_{\Sigma_t}=\pm 4\sqrt{-1}$ is not real, we have
 a single cuspidal edge along a vertical line on $\C$.
 This cuspidal edge travels out to the end, hence we have a simple
 example for which every open neighborhood of the end contains singular
 points, in particular, each $f_t$ is not complete.
 As noted in \cite{R}, the corresponding caustic is the horosphere.
 See Figure~\ref{fig:revolution} (d).
\end{example}
\begin{figure}
 \begin{center}
  \begin{tabular}{c@{\hspace{6em}}c}
   \includegraphics[width=2.5cm]{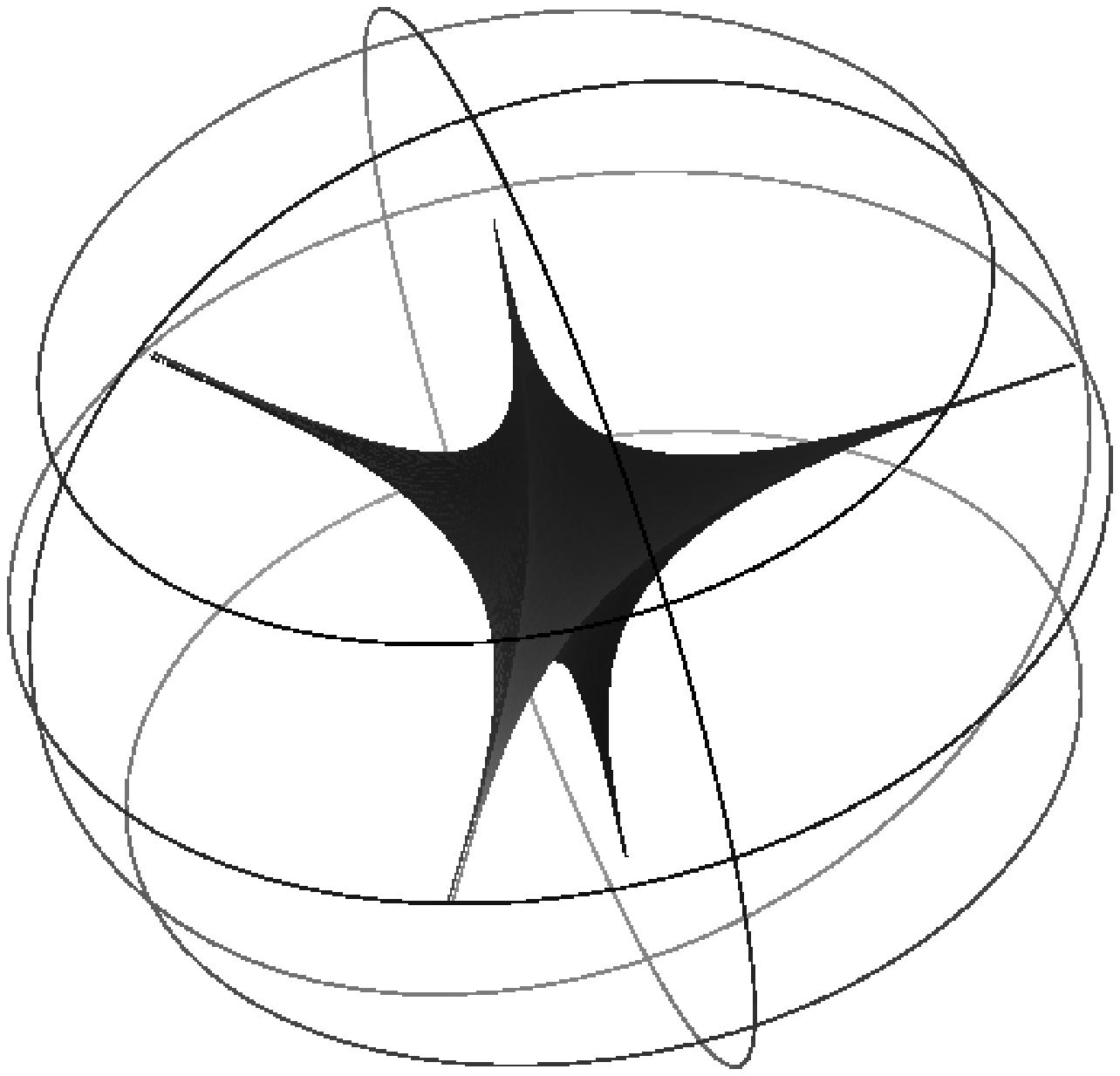} &
   \includegraphics[width=2.5cm]{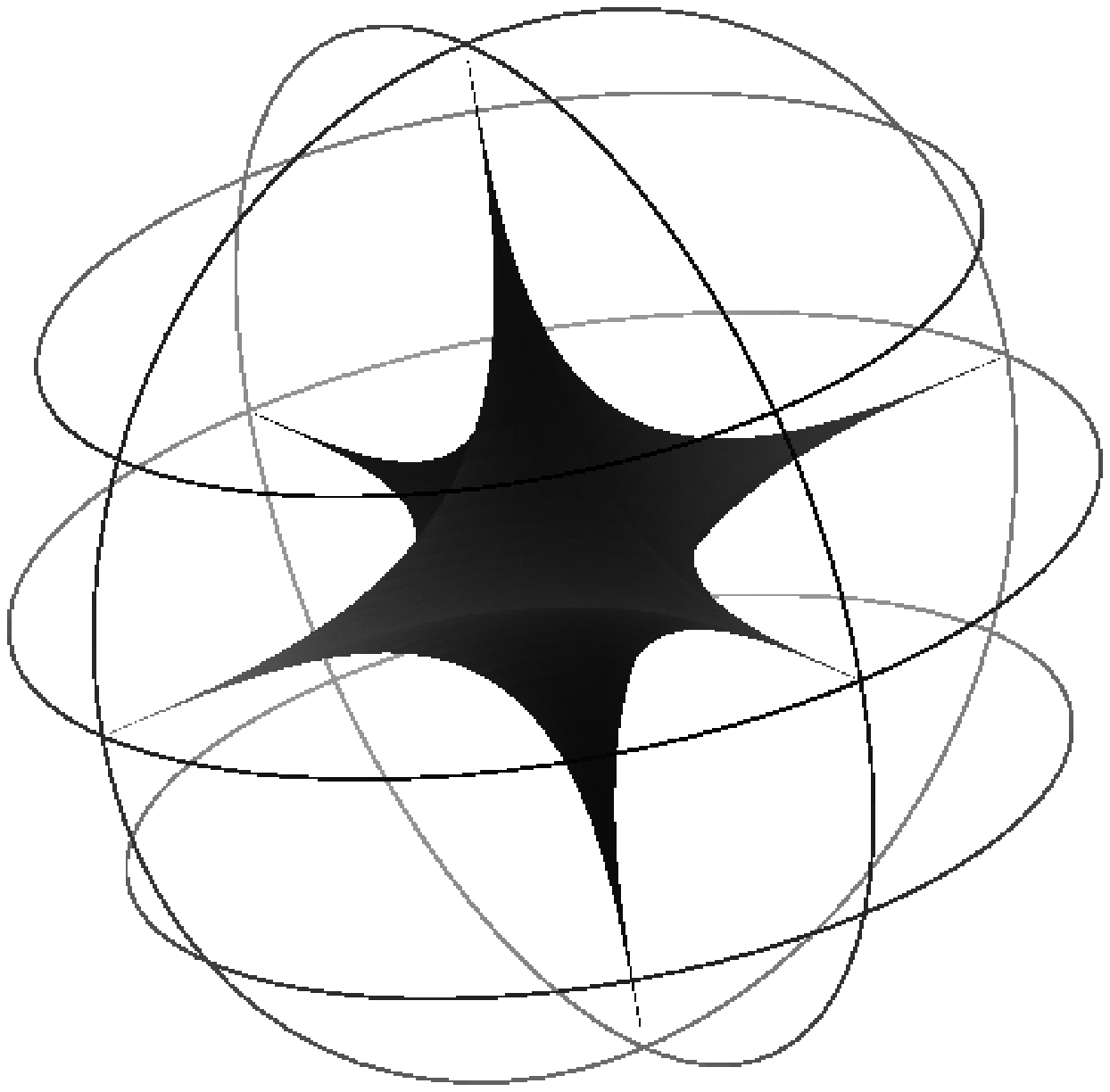}\\
   \small The caustic of a $3$-noid &
   \small The caustic of a $4$-noid 
  \end{tabular}
 \end{center}
 \caption{Example~\ref{ex:nnoid}: caustics for $3$-noids and
  $4$-noids.}
 \label{fig:nnoid}
\end{figure}
\begin{example}[$n$-noid flat fronts]\label{ex:nnoid}
 $n$-ended flat fronts for $n\geq 3$ can be made with 
 $G=z$, $G_*=z^{1-n}$ on $M^2=\C\cup\{\infty\}\setminus\{z\,;\,z^n=1\}$.
 For all $t$, the points $z=0$, $\infty$ are finite and non-singular.
 We have $\varDelta=e^{-t/2}\sqrt[n]{z^n-1}$,
 \begin{align*}
  \rho(z) &= (n-1)e^{-2t}z^{n-2}(z^n-1)^{\frac{4-2n}{n}},\\
  \sqrt{\zeta_c} &= \pm \frac{(n-2)(z^n+1)}{\sqrt{n-1}\,
  z^{\frac{n}{2}}},\\
  \zeta_s &= \frac{n(2-n)}{2(n-1)}\frac{(z^n-1)^2}{z^n}.
 \end{align*}%
 Then $\sqrt{\zeta_c}$ is real when $z\in S^1$ or
 $z^n\in\R^+\setminus\{0\}$.
 Since $\frac{1}{n-2}\zeta_c+\frac{2}{n}\zeta_s=\frac{4(n-2)}{n-1}$,
 $Z_c\cap Z_s=\{z\,;\,z^n=1\}\cap M^2=\emptyset$, 
 and one can easily show that  $Z_0=\{z\,;\,z^n=-1\}$.
 In the case $n=3$ (resp.\ $4$), 
 if $6t\log 2$ (resp.\ $2t \log(3/2)$), there are twelve
 (resp.\ sixteen) swallowtails, and all other singularities are
 cuspidal.
 If $6t=\log 2$ (resp.\ $2t=\log(3/2)$), there are three 
 (resp.\ four)
 degenerate singularities at the points $z^3=-1$ (resp.\ $z^4=-1$).
 If $6t<\log 2$ (resp.\ $2t<\log(3/2)$), then there are six
 (resp.\ eight)
 swallowtails, and otherwise cuspidal edges.
 For figures of $3$-noids, see \cite{KUY1}.
 The caustics corresponding to $3$-noids and $4$-noids are shown in
 Figure~\ref{fig:nnoid}.
\end{example}
\begin{figure}
\begin{center}
\begin{tabular}{c@{\hspace{3em}}c@{\hspace{3em}}c}
 \includegraphics[width=2.5cm]{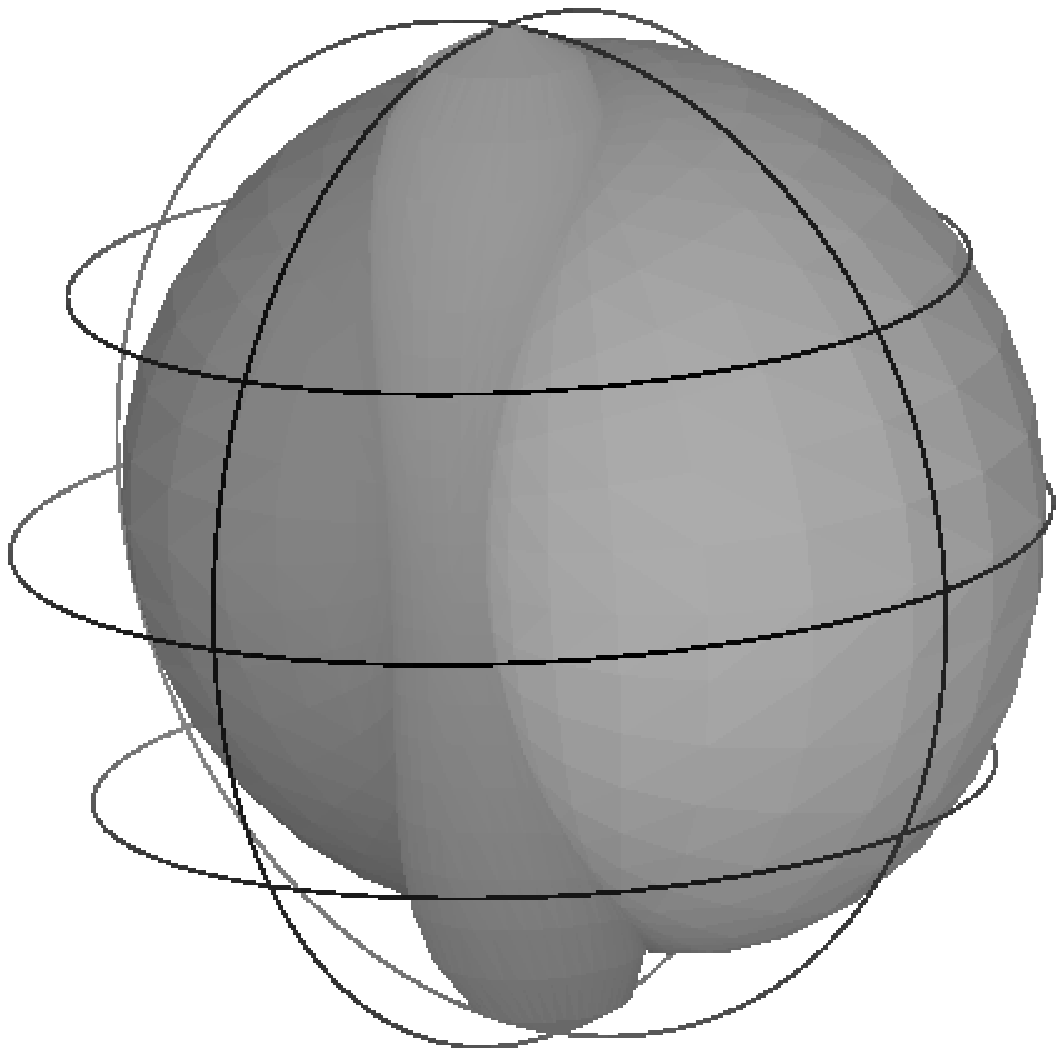} &
 \includegraphics[width=2.5cm]{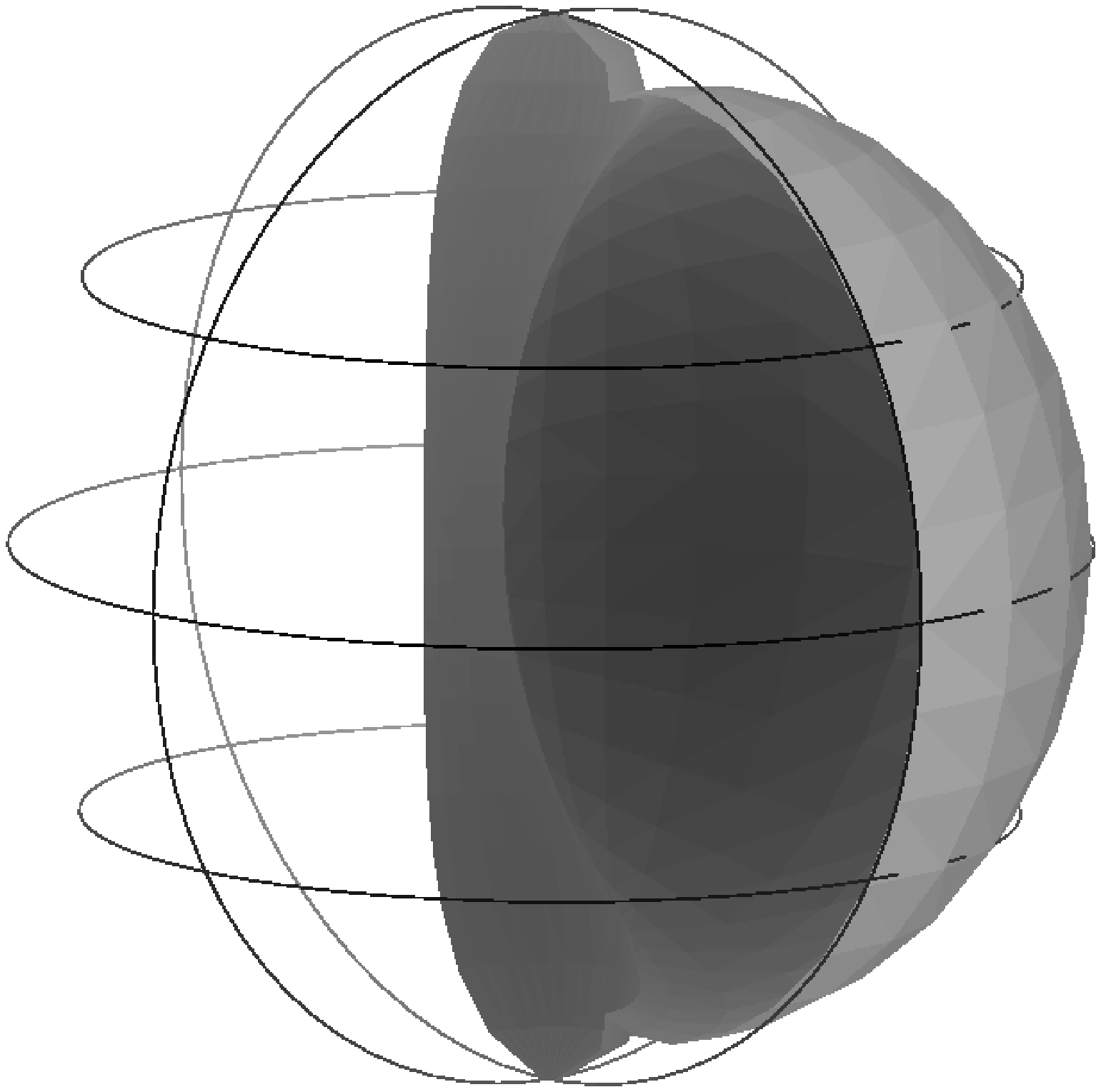} &
 \includegraphics[width=2.5cm]{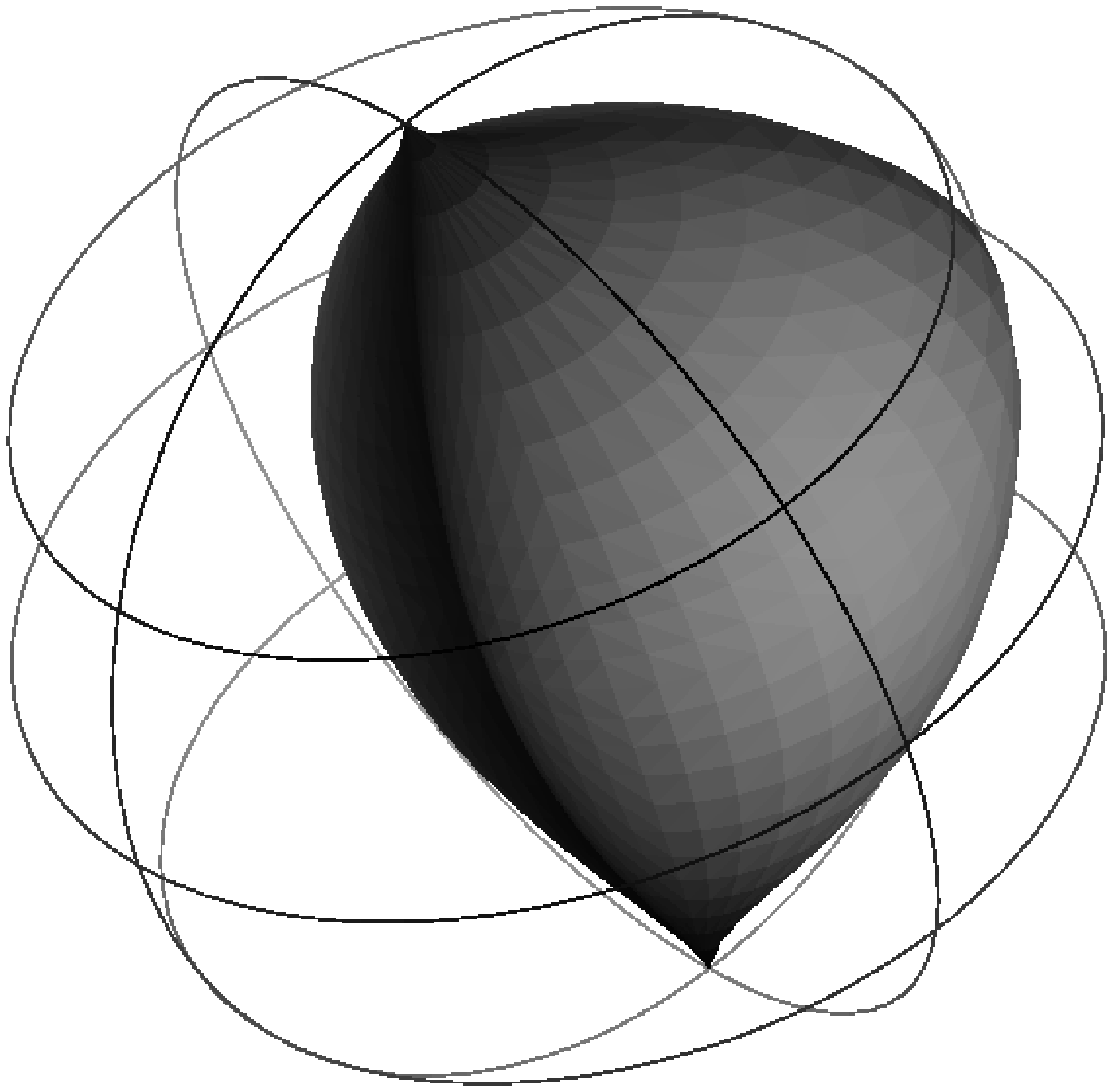} \\
 \small $(G,G_*)=(z,z^2)$ &
 \small $(G,G_*)=(z,z^2)$ (half cut) &
 \small $(G,G_*)=(z,z^2)$ (caustic) \\[1ex]
 \includegraphics[width=2.5cm]{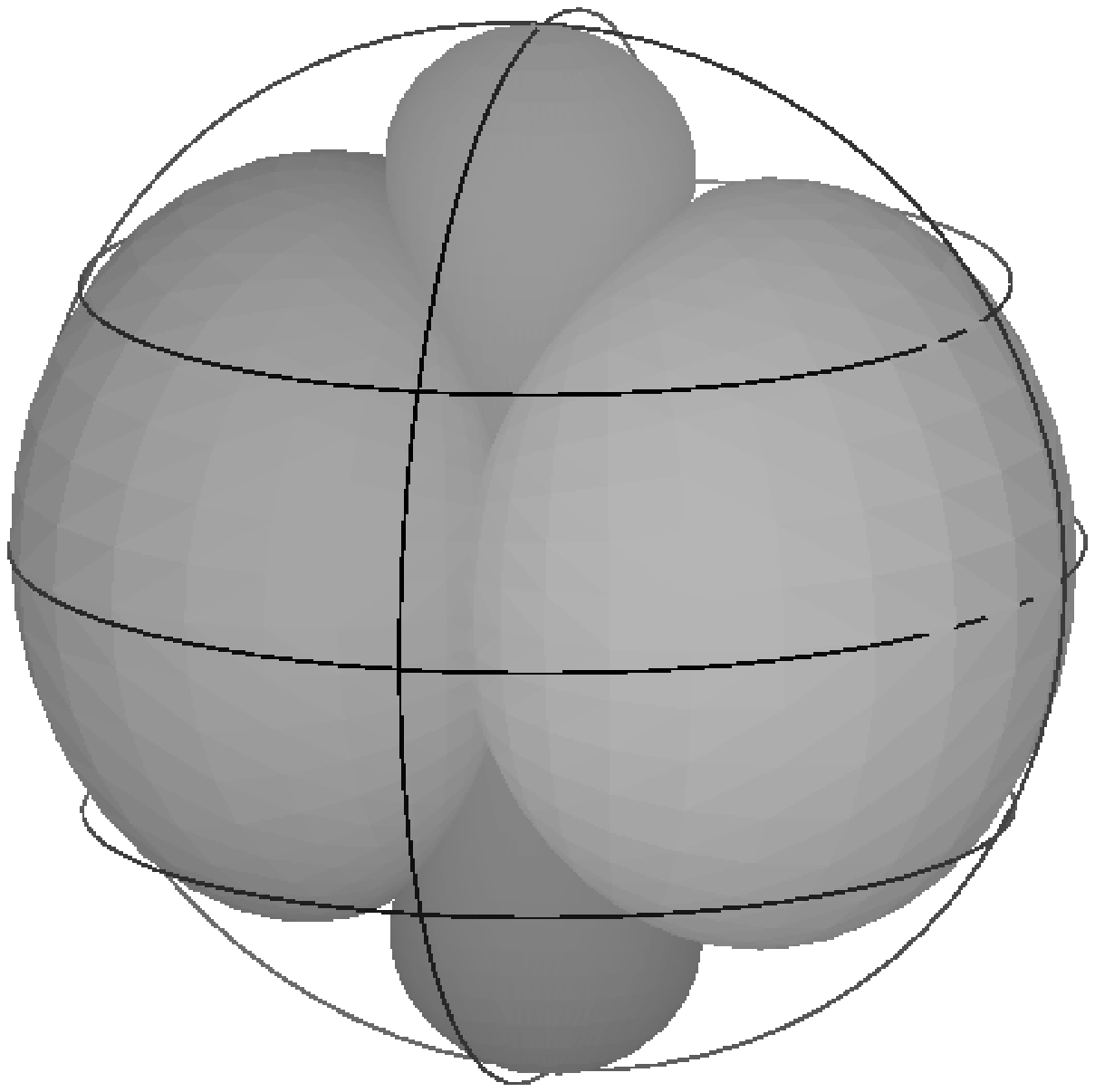} &
 \includegraphics[width=2.5cm]{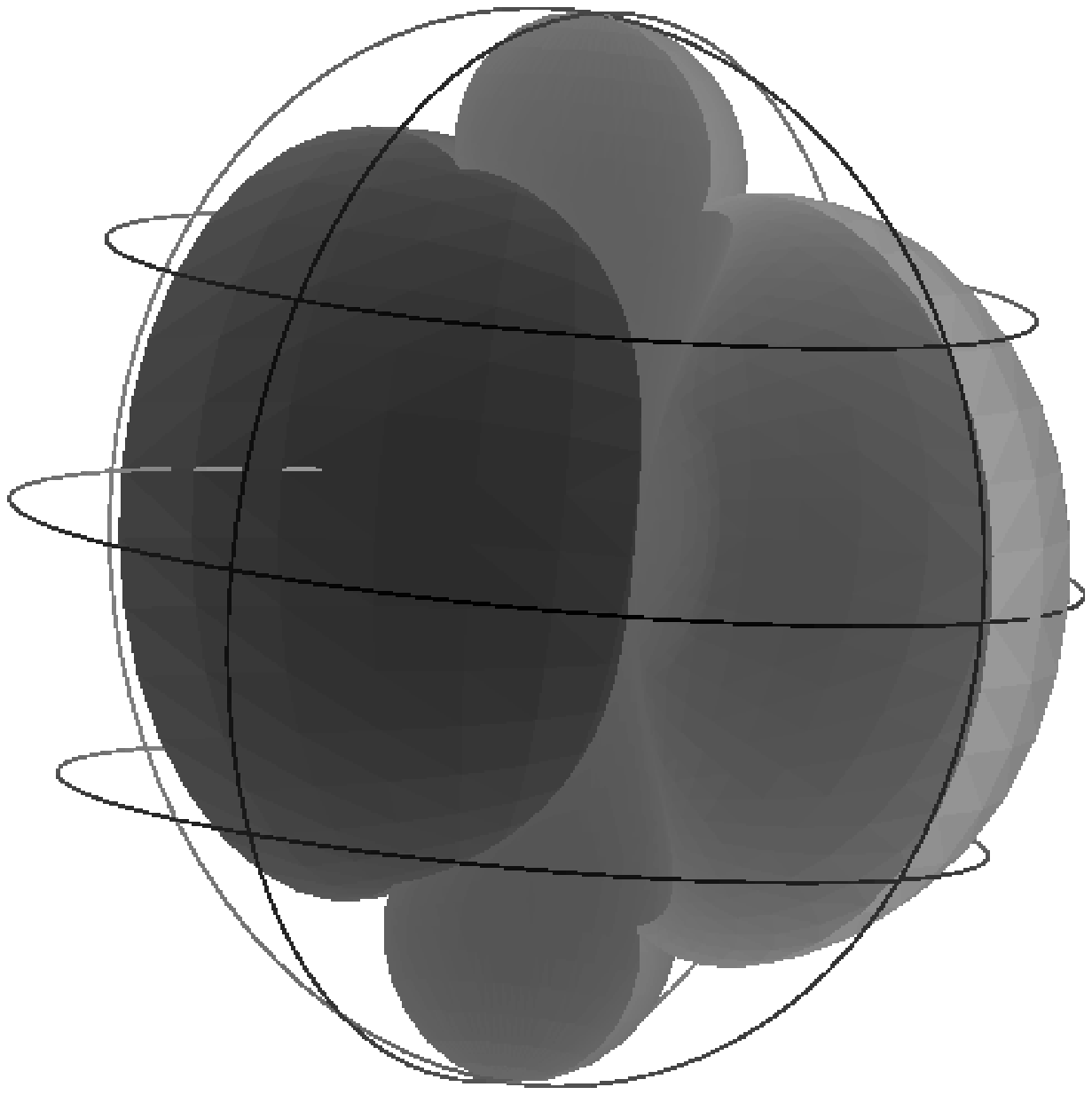} &
 \includegraphics[width=2.5cm]{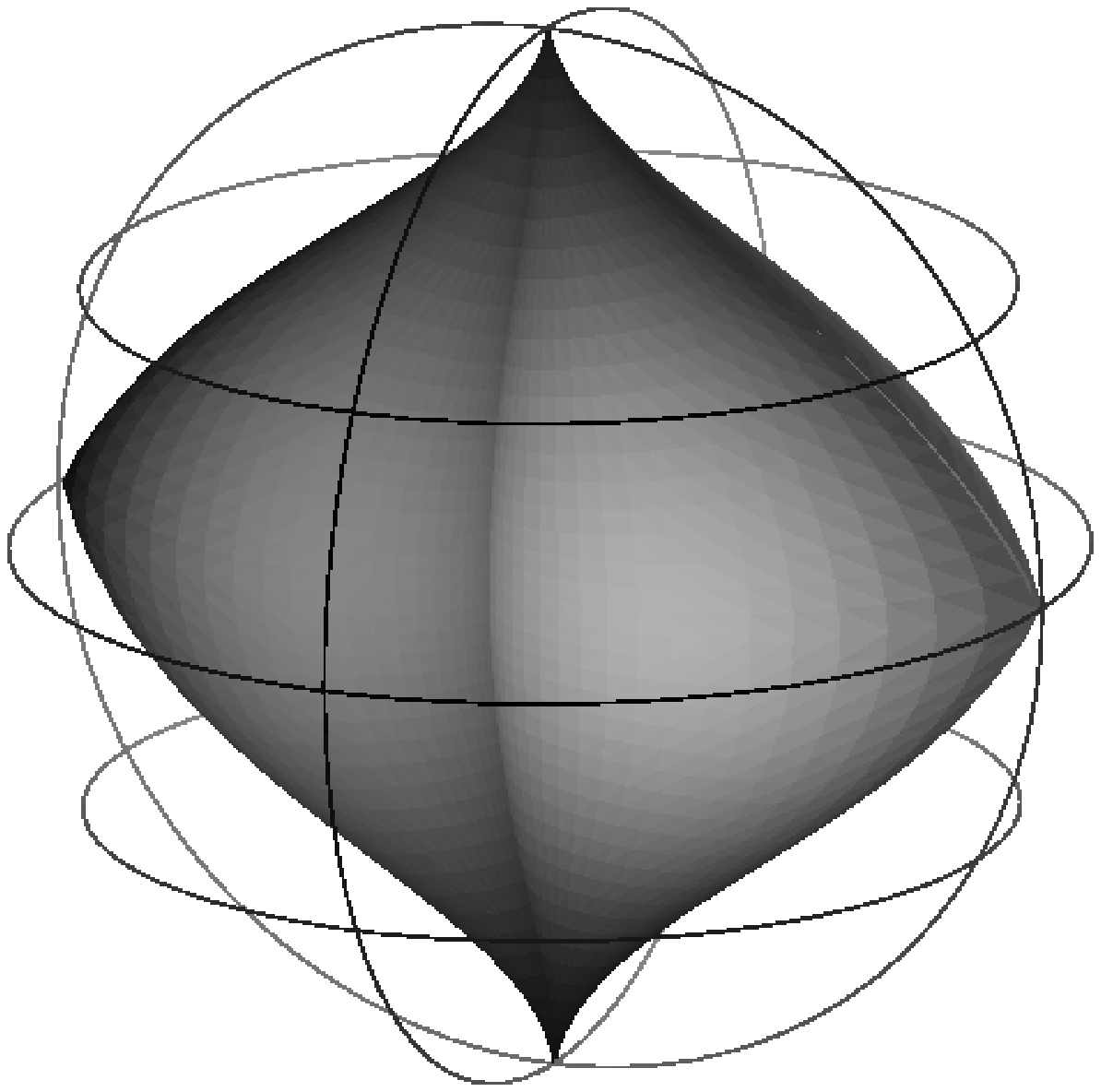} \\
 \small $(G,G_*)=(z,z^3)$ &
 \small $(G,G_*)=(z,z^3)$ (half cut) &
 \small $(G,G_*)=(z,z^3)$ (caustic) \\
\end{tabular}
\end{center}
 \caption{Example~\ref{ex:mn}}
 \label{fig:12}
\end{figure}
\begin{figure}
\begin{center}
\begin{tabular}{c@{\hspace{3em}}c@{\hspace{3em}}c}
 \includegraphics[width=3cm]{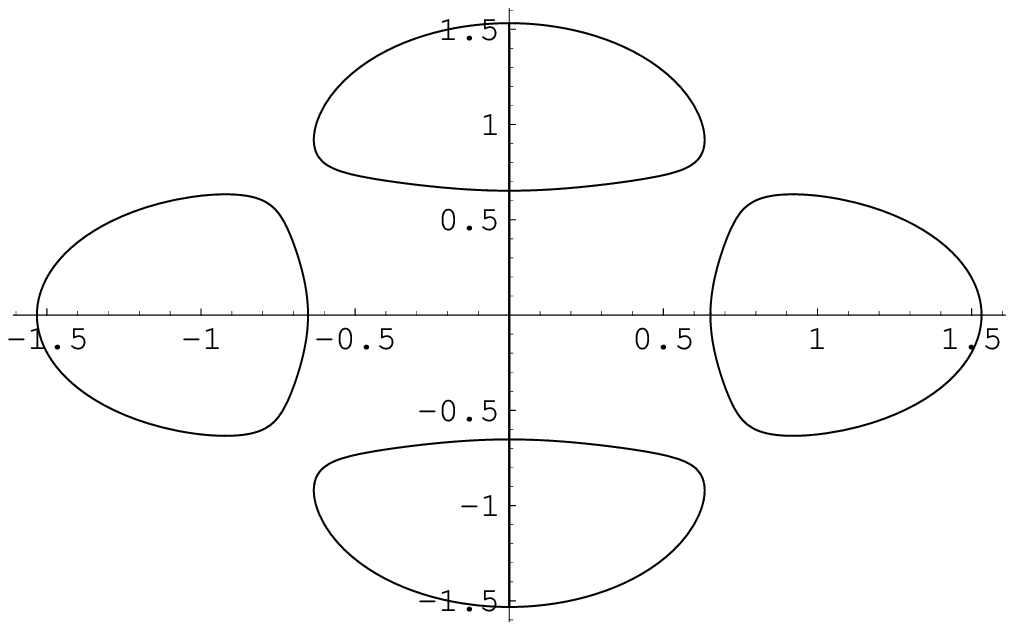} &
 \includegraphics[width=3cm]{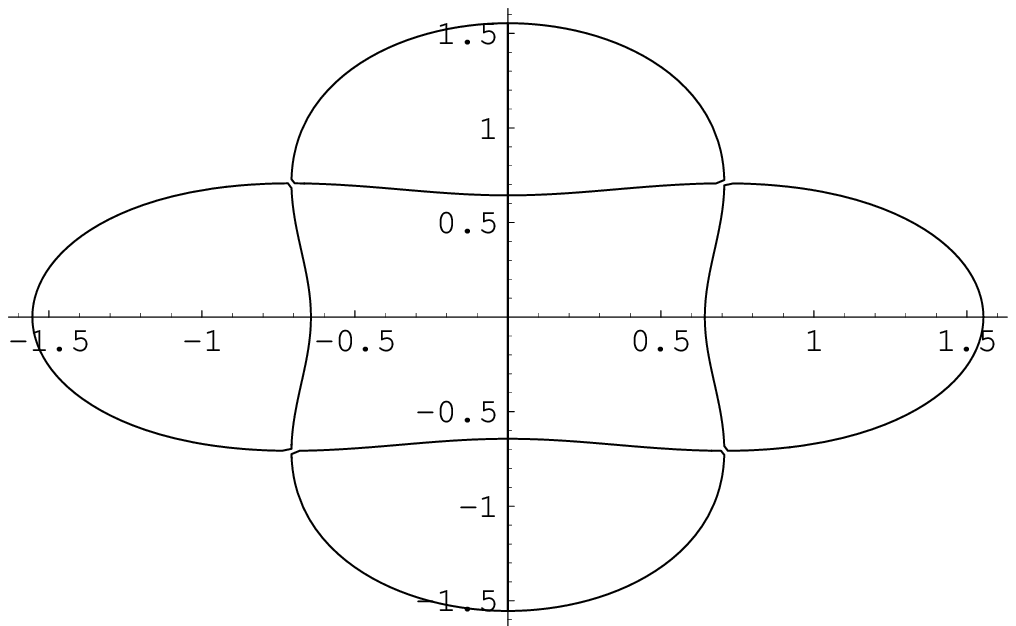} &
 \includegraphics[width=3cm]{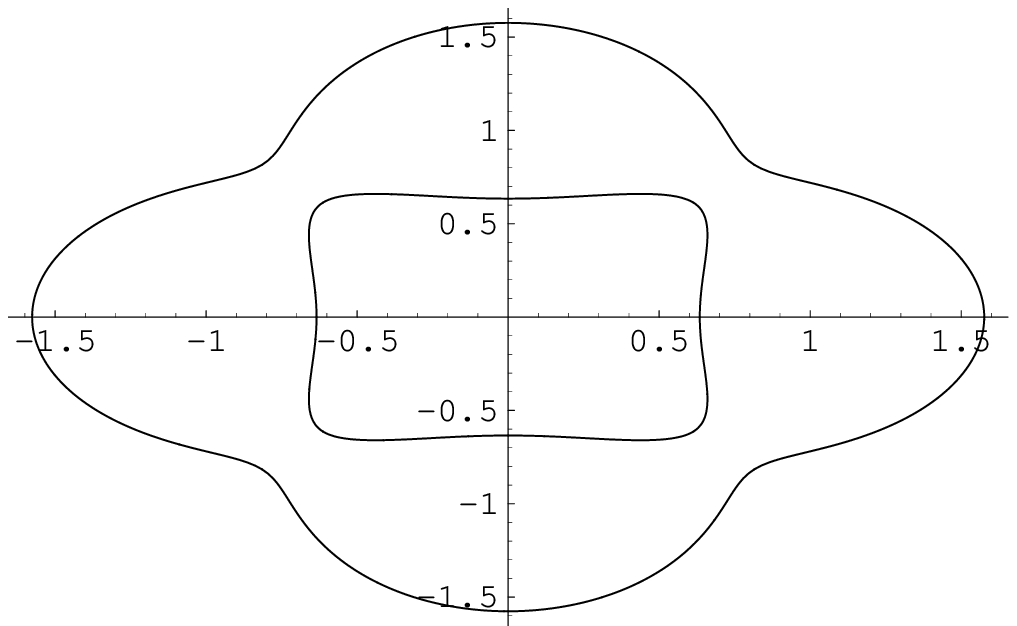} \\
 \includegraphics[width=3cm]{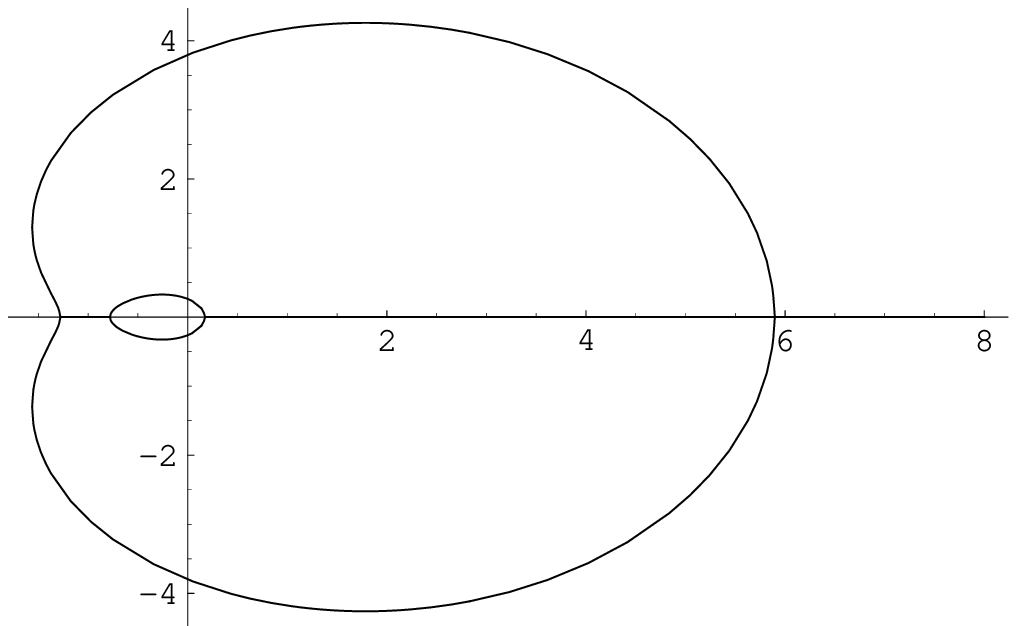} &
 \includegraphics[width=3cm]{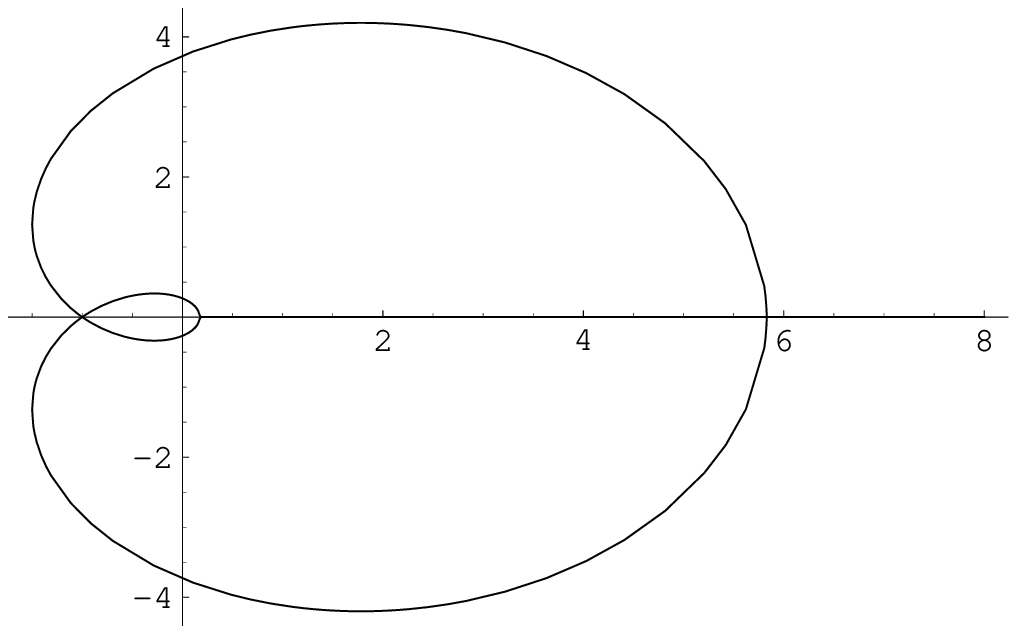} &
 \includegraphics[width=3cm]{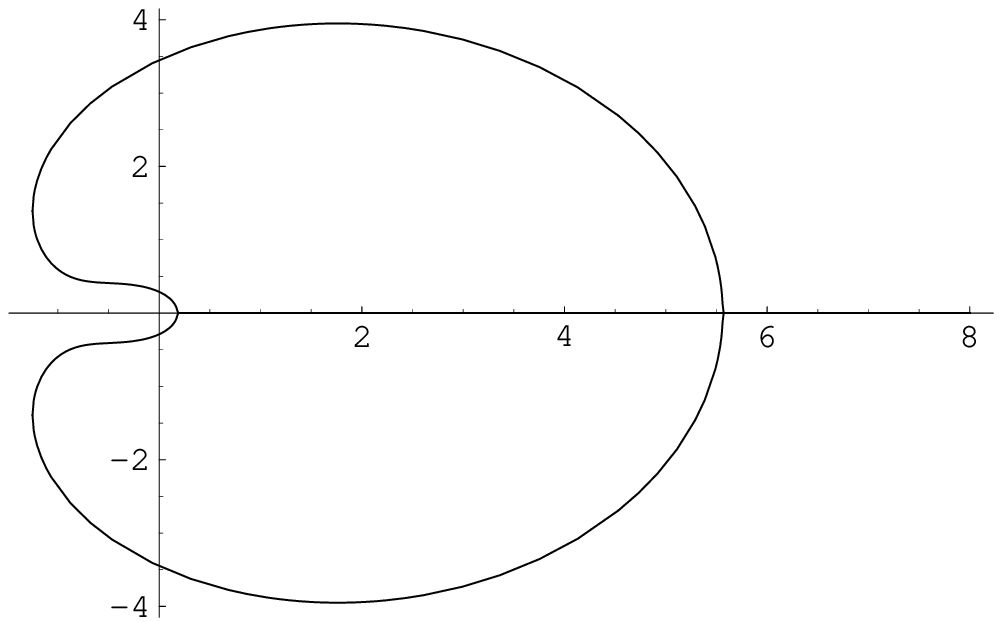} \\
\end{tabular}
\\
\begin{quote}
{\small 
  The upper row shows (from left to right) the three 
  cases $2t\log(3/2)$, $2t=\log(3/2)$ and $2t<\log(3/2)$
  for $n=4$ in Example~\ref{ex:nnoid}.
  The lower row shows (from left to right) the three cases
  $e^{2t}<1/32$, $e^{2t}=1/32$ and $e^{2t}1/32$ for $n=1$ 
  and $m=2$ in Example~\ref{ex:mn}.
}
\end{quote}
\end{center}
 \caption{
  Singular curves in the domains $M^2$ in 
  Examples~\ref{ex:nnoid} and \ref{ex:mn}.}
 \label{fig:singular-locus}
\end{figure}
\begin{example}[Flat fronts with $G=z^n$ and $G_*=z^m$ for %
$1\leq n < m$]%
\label{ex:mn}
 In this case, $M^2=\C\setminus\{z;z^m=z^n\}$,
 and we have $m-n+2$ ends.
 Then 
 \begin{align*}
   \varDelta & =e^{-t/2}z^n(1-z^{m-n})^{\frac{n}{n-m}},\\
   \rho(z) &=  -\frac{m}{n}e^{-2t}z^{m+n}
          (1-z^{m-n})^{\frac{2(m+n)}{n-m}},\\
   \sqrt{\zeta_c}&=
       \pm\sqrt{-1}\frac{(m+n)(z^m+z^n)}{\sqrt{mn}}
       z^{-\frac{m+n}{2}},\\
   \zeta_s&=
      \frac{m^2-n^2}{2mn} z^{-m-n}(z^m-z^n)^2.
 \end{align*}
 For small values of $m$ and $n$, we can easily investigate the
 singularities.
 
 For $n=1$ and $m=2$ (resp.~$m=3$), 
 for all $t$, all singularities are always cuspidal edges except 
 two (resp.\ four) swallowtails when $e^{2t}<1/32$
 (resp.\ $e^{2t}<3/16$) and at one (resp.\ two) degenerate
 singularity (resp.\ singularities)
 when $e^{2t}=1/32$ (resp.\ $e^{2t}=3/16$).
 As the value $e^{2t}$ increases through $1/32$ (resp.\ $3/16$),
 the two (resp.\ four) swallowtails come together into a single
 (resp.\ two) degenerate singularity (resp.\ singularities)
 and then disappear, leaving only cuspidal
 edges.
 Surfaces for $n=1$ and $m=2$ and $3$, and their corresponding caustics,
 are shown in Figure~\ref{fig:12}.

 For $n=2$ and $m=3$, and for all $t$, the singular points are
 always cuspidal edges or swallowtails.
 All singular points are cuspidal edges, except for one swallowtail when
 $t<0$.
 When $t\geq 0$, there are no swallowtails.
 As $t$ increases to $0$, the swallowtail moves out to an end and
 disappears when $t=0$.
 For $t=0$, the singular set is the line $\Re z=1/2$, and hence 
 the cuspidal edge travels out to the end $z=\infty$.%
\end{example}

\appendix
\section{Proof of Lemma \ref{fact:Z}}\label{app}
In this appendix we prove Lemma \ref{fact:Z},
as only a sketch of the proof  given in \cite{Z}.
The authors hope this will help readers 
who are not familiar with singularity theory.

We use the following three well known facts:

\begin{fact}\label{fact:A1}
 Let $f:M^n\to \R^N$ be an immersion of an $n$-manifold $M$.
 Then for each point $p\in M$, there exists a neighborhood
 $U$ of $p$ such that the restriction $f|_U$ is an embedding.
\end{fact}

\begin{fact}\label{fact:A2}
 Let $U_1$ and $U_2$ be neighborhoods of the origin $o$ in $\R^n$.
 Let $f_i:U_i\to \R^N$ $(i=1,2)$ be two embeddings such that 
 $f_1(U_1)\subset f_2(U_2)$ and $f_1(o)=f_2(o)$.
 Then there exists a local diffeomorphism
 $\phi:U_1\to U_2$ such that $f_1=f_2\circ \phi$ holds.
\end{fact}

It is well-known that a front can be considered 
as a projection of a Legendre immersion 
$L:U\to P(T^*\R^3)$, where $U$ is a domain in $\R^2$
and $P(T^*\R^3)$ is the projective cotangent bundle.
The canonical contact structure of the unit cotangent 
bundle $T^*_1\R^3$ is the pull-back of that of $P(T^*\R^3)$.
We remark that this contact structure on $P(T^*\R^3)$
does not depend on the Riemannian metric
on $\R^3$ (see \cite{A}).
So we have the following:

\begin{fact}\label{fact:A3}
 Let $f:U\to \R^3$ be a front, 
 where $U$ is an open subset of $\R^2$ and
 \[
    \Phi:\R^3\longrightarrow \R^3
 \]
 is a diffeomorphism. 
 Then the composition $\Phi\circ f$ is also a front.
\end{fact}

By the above three facts, the theorem reduces into the
following proposition:

\begin{proposition}\label{prop:A4}
 Let $f_i:U_i\to \R^3$ $(i=1,2)$ 
 be two fronts satisfying $f_1(o)=f_2(o)$,
 whose associated Legendrian immersions
 $L_{f_i}:U_i\to T^*_1\R^3$ are
 embeddings, where $U_i$ are 
 neighborhoods of the origin $o$ in $\R^2$. 
 Suppose that there exists a relatively compact neighborhood
 $V_i$ of $o$ $(i=1,2)$ such that
 \begin{enumerate}
  \item\label{item:propA4-1} 
       The closure $\overline{V_i}$ is contained in $U_i$ for $i=1,2$.
  \item\label{item:propA4-2} 
       The set of regular points of $f_i$ in $\overline{V_i}$
       is dense in $\overline{V_i}$ $(i=1,2)$.
  \item\label{item:propA4-3} 
       $f_1(\overline{V_1})= f_2(\overline{V_2})$
 \end{enumerate}
Then $L_{f_1}(\overline{V_1})=L_{f_2}(\overline{V_2})$ holds.
\end{proposition}

Before proving this proposition, 
we give the proof of Lemma \ref{fact:Z}:

\begin{proof}[Proof of Lemma \ref{fact:Z}]
\ref{item:fact-Z-1} follows from \ref{item:fact-Z-2} immediately.
So it is sufficient to show \ref{item:fact-Z-1} implies \ref{item:fact-Z-2}.
By Fact~\ref{fact:A3},
we may assume $f_1(V_1)=f_2(V_2)$.
Without loss of generality we may assume that
$V_1$ and $V_2$ are relatively compact 
and $\overline{V_1},\overline{V_2}\subset U$.
By Fact A.1, we may assume that
the associated Legendrian immersion
$L_{f_i}:U\to T^*_1\R^3$ is an embedding.
Since $V_1$ and $V_2$ are relatively compact, we have
\[
   f_1(\overline{V_1})=\overline{f_1(V_1)}=\overline{f_2(V_2)}
                =f_2(\overline{V_2}).
\]
Thus by Proposition~\ref{prop:A4}, we have
$L_{f_1}(\overline{V_1})=L_{f_2}(\overline{V_2})$,
in particular we have
\[
    L_{f_1}(V_1)\subset L_{f_2}(U).
\]
By Fact~\ref{fact:A2}, there exists a local diffeomorphism  $\phi$ on
$\R^2$ such that 
$L_{f_2}=L_{f_1}\circ \phi$, which proves the assertion.
\end{proof}

To prove the Proposition~\ref{prop:A4}, we set
\begin{align*}
 S&=f_1(\overline{V_1})=f_2(\overline{V_2}), \\
 Z_i&=\{f_i(p)\in S\,;\, 
          \text{$p\in \overline{V_i}$ is a singular point of $f_i$} \}
                    \qquad (i=1,2), \\
 Z&=Z_1\cup Z_2, \qquad  R=S\setminus Z,
\end{align*}
and first prove the following simple lemma:
\begin{lemma}\label{lem:A5}
For each $a\in S\setminus Z_i$, $f^{-1}_i(a)$  is a finite set.
\end{lemma}

\begin{proof}
Suppose that  $f^{-1}(a)$ is not a finite set.
Without loss of generality, we can take a sequence
$\{p_n\}$ such that
\[
   f_i(p_n)=a \qquad (n=1,2,\dots).
\]
Moreover, by taking a subsequence we may assume $\{p_n\}$ converges to a
point $p\in \overline{V_i}$.
Then by continuity, we have $f_i(p)=a$.
Since $a\in S\setminus Z_i$, $p$ is a regular point of $f_i$.
Thus, there exists a neighborhood $V$ of $p$
such that $f_i|_V$ is an embedding,
which contradicts
\[
 f_i(p_n)=a=f_i(p),
\]
since $p_n\in V$ for sufficiently large $n$.
\end{proof}

\begin{proof}[Proof of Proposition~\ref{prop:A4}]
We fix $a\in R$ arbitrarily.
By the previous lemma, we may set
\[
 f_1^{-1}(a)=\{p_1,\dots,p_m\},\qquad
 f_2^{-1}(a)=\{q_1,\dots,q_l\}.
\]
We identify $T_1^*\R^3$ with $T_1\R^3=\R^3\times S^2$.
Then $L_{f_i}$ ($i=1,2$) is considered as a map into
$\R^3\times S^2$, and
there exist unit vectors $\nu_1,\dots,\nu_m$
and $\xi_1,\dots,\xi_\ell$
such that
\[
 L_{f_1}(p_j)=(a,\nu_j),\quad
 L_{f_2}(q_k)=(a,\xi_k)\qquad (j=1,\dots,m,\ k=1,\dots,l).
\]
Since $L_{f_1}$ and $L_{f_2}$ are embeddings, 
$\nu_1,\dots,\nu_m$ (resp. $\xi_1,\dots,\xi_l$)
are mutually distinct.
Thus the image of $f_i$ at $a$ consists of a mutually transversal
finite number of components of surfaces.
Since $f_1(\overline{V_1})=f_2(\overline{V_2})$,
we can conclude that $m=l$ and
\begin{equation}\tag{$*$}\label{eq:lift-coincide}
 L_{f_1}(p_j)=(a,\nu_j)=(a,\xi_j)=L_{f_2}(q_j) \qquad 
     (j=1,2,\dots,m) 
\end{equation}
for a suitable permutation of $p_1,\dots,p_m$.
Now we set
\[
 W_1=(f_1|_{\overline{V_1}})^{-1}(R),\qquad 
 W_2=(f_2|_{\overline{V_2}})^{-1}(R).
\]
By \eqref{eq:lift-coincide}, we have
\[
 L_{f_1}(W_1)= L_{f_2}(W_2).
\]
Then by the continuity of $L_{f_1}$ and $L_{f_2}$, we have
\[
   L_{f_1}(\overline{W_1})= L_{f_2}(\overline{W_2}).
\]
Thus it is sufficient to show that $W_i$ is dense
in $\overline{V_i}$.
In fact, 
suppose that $(f_i|_{\overline{V_i}})^{-1}(Z)$ has an interior point. 
By the assumption (2) of Proposition A.4,
there exists an open subset $O_i \,\,(\subset V_i)$ such that
$f_i(O_i)\subset Z$ and $f_i$ is an immersion on $O_i$. 
Take a point $q_i\in O_i$.
Let $T_i$ be the tangent plane (as a two dimensional 
affine plane in $\R^3$) of the regular surface
$f_i(O_i)$ at $q_i$, and
\[
   \pi_i:\R^3\to T_i \qquad (i=1,2)
\]
the orthogonal projection.
Since $\pi_i\circ f_i$ has a regular point, 
$\pi_i\circ f_i(O_i)$ contains an interior point.
On the other hand, by Sard's theorem,
the critical value set $\pi_i(Z)$ of $\pi_i\circ f_i$ 
is a measure zero set.
Since $\pi_i\circ f_i(O_i)\subset \pi_i(Z)$, 
this makes a contradiction.
Hence $(f_i|_{\overline{V_i}})^{-1}(Z)$ does not have any interior points.
Since
\[
   (f_i|_{\overline{V_i}})^{-1}(\overline{V_1})
           =(f_i|_{\overline{V_i}})^{-1}(R\cup Z)
           =(f_i|_{\overline{V_i}})^{-1}(R)\cup 
               (f_i|_{\overline{V_i}})^{-1}(Z),
\]
$W_i=(f_i|_{\overline{V_i}})^{-1}(R)$ is dense in $\overline{V_i}$.
\end{proof}

The authors' original proof of Proposition A.4
used the Hausdorff dimension of $f_i(O_i)$. 
Go-o Ishikawa pointed out to us a 
simplification of the proof that requires
only the classical Sard's theorem.


\end{document}